\documentclass[10pt]{article}

\usepackage{amssymb}
\usepackage{latexsym}
\usepackage{euscript}
\usepackage{amsmath,amsthm,amsfonts}
\usepackage{amscd}
\usepackage{epsfig}
\usepackage{xspace}

\newcounter{Universal}
\renewcommand{\theUniversal}{\arabic{Universal}}

\newenvironment{Italic}{\refstepcounter{Universal} \par \vspace{1em}
\noindent {\bf (\theUniversal)}\ \it}{\par}

\newenvironment{Theorem}{\begin{Italic}{\sc Theorem: }}{\end{Italic}}
\newenvironment{Proposition}{\begin{Italic}{\sc Proposition:} }{\end{Italic}}
\newenvironment{Definition}{\begin{Italic}{\sc Definition: } }{\end{Italic}}

\newenvironment{Lemma}{\begin{Italic}{\sc Lemma: }}{\end{Italic}}
\newenvironment{Conjecture}{\begin{Italic}{\sc Conjecture: }}{\end{Italic}}
\newenvironment{Proof}{\par \noindent{\sc Proof:} }{\xspace\qed \par}{}
\newenvironment{Example}{\begin{Italic} \noindent{\sc Example: }}
{\end{Italic}}
\newenvironment{Remark}{\begin{Italic} \noindent{\sc Remark: }}{\end{Italic}}
\newenvironment{Remarks}{\begin{Italic} \noindent{\sc Remarks:
}}{\end{Italic}}

\newenvironment{Background}{\begin{Italic} \noindent{\sc Background:
}}{\end{Italic}}

\newcounter{enum}

\renewcommand{\qed}{$\blacksquare$}
\newcommand{\C}{\mathbb{C}}
\newcommand{\R}{\mathbb{R}}
\newcommand{\Q}{\mathbb{Q}}
\newcommand{\iso}{\cong}
\newcommand{\htp}{\simeq}

\newcommand{\JM}{\scriptscriptstyle{\mathrm{JM}}}
\newcommand{\Z}{\mathbb{Z}}
\newcommand{\reg}{\scriptscriptstyle{\mathrm{reg}}}
\newcommand{\mult}{\scriptscriptstyle{\mathrm{mult}}}
\newcommand{\sub}{\scriptscriptstyle{\mathrm{sub}}}
\newcommand{\resc}{\scriptscriptstyle{\mathrm{resc}}}
\newcommand{\half}{\frac{1}{2}}

\newcommand{\Kh}{\mathit{Kh}}

\renewcommand{\gg}{\mathfrak{g}}
\renewcommand{\sl}{\mathfrak{sl}}
\newcommand{\gl}{\mathfrak{gl}}
\newcommand{\hh}{\mathfrak{h}}
\newcommand{\uu}{\mathfrak{u}}

\newcommand{\symp}{\mathrm{symp}}
\newcommand{\Conf}{\mathit{Conf}}

\newcommand{\Slice}{\EuScript{S}}
\newcommand{\TT}{\EuScript{T}}
\newcommand{\LL}{\EuScript{L}}

\newcommand{\Y}{\EuScript{Y}}
\newcommand{\OO}{\EuScript{O}}
\newcommand{\FF}{\EuScript{F}}

\newcommand{\CC}{\EuScript{C}}

\newcommand{\minushorizresolution}{\;\;\includegraphics{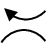}\;}
\newcommand{\plushorizresolution}{\;\;\includegraphics{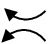}\;}
\newcommand{\minusvertresolution}{\;\;\includegraphics{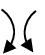}\;}
\newcommand{\plusvertresolution}{\;\;\includegraphics{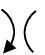}\;}
\newcommand{\minuscrossing}{\;\;\includegraphics{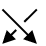}\;}
\newcommand{\pluscrossing}{\;\;\includegraphics{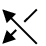}\;}

\title{A link invariant from the symplectic\\ geometry of nilpotent slices}
\author{Paul Seidel, Ivan Smith}
\date{Final version, April 2006}

\begin{document}
\maketitle
\begin{abstract}
\noindent We define an invariant of oriented links in $S^3$ using the 
symplectic geometry of certain spaces which arise naturally in Lie 
theory. More specifically, we present a knot as the closure of a braid, 
which in turn we view as a loop in configuration space. Fix an affine 
subspace $\Slice_m$ of the Lie algebra $\frak{sl}_{2m}(\C)$ which is a 
transverse slice to the adjoint action at a nilpotent matrix with two 
equal Jordan blocks. The adjoint quotient map restricted to $\Slice_m$ 
gives rise to a symplectic fibre bundle over configuration space. An 
inductive argument constructs a distinguished Lagrangian submanifold
$L_{\wp_{\pm}}$  
of a fibre $\Y_{m,t_0}$ of this fibre bundle; we regard the braid $\beta$ as 
a symplectic automorphism of the fibre, and apply Lagrangian Floer 
cohomology to $L_{\wp_{\pm}}$ and $\beta(L_{\wp_{\pm}})$ inside
$\Y_{m,t_0}$. The main theorem  
asserts that this group is invariant under the Markov moves, hence 
defines an oriented link invariant. We conjecture that this invariant 
co-incides with Khovanov's combinatorially defined link homology theory, 
after collapsing the bigrading of the latter to a single grading.
\end{abstract}


\section{Introduction}

From its first introduction in \cite{jones85}, the Jones polynomial
has played a decisive role in knot theory. This invariant associates
to an oriented link $\kappa \subset S^3$ a Laurent polynomial
$V_\kappa(t^{1/2})$. It is completely characterised by $V_{unknot} =
1$ and a relation obtained from the Kauffman bracket calculus:
\begin{equation} \label{eq:skein}
\begin{aligned}
 &
 t^{-1/2} V_{\plushorizresolution}
 + t^{3v/2} V_{\plusvertresolution}
 + t^{-1} V_{\pluscrossing} = 0, \\
 &
 t^{3v/2} V_{\minushorizresolution}
 + t^{1/2} V_{\minusvertresolution}
 + t V_{\minuscrossing} = 0.
\end{aligned}
\end{equation}
The two equations are for a positive and negative crossing,
respectively. In the complement of the crossing under consideration,
take the arc which ends at the top left corner of the crossing. Then
$v$ is the signed number of crossings between this arc and the other
connected components of the complement (this term compensates for
the non-local change of orientation that occurs in one of the two
ways of resolving the crossing). \eqref{eq:skein} allows one to
successively reduce the number of crossings, which in principle is
suitable for algorithmic computation. However, despite arising in a
farrago of contexts, the geometric meaning of the Jones polynomial
has remained somewhat mysterious.

In \cite[Section 7]{khovanov98} Mikhail Khovanov introduced a
categorified Jones polynomial, which is a bigraded abelian group
$\Kh^{*,*}(\kappa)$. $\Kh^{i,j}($\textit{unknot}$)$ is $\Z$ for $i =
0$, $j = \pm 1$, and zero otherwise. The role of equations
\eqref{eq:skein} is played by long exact sequences
\begin{multline}
 \label{eq:plustriangle}
 \qquad \cdots \longrightarrow \Kh^{i,j}(\!\pluscrossing) \longrightarrow
 \Kh^{i,j-1}(\!\plushorizresolution) \longrightarrow
 \Kh^{i-v,j-3v-2}(\plusvertresolution) \qquad\qquad \\ \longrightarrow
 \Kh^{i+1,j}(\!\pluscrossing) \longrightarrow \cdots \qquad \qquad \qquad
 \qquad \quad \qquad \qquad \qquad
\end{multline}
and
\begin{multline}
 \label{eq:minustriangle}
 \quad \cdots \longrightarrow  \Kh^{i,j}(\minuscrossing) \longrightarrow
 \Kh^{i-v+1,j-3v+2}(\minushorizresolution) \longrightarrow
 \Kh^{i+1,j+1}(\minusvertresolution) \quad\quad \\ \longrightarrow
  \Kh^{i+1,j}(\minuscrossing) \longrightarrow \cdots \qquad \qquad \qquad
  \qquad \quad \qquad \qquad \qquad
\end{multline}
Starting from these, an obvious computation shows that the graded
Euler characteristic $\chi_\kappa(q) = \sum_{i,j} (-1)^i q^j
\dim(\Kh^{i,j}(\kappa) \otimes \Q)$ is, up to normalisation and
change of variables, the Jones polynomial:
\[
V_{\kappa} = \left.
\frac{\chi_{\kappa}(q)}{q+q^{-1}}\right|_{q=-t^{1/2}}
\]
$\Kh^{*,*}$ is known to be a strictly stronger invariant than
$V_\kappa$ \cite{barnatan02};  by definition it remains
algorithmically computable.
%
%

The groups $\Kh^{*,*}$ are only the starting point for a very rich
theory. First of all, one can vary their definition in many ways,
giving rise to potentially useful additional invariants, such as
Lee's spectral sequence \cite{lee02b}. More spectacularly, they can
be shown to fit into a topological quantum field theory for
two-knots in four-space \cite{khovanov00,jacobsson02,khovanov02b}.
Very recently, Rasmussen \cite{rasmussen04} has used both properties
to give a proof of Milnor's conjecture on the slice genus of torus
knots. This was previously accessible only via gauge theory
(instanton invariants \cite{kronheimer-mrowka93} originally, then
Seiberg-Witten theory \cite{ozsvath-szabo00}, and finally the
Ozsv\'ath-Szab\'o rebirth of the latter in terms of pseudoholomorphic
curves \cite{ozsvath-szabo03b}). Rasmussen's argument is explicitly
modelled on Ozsv\'ath-Szab\'o theory, and in fact he conjectures the
equality of a certain numerical invariant obtained from $\Kh^{*,*}$
and its geometric counterpart \cite[p.\ 2]{rasmussen04}. This use of
Khovanov homology as a combinatorial substitute for gauge theory
does not come as a complete surprise. The structural resemblance of
the two theories, for instance looking at
\eqref{eq:plustriangle},\eqref{eq:minustriangle} versus Floer's
exact triangle, had been noticed for some time, and has found a
concrete expression in the spectral sequence from
\cite{ozsvath-szabo03}, which goes from a variant of Khovanov
homology to the Heegard Floer homology of the branched double cover.
In that picture, Khovanov homology appears as a combinatorial
approximation to Heegard Floer homology, and does not itself take on
a geometric meaning.

The approach proposed in this paper is different, in that we give a
tentative symplectic geometry description of Khovanov homology
itself. The construction is fairly involved, but we can give a very
rapid and superficial sketch at once: following Jones and others, we
present an oriented link $\kappa$ as a closure of an $m$-stranded
braid $b \in Br_m$. Adding another $m$ trivial strands gives $b
\times 1^m \in Br_{2m}$, which can be represented as a loop $\beta$
in configuration space $\Conf_{2m}(\C)$ with respect to some base
point $t_0$. We introduce a $4m$-dimensional noncompact symplectic
manifold $M = \Y_{m,t_0}$ (here and below, the more complicated
notation is the one used in the body of the paper). This is the
fibre at $t_0$ of a symplectic fibration over configuration space,
whose monodromy along $\beta$ yields a symplectic automorphism $\phi
= h^{\resc}_{\beta}$. Our $M$ also contains a canonical (up to
isotopy) Lagrangian submanifold $L = L_{\wp_{\pm}}$ diffeomorphic to
$(S^2)^m$. We apply Lagrangian Floer cohomology to these geometric
data, and set
\[
 \Kh_{\symp}^*(\kappa) = HF^{*+m+w}(L,\phi(L))
 = HF^{*+m+w}(L_{\wp_{\pm}},h^{\resc}_{\beta}(L_{\wp_{\pm}}))
\]
where $m$ is as before, and $w$ is the writhe of the braid
presentation (the number of positive minus the number of negative
crossings).

\begin{Theorem} \label{th:main}
Up to isomorphism of graded abelian groups, $\Kh_{\symp}^*(\kappa)$
depends only on the oriented link $\kappa$, and not on its
presentation as a braid closure.
\end{Theorem}

Reversing the orientation of all components leaves $\Kh_{\symp}^*$
unchanged, in particular we get an invariant of unoriented knots.
The proposed relation with Khovanov homology is that our invariant
should be obtained from it by collapsing the bigrading (actually,
the sign of the $j$-grading should be reversed first, which is a
simple change of conventions already applied in \cite{khovanov00}).

\begin{Conjecture} \label{th:conj}
\[
\Kh_{\symp}^k(\kappa) \iso \bigoplus_{i-j = k} \Kh^{i,j}(\kappa).
\]
\end{Conjecture}

We admit at once that, at least for the time being, the symplectic
theory does not come with a bigrading corresponding to the one in
$\Kh^{*,*}$. This prevents us from seeing the connection to the
Jones polynomial geometrically, since the Euler characteristic only
recovers the uninteresting specialization $V_{\kappa}(t^{1/2} = 1)$
which counts components of the link. Evidence for Conjecture
\ref{th:conj} comes from various sources. The two sides have the
same value for the unknot, and more generally they behave in the
same way under adding an unlinked unknot component. Another example
is provided by the trefoil knot, whose $\Kh^*_{\symp}$ we will
compute by direct geometric means. More speculatively, a
generalization of \cite{seidel01} should yield the counterpart for
$\Kh^*_{\symp}$ of the long exact sequences \eqref{eq:plustriangle},
\eqref{eq:minustriangle} (but there are many details still to be
carried out). Starting from this, one could follow the construction
of the ``hypercubes of crossing resolutions'' of \cite{khovanov98}
in $\Kh^*_{\symp}$ and thereby obtain a spectral sequence which
starts with $E_2= \Kh^{*,*}$ and converges to $\Kh^*_{\symp}$ (this
is precisely the approach used by Ozsv\'ath-Szab\'o
\cite{ozsvath-szabo03}). The conjectural vanishing of the higher
order differentials in this sequence seems more difficult to explain
at present, even though the fact that both theories are $\Z$-graded
restricts the possibilities somewhat. A possible more fundamental
explanation for Conjecture \ref{th:conj} would arise from a
relation, on the derived level, between the Fukaya categories of $M
= \Y_{m,t_0}$ and differential graded modules over the arc algebras
$H_m$ from \cite{khovanov00}. This was one of Khovanov's motivations
when he proposed that these particular manifolds should be relevant
for understanding $\Kh^{*,*}$ \cite{khovanov02c}. Because of its
abstract homological algebra nature, this approach may seem
far-fetched, but it has been successfully carried out in a toy model
case \cite{khovanov-seidel98,seidel-thomas99}.

After this rather tentative discussion, we return to the concrete
geometry underlying the definition of $\Kh^*_{\symp}$. The manifolds
$M = \Y_{m,t_0}$ can be described in elementary terms (as given by
matrices of a certain special form, and with prescribed
eigenvalues), but the proper framework for understanding them is
provided by Lie theory. For any semisimple complex Lie algebra
$\gg$, one can consider the adjoint quotient map $\chi: \gg
\rightarrow \hh/W$. In the case of $\gg = \sl_{2m}$, which is the
one relevant to us here, this associates to each matrix its
characteristic polynomial. Thinking of $\hh/W$ as the set of
unordered eigenvalues with multiplicities, we can identify an open
dense subset $\hh^{\reg}/W$ with the space $\Conf_{2m}^0(\C)$ of
configurations having zero center of mass, and the restriction of
$\chi$ to that subspace is a differentiable fibre bundle. We
actually want to restrict $\chi$ to a so-called transverse slice,
which is an affine subspace of $\gg$ intersecting all orbits of the
natural adjoint $G$-action transversally. There is a well-known
general construction of such slices $\Slice^{\JM} \subset \gg$,
which starts with a nilpotent element of $\gg$ and invokes the
Jacobson-Morozov theorem. The restrictions $\chi|\Slice^{\JM}:
\Slice^{\JM} \rightarrow \hh/W$ are still differentiable fibre
bundles over the subset $\hh^{\reg}/W$, and their topological
monodromy has been used by Slodowy \cite{slodowy} and others to give
an alternative construction of Springer's Weyl group
representations. For our particular slices $\Slice_m$, we take the
nilpotent $n^+ \in \gg$ which has two Jordan blocks of size $m$, and
use a slight generalization of the Jacobson-Morozov construction
(this is purely for technical reasons: $\Slice_m$ is isomorphic to
the corresponding Jacobson-Morozov slice $\Slice^{\JM}$, but is
slightly easier to use). The fibre of $\chi|\Slice_m$ at the point
$t_0$ (translated by a constant to put it into $\Conf_{2m}^0(\C)$,
to be precise) is our $\Y_{m,t_0}$.

Crucially, if we take a point $t \in \hh/W$ where two eigenvalues
come together, the fibre of $\chi|\Slice_m$ over $t$ has a fibered
$(A_1)$ singularity. By this we mean that the stratum of singular
points is itself smooth, and that the normal structure of this
stratum is that of an ordinary double point singularity $a^2 + b^2 +
c^2 = 0$. Moreover, up to symplectic isomorphism the singular
stratum itself can be identified with a regular fibre of
$\chi|\Slice_{m-1}$. This leads to an inductive scheme, where we
construct Lagrangian submanifolds in $\Y_{m,t_0}$ by bringing the
eigenvalues together in successive pairs. More precisely, each
``crossingless matching'' $\wp$ of the configuration $t_0$ in the
plane gives rise to a Lagrangian submanifold $L_\wp$, and the $L =
L_{\wp_{\pm}}$ which appeared above is just the example obtained
from a certain standard choice of matching.  To prove Theorem
\ref{th:main} one needs to check the invariance of our symplectic
Floer cohomology group under the Markov moves which relate different
braid presentations of the same link. For the most difficult
(because it changes the number $m$ of strands) type II move this
will follow from the observation that if $t$ is a point where three
eigenvalues are brought together, the fibre of $\chi|\Slice_m$ at
that point has a fibered $(A_2)$ type singularity. The main lesson
to be learnt from this is that the data which enter into the
definition of $\Kh^*_{\symp}$, and the properties which make it an
oriented link invariant, are all derived from the basic geometry of
the adjoint quotient, and general facts from symplectic geometry. To
emphasize that, the paper alternates sections of general exposition
with others more specifically tailored to our needs.

Alternatively, one can think of $\Y_{m,t_0}$ as a space of solutions
of Nahm's equations with certain boundary data (these are basically
the equations for $\R^3$-invariant instantons on $\R^4$, with gauge
group $SU_{2m}$, cf. \cite{atiyah-bielawski02} generalizing
\cite{kronheimer89}). Via an ADHM transform $\Y_{m,t_0}$ can be
viewed as a quiver variety, cf. \cite[Theorem 8.4]{nakajima94}
(strictly this is partly conjectural, since the cited work does not
quite cover our case). These descriptions show geometric properties
of $\Y_{m,t_0}$ which may not be immediately apparent from the point
of view taken here, and lead to several possible avenues for further
development (other than the obvious question of trying to prove
Conjecture \ref{th:conj}). For instance, there is an involution on
$\Y_{m,t_0}$ whose fixed point set is related to the Jacobian of the
double cover of $\C\mathbb{P}^1$ branched along the configuration of
points $t_0\in\C$, and that leads to a geometric
interpretation of the relation between Khovanov and Ozsv\'ath-Szab\'o
theory. This will be discussed in detail in a sequel.

\textbf{Acknowledgements.} We are enormously indebted to Mike
Khovanov for generously sharing his ideas; this project would not
have got off the ground without his help.

\section{Geometry of the adjoint quotient}

This section gathers some background material from geometric Lie
theory. All of it is essentially extracted from the textbooks
\cite{collingwood-mcgovern,chriss-ginzburg} and more advanced works
\cite{slodowy,slodowy82,kronheimer90,rossmann95}. To make the
exposition more focused and easier for non-specialists to read, we
only deal with the Lie algebras $\gg = \sl_n$. In this case, most
proofs are within the reach of elementary linear algebra, and we
have omitted some of them. However, each section ends with some
brief remarks which outline the situation for general semisimple
$\gg$, and gives more references to the literature.
%

\subsection{Local slices}
Let $G = SL_n(\C)$, and $\gg = \sl_n$ its Lie algebra. The adjoint
quotient map $\chi$ associates to each matrix $x \in \gg$ its
characteristic polynomial. Since that polynomial has $(n-1)$
nontrivial coefficients, the adjoint quotient is a holomorphic map
$\chi: \gg \rightarrow \C^{n-1}$. Equivalently, one can think of it
as giving the eigenvalues of the matrix (with multiplicities). From
that perspective, $\chi(x) \in \hh/W$, where $\hh \subset \gg$ is
the subspace of diagonal matrices, and $W = S_n$ the permutation
group acting on it. The connection between the two points of view is
established by the elementary symmetric functions, which give a
holomorphic isomorphism $\hh/W \iso \C^{n-1}$.

\begin{Example} \label{ex:sl2}
Take $\gg = \sl_2$. The adjoint quotient is just the determinant,
which for a suitable choice of coordinates on $\gg$ can be written
as
\[
 \chi: \C^3 \longrightarrow \C, \quad (a,b,c) \longmapsto
 a^2+b^2+c^2.
\]
It has a single nondegenerate critical point at the origin, so the
fibre $\chi^{-1}(0)$ has an ordinary double point singularity.
Following common terminology, we will call these $(A_1)$ type
critical point and singularity, respectively.
\end{Example}

$G$ acts on $\gg$ by conjugation, $Ad(g)y = gyg^{-1}$. The
corresponding infinitesimal action of $\gg$ on itself is $ad(x)y =
[x,y]$. Either of these is usually called the adjoint action. $\chi$
is constant along $G$-orbits (and in fact can be thought of as the
projection to the algebro-geometric quotient $\gg/G \iso \hh/W$).

\begin{Example} \label{ex:different-ev}
Suppose that $x \in \gg$ has $n$ pairwise distinct eigenvalues.
Since it is semisimple and $\chi$ is $G$-invariant, we may assume
that $x \in \hh$. Moreover, since the eigenvalues are distinct, $x$
is a regular point of the projection $\hh \rightarrow \hh/W$. This
implies that $x$ is a regular point of $\chi$.
\end{Example}

Because of the $G$-symmetry, the local geometry of $\chi$ can be
studied in terms of transverse slices, which we will now define. Fix
$x \in \gg$. By definition, the tangent space to the adjoint orbit
$Gx$ is
\[
T_x(Gx) = ad(\gg)x = [x,\gg].
\]
A local transverse slice to the orbit is simply a local complex
submanifold $\Slice \subset \gg$, $x \in \Slice$, whose tangent
space at $x$ is complementary to $T_x(Gx)$. Take such a slice, and
in addition, let $K \subset G$ be a local submanifold containing the
identity $e \in G$, such that $T_eK$ is complementary to the
stabilizer $\gg_x = \{y \;:\; [y,x] = 0\}$. Then we get a
commutative diagram
\begin{equation} \label{eq:basic}
\begin{CD}
K \times \Slice @>{Ad|K \times \Slice}>> \gg \\
@V{\text{projection}}VV @VV{\chi}V \\
\Slice @>{\chi|\Slice}>> \hh/W
\end{CD}
\end{equation}
where the top $\rightarrow$ map, $(g,y) \mapsto Ad(g)y$, is a local
holomorphic isomorphism at $(e,y)$. In words, this means that $\chi$
looks locally like $\chi|\Slice$ times a constant map in the
remaining coordinates.

\begin{Lemma} \label{th:stupid}
(i) For all $y \in \Slice$ sufficiently close to $x$, the
intersection $\Slice \cap Gy$ is transverse at $y$. (ii) For all $y
\in \Slice$ sufficiently close to $x$, we have that $y$ is a
critical point of $\chi|\Slice$ iff it is a critical point of
$\chi$. (iii) Any two local transverse slices at $x$ are locally
isomorphic, by an isomorphism which moves points only inside their
$G$-orbits.
\end{Lemma}

\proof (i) is clear from the definition. (ii) follows from
\eqref{eq:basic}. For (iii), if $\Slice'$ is another slice, then the
desired isomorphism is
\begin{equation}
\label{eq:stupid-id}
 \Slice' \xrightarrow{\text{inclusion}}
 \gg \xleftarrow{Ad|K \times \Slice} K \times \Slice
 \xrightarrow{\text{projection}} \Slice. \qquad \text{\qed}
\end{equation}
%

\begin{Example} \label{ex:regular1}
Let $x$ be the nilpotent consisting of a single maximal Jordan
block,
\begin{equation} \label{eq:standard-jordan}
 x = \begin{pmatrix}
 0 & 1 &&& \\
  &&& \dots & \\
  &&& \dots & 1 \\
 &&& & 0
 \end{pmatrix}.
\end{equation}
Then the space $\Slice$ of matrices
\begin{equation} \label{eq:standard-slice}
 y = \begin{pmatrix}
 0 & 1 &&& \\
 y_{21} && 1 && \\
 \dots &&& \dots & \\
 y_{n-1,1} &&&& 1 \\
 y_{n1} &&&& 0
 \end{pmatrix}
\end{equation}
is a slice at $x$, and $\chi|\Slice = id_{\C^{n-1}}$. In view of
Lemma \ref{th:stupid}(ii) above, it follows that $y$ is a regular
point of $\chi$ itself.
\end{Example}

\begin{Example} \label{ex:sl3}
In $\gg = \sl_3$, the space $\Slice$ of matrices
\begin{equation} \label{eq:sl3-slice}
 y = \begin{pmatrix} \alpha & 0 & 1 \\ \beta & -2\alpha & 0 \\
 \delta & \gamma & \alpha \end{pmatrix}
\end{equation}
is a slice at $x = \{\alpha = \beta = \gamma = \delta = 0\}$. After
changing coordinates to $\alpha = \half a$, $\beta = b$, $\gamma =
-c$, $\delta = \frac{3}{4}a^2 - d$, the characteristic polynomial of
$y$ is $t^3 - td + (a^3 - ad + bc)$, so one can write
\begin{equation} \label{eq:a2-map}
 \chi|\Slice: \C^4 \longrightarrow \C^2, \quad \chi(a,b,c,d) =
 (d,a^3-ad+bc).
\end{equation}
This map is known to singularity theorists as the miniversal
unfolding of the $(A_2)$ type surface singularity $a^3 + bc = 0$.
\end{Example}

Both examples above describe slices at nilpotent points. These are
particularly important because the geometry of the general case can
be reduced to the nilpotent one, as we will now explain. Consider an $x
\in \gg$ whose eigenvalues are $(\mu_1,\dots,\mu_n)$. For
simplicity, assume that the first $k$ eigenvalues are equal, and
that there are no other coincidences between them. Write $x$ as the
sum of its semisimple and nilpotent parts, $x = x_s + x_n$. The
stabilizer $\gg_{x_s} = \{y \;:\; [y,x_s] = 0\} \subset \gg$ is a
Lie subalgebra of block diagonal matrices. To write this down
explicitly, let $E$ be the $\mu_1$-eigenspace of $x_s$, and
$L_{k+1},\dots,L_n$ the remaining eigenspaces (which are
one-dimensional by assumption). Then $\gg_{x_s}$ is the trace-free
part of $\gl(E) \oplus \gl(L_{k+1}) \cdots \oplus \gl(L_n)$. One can
further decompose this as
\begin{equation} \label{eq:zplit}
\gg_{x_s} = \hat\gg \oplus {\mathfrak z}.
\end{equation}
Here, the first factor is $\hat\gg = \sl(E)$, while the second one
${\mathfrak z}$ is the center, consisting of the trace-free part of
$\{\C \cdot 1 \subset \gl(E)\} \oplus \gl(L_{k+1}) \oplus \cdots
\oplus \gl(L_n)$. Note that ${\mathfrak z}$ can be identified with
$\C^{n-k} \subset \C^{n-k+1}$ in an obvious way, without choosing
bases for our eigenspaces. The nilpotent part $x_n$ naturally lies
in $\hat\gg$. Suppose that we are given a slice $\hat\Slice \subset
\hat\gg$ to $x_n$, with respect to the adjoint action on that
smaller Lie algebra. An explicit comparison of $[x,\gg]$ with
$[x_n,\hat\gg]$ shows that

\begin{Lemma} \label{th:make-slice}
$\Slice = x_s + \hat\Slice + {\mathfrak z} \subset \gg$ is a local
transverse slice for the adjoint action at $x$. \qed
\end{Lemma}

Write $\hat\chi: \hat\gg \rightarrow \hat\hh/\hat{W}$ for the
adjoint quotient map on $\hat\gg$. Take the isomorphism $\hat\Slice
\times {\mathfrak z} \rightarrow \Slice$, $(y,z) \mapsto x_s+y+z$.
On the bases, consider the map $\hat{\hh}/\hat{W} \times {\mathfrak
z} \rightarrow \hh/W$ which takes $\lambda,\rho \in \hat{\hh} \times
{\mathfrak z}$ and adds up the three collections of eigenvalues
$(\lambda_1,\dots,\lambda_k,0,\dots,0)$, $(\rho_1,\dots,\rho_1,
\rho_{k+1},\dots,\rho_n)$ and $(\mu_1,\dots,\mu_n)$. This is a local
isomorphism near the origin, since the subgroup of $W$ fixing
$(\mu_1,\dots,\mu_n)$ is precisely $\hat{W}$. Together, these maps
fit into a commutative diagram
\begin{equation} \label{eq:reduction}
\begin{CD}
 \hat\Slice \times {\mathfrak z} @>{\text{isomorphism}}>>
 \Slice \\
 @V{(\hat\chi|\hat\Slice)\times Id}VV @VV{\chi|\Slice}V \\
 {\hat{\hh}/\hat{W} \times {\mathfrak z}} @>{\text{local isomorphism}}>> {\hh/W.}
\end{CD}
\end{equation}
As promised, this means that the local structure of $\chi|\Slice$
reduces to that of $\hat\chi|\hat\Slice$ times the identity map in
the remaining coordinates.

\begin{Example} \label{th:2-block}
Suppose that $x$ as above, with eigenvalues $(\mu_1,\ldots,\mu_n)$, is itself semisimple, so that $x_s = x$ and $x_n =
0$. Then $\hat\Slice = \hat\gg = \sl(E)$, hence $\Slice = \gg_x$.
Now assume in addition that $k = 2$. By choosing coordinates on
$\sl(E)$ as in Example \ref{ex:sl2}, and using \eqref{eq:reduction},
we get the following local picture of $\chi|\Slice$:
\begin{equation} \label{eq:fibered-a1-trivial}
 \chi(a,b,c,\rho_1,\dots,\rho_{n-2}) = (a^2+b^2+c^2,\rho_1,\dots,\rho_{n-2}).
\end{equation}
To make this even simpler, take a disc $D \subset \hh/W$
corresponding to eigenvalues
$(\mu_1-\sqrt{\epsilon},\mu_2+\sqrt{\epsilon},\mu_3,\dots,\mu_n)$
for small $\epsilon$. By \eqref{eq:reduction}, $\chi^{-1}(D) \cap
\Slice$ just singles out the $\hat\Slice$-factor of the slice.
Hence, the restriction of $\chi$ to $\chi^{-1}(D) \cap \Slice$ has
an $(A_1)$ type critical point at $x$.
\end{Example}

The case of a general $x$, where several eigenvalues may have
nontrivial multiplicities $k_1,\dots,k_r > 1$, is analogous. One
splits $\gg_{x_s}$ into $\hat\gg \iso \sl_{k_1} \times \cdots
\sl_{k_r}$ and the center ${\mathfrak z} \iso
\C^{n-k_1-\cdots-k_r+r-1}$. For each $j = 1,\dots,r$, choose a slice
$\hat\Slice_j \subset \sl_{k_j}$ at the corresponding component of
$x_n$. By combining these with ${\mathfrak z}$, one again obtains a
slice $\Slice$ for $x$. The outcome is that $\chi|\Slice$ looks like
the product of the maps $\hat\chi_j|\hat\Slice_j$ for $j =
1,\dots,k$, together with the identity map on ${\mathfrak z}$.

\begin{Example} \label{th:all-regular}
Suppose that $x \in \gg$ is a matrix which has a single Jordan block
for each eigenvalue. This means that for each $j$, the
$\sl_{k_j}$-component of $x_n$ is as in Example \ref{ex:regular1}.
We know that these are regular points of the adjoint quotient maps
$\hat\chi_j$, hence $x$ itself is a regular point of $\chi$ (these
are in fact all the regular points).
\end{Example}

\begin{Example} \label{th:2-2-block}
Let $x \in \gg$ be a semisimple matrix with eigenvalues
$(\mu_1,\dots,\mu_n)$, where the first $2k$ form equal pairs $(\mu_1
= \mu_2, \mu_3 = \mu_4,\cdots,\mu_{2k-1} = \mu_{2k})$, and with no
other coincidences. This is quite similar to Example
\ref{th:2-block}: $\Slice = \gg_x$, and each $\hat\Slice_j$ is the
whole of $\sl_2$. Consider the polydisc $P \subset \hh/W$ formed by
the sets of eigenvalues
$(\mu_1-\sqrt{\epsilon_1},\mu_2+\sqrt{\epsilon}_1,\dots,
\mu_{2k-1}-\sqrt{\epsilon_k},\mu_{2k}+\sqrt{\epsilon}_k,\mu_{2k+1},\dots,
\mu_n)$. Then, the restriction of $\chi$ to $\chi^{-1}(P) \cap
\Slice$ looks locally like the product of $k$ copies of the $(A_1)$ type map.
\end{Example}

\begin{Remarks}
For a general semisimple Lie group $G$, with Lie algebra $\gg$, one
defines the adjoint quotient as the projection $\chi: \gg
\rightarrow \gg/G \iso \hh/W$, where $\hh$ is a Cartan subalgebra
and $W$ the Weyl group (the isomorphism is Chevalley's theorem
$\C[\gg]^G = \C[\hh]^W$, see e.g.\ \cite[Theorem
3.1.38]{chriss-ginzburg}). The basic facts about slices, notably
\eqref{eq:basic} and Lemma \ref{th:stupid}, continue to hold, since
they actually apply to general holomorphic $G$-actions and
$G$-invariant maps. Examples \ref{ex:sl2} and \ref{ex:sl3} are
special cases of a fundamental result of Brieskorn and Slodowy
\cite{brieskorn70b,slodowy}, which says that if $\Slice$ is a slice
to a subregular nilpotent element inside a simply-laced $\gg$, then
$\chi|\Slice$ is the miniversal unfolding of the corresponding
simple (ADE type) singularity.

The general version of Lemma \ref{th:make-slice} looks as follows.
Any $x \in \gg$ has a canonical decomposition $x = x_s + x_n$ into
mutually commuting semisimple and nilpotent parts. Since $x_s$ is
semisimple, $\gg_{x_s}$ is a reductive Lie algebra \cite[Lemma
2.1.2]{collingwood-mcgovern}, hence splits into its center
${\mathfrak z}$ and the semisimple derived Lie algebra $\hat\gg =
[\gg_{x_s},\gg_{x_s}]$, as in \eqref{eq:zplit}. Moreover, $ad(x_s)$
is a semisimple endomorphism of $\gg$, so $\gg = [x_s,\gg] \oplus
\gg_{x_s}$. Finally, $ad(x_s)$, $ad(x_n)$ are polynomials in $ad(x)$
with zero constant terms, hence $[x,\gg] = [x_s,\gg] + [x_n,\gg]$.
Taking these three facts together, one finds that if $\hat\Slice
\subset \hat\gg$ is a $\hat{G}$-slice at the point $x_n$, then
$\Slice = x_s + {\mathfrak z} + \hat\Slice$ is a $G$-slice at $x$.
Restricting the adjoint quotient to this slice, one obtains a
diagram like \eqref{eq:reduction}. For instance, $x$ is a regular
point of $\chi$ iff $x_n \in \hat\gg$ is a regular nilpotent, which
is the general version of Example \ref{ex:regular1}.
\end{Remarks}

\subsection{Homogeneous slices}

Let $x \in \gg$ be nilpotent. The Jacobson-Morozov Lemma says that
one can find a triple $(n^+,n^-,h)$ of elements of $\gg$, where $n^+
= x$, which satisfy
\begin{equation} \label{eq:jm}
[h,n^+]=2n^+, \quad [h,n^-]=-2n^-, \quad [n^+,n^-]=h.
\end{equation}
There are many different choices of $(n^-,h)$ for a fixed $n^+$.
However, Kostant's uniqueness theorem says that any two are
conjugate by an element of the stabilizer $G_{n^+} \subset G$.

The elements of a Jacobson-Morozov (JM) triple define a homomorphism
$\sl_2 \rightarrow \gg$, which in combination with the adjoint
action makes $\gg$ into an $\sl_2$-module. This allows one to apply
elementary facts about $\sl_2$-representations. For instance, $ad(h)
\in End(\gg)$ is necessarily semisimple and has integer eigenvalues.
Hence, the vector field on $\gg$ defined by
\begin{equation} \label{eq:k-field}
K_y = -[h,y]+2y
\end{equation}
generates a linear $\C^*$-action $\lambda_r(y) = r^2 exp(-log(r)h)
\,y\,exp(log(r)h)$. Via the adjoint quotient map, this is compatible
with the $\C^*$-action on $\hh/W$ which multiplies all eigenvalues
by $r^2$. By definition $K_x = 0$, so $x$ is a fixed point of
$\lambda$. We define a homogeneous slice at $x$ to be an affine
subspace $\Slice \subset \gg$ invariant under $\lambda$, which is a
local transverse slice for the adjoint action at $x$. There is
actually a canonical choice of homogeneous slice, namely the JM
slice $\Slice^{\JM} = x + \gg_{n^-}$. The fact that this is a slice,
or equivalently that $\gg = [n^+,\gg] \oplus \gg_{n^-}$, is another
easy observation from $\sl_2$-representation theory.

\begin{Example} \
If $x$ is as in Example \ref{ex:regular1}, one can take
\begin{equation} \label{eq:standard-jm}
\begin{aligned}
 h & = \begin{pmatrix}
 n-1  &&&& \\
 & n-3 &&& \\
 && n-5 && \\
 && & \dots & \\
 &&&& -n+1 \end{pmatrix}, \\
 n^- & = \begin{pmatrix}
 0 &&&&& \\
 n-1 &&&&& \\
 & 2(n-2) &&&& \\
 && 3(n-3) &&& \\
 &&& \dots && \\
 &&&& n-1 & \;\; 0 \;\;\;
 \end{pmatrix}.
\end{aligned}
\end{equation}
The associated $\C^*$-action is
\begin{equation} \label{eq:standard-lambda}
 \lambda_r: y \longmapsto \begin{pmatrix}
 r^2 y_{11} & y_{12} & \dots && r^{4-2n} y_{1n} \\
 r^4 y_{21} & r^2 y_{22} \\
 \dots && \dots &&  \\
 &&&& y_{n-1,n} \\
 r^{2n} y_{n1} && \dots & r^4 y_{n,n-1} & r^2 y_{nn}
 \end{pmatrix}.
\end{equation}
By listing the eigenvalues of $ad(h)$, one sees that that the
$\sl_2$-module $\gg$ breaks up into indecomposables of rank
$3,5,\dots,2n-1$ (one each). The subspace $\gg_{n^-}$ is therefore
$(n-1)$-dimensional, which means that it is spanned by powers of
$n^-$, so $\Slice^{JM} = n^+ + (\C n^- \oplus \C (n^-)^2 \oplus
\cdots)$. One should compare this with the slice considered in
Example \ref{ex:regular1}, which is homogeneous for the same choice
of $(h,n^-)$, but not a JM slice.
\end{Example}

\begin{Example} \label{ex:sl3-again}
Take the nilpotent $x \in \gg = \sl_3$ from Example \ref{ex:sl3}.
The slice constructed there is a JM slice, obtained by taking
\begin{equation} \label{eq:sl3-triple}
h = \begin{pmatrix} 1 && \\ & 0 &\!\! \\ && \!\!-1 \end{pmatrix},
\;\; n^- =
\begin{pmatrix} && \\ && \\ 1 & \; & \; \end{pmatrix}.
\end{equation}
By Kostant's theorem, any other JM triple with $n^+ = x$ is the
conjugate of this one by an element of $G_{n^+}$, which is the group
of matrices
\begin{equation} \label{eq:g}
 g = \begin{pmatrix} \delta & \epsilon & \kappa
 \\ & \delta^{-2} & \tau \\ & & \delta
 \end{pmatrix}.
\end{equation}
Note that those $g$ which are diagonal preserve $n^-$, hence the
entire triple. Thus, the space of all JM triples with $n^+ = x$
becomes isomorphic to $G_{n^+}/G_{n^+} \cap G_{n^-}$, or
equivalently, to the subgroup $U \subset G_{n^+}$ of unipotent
matrices. To interpret this in a more geometric way, note that $x$
gives rise to a flag in $\C^3$,
\[
0 \subset F^1 = im(x) \subset F^2 = ker(x) \subset \C^3.
\]
Given a JM triple, the eigenspaces of $h$ yield a splitting of this
flag into one-dimensional spaces. This is obvious for
\eqref{eq:sl3-triple}, and follows for general triples by Kostant's
theorem. $U$ acts simply transitively on the space of such
splittings, and this proves that choices of splittings and JM
triples correspond to each other bijectively.
\end{Example}

We will now look at the general properties of homogeneous slices
$\Slice$. Decompose $\gg$ and $T_x\Slice$ into $ad(h)$-eigenspaces
$\gg^{(j)}$, $(T_x\Slice)^{(j)}$, $j \in \Z$. Again appealing to
basic facts about $\sl_2$-representations, we see that the map
$ad(x)^{(j)}: \gg^{(j)} \rightarrow \gg^{(j+2)}$ is injective for
$j<0$, and surjective for $j>-2$. By definition, $(T_x\Slice)^{(j)}$
is a complementary subspace to the image of $ad(x)^{(j-2)}$, hence
is zero for all $j>0$. This implies that every homogeneous slice is
necessarily contained in the subset
\[
 \TT = x + \bigoplus_{j < 2} \gg^{(j)}
\]
of those $y \in \gg$ which satisfy $\lim_{r \rightarrow 0}
\lambda_r(y) = x$. In words, the $\C^*$-action shrinks the slice to
the point $x$. This immediately leads to an improved version of the
first two parts of Lemma \ref{th:stupid}:

\begin{Lemma} \label{th:global-slice}
Let $\Slice$ be a homogeneous slice for $x$. Then, (i) the
intersection of $\Slice$ with any adjoint orbit is transverse. (ii)
A point of $\Slice$ is a critical point of $\chi$ iff it is a
critical point of $\chi|\Slice$. \qed
\end{Lemma}

Define $\uu = \bigoplus_{j<0} \gg^{(j)}$, and let $U = exp(\uu)
\subset G$ be the corresponding subgroup. For any $u \in \uu$, the
vector field defined by $L_y = [u,y]$ is tangent to $\TT$. Hence,
the adjoint action of $U$ preserves $\TT$, in particular we get a
map
\begin{equation} \label{eq:ad-u-s}
Ad|(U \times \Slice): U \times \Slice \longrightarrow \TT.
\end{equation}
The derivative of this at the point $(e,x)$ is $\uu \oplus T_x\Slice
\longrightarrow T_x\TT$, $(u,y) \longmapsto [u,x] + y$. In view of
the observations made above, this is invertible. Equip $U \times
\Slice$ with the $\C^*$-action $r \cdot (g,y) = (exp(-log(r)h)\cdot
g\cdot exp(log(r)h),\lambda_r(y))$, and $\TT$ with $\lambda$, so
that the map between the two becomes equivariant. Since both actions
contract the relevant spaces to the point $(e,x)$ respectively its
image $x$, it follows that \eqref{eq:ad-u-s} is a global
$\C^*$-equivariant isomorphism.

\begin{Lemma} \label{th:invariant-slice}
Let $\Slice,\Slice'$ be two homogeneous slices at $x$, possibly
defined using different JM triples. Then, there is an isomorphism
$\Slice \rightarrow \Slice'$ which is $\C^*$-equivariant, and which
moves points only inside their adjoint orbits.
\end{Lemma}

\begin{Proof}
Suppose first that our two slices share the same underlying JM
triple. From \eqref{eq:ad-u-s} we then get isomorphisms
\begin{equation} \label{eq:kronheimer}
\Slice \longrightarrow \TT/U \longleftarrow \Slice'.
\end{equation}
On the other hand, if we are considering Jacobson-Morozov slices
associated to different choices of $(n^-,h)$, then the conjugating
element provided by Kostant's theorem directly yields the
isomorphism in question. In both cases, all the desired properties
are obvious; and by combining them, the general result follows.
\end{Proof}

\begin{Remark}
For the Jacobson-Morozov Lemma in the setting of general semisimple
$\gg$, see \cite[Theorem 3.3.1]{collingwood-mcgovern}. Entirely in
parallel with the case $\gg = \sl_n$, this leads to a definition of
the $\C^*$-action $\lambda$, and to the notion of homogeneous
slices.
Among these, Jacobson-Morozov slices are the most commonly used ones
in the literature.
The description of these slices as quotients \eqref{eq:kronheimer},
which uses $h$ but not $n^-$, is due to Kronheimer \cite[Lemma
11]{kronheimer90}.

The general version of Kostant's theorem can be found e.g.\ in
\cite[Theorem 3.4.10]{collingwood-mcgovern}. The statement is
actually a little better than the version given here, since it says
that the conjugating elements can be taken to lie in a certain
unipotent subgroup of $G_{n^+}$ (this is quite visible in Example
\ref{ex:sl3-again}, for instance).
\end{Remark}

\subsection{Simultaneous resolution}

Consider the open subset $\hh^{reg}/W \subset \hh/W$ corresponding
to $n$-tuples of pairwise different eigenvalues. We will identify
this with the subspace $\Conf_n^0(\C) \subset \Conf_n(\C)$ of point
configurations with zero center of mass. Each $t \in \Conf_n^0(\C)$
is a regular value of $\chi$, by Example \ref{ex:different-ev} (in
fact, these are all the regular values). Therefore, the part of the
adjoint quotient lying over $\Conf_n^0(\C)$ is a submersion. We need
to show that it is in fact a differentiable fibre bundle. The
technical difficulty is that the fibres are not compact, and we will
resolve this by using Grothendieck's simultaneous resolution and a
suitable $\C^*$-action.

Let $\tilde{\gg}$ be the space of pairs $(x,F)$, where $x \in \gg$
and $F$ is a complete flag such that $x(F^i) \subset F^i$ for all
$i$. Since the flag manifold is compact, projection $\tilde\gg
\rightarrow \gg$ is a proper map. Next, note that for $(x,F) \in
\tilde\gg$, we can consider the endomorphism of each quotient
$F^i/F^{i-1}$ induced by $x$, which is multiplication by some
$\tilde{t}_i \in \C$. The $\tilde{t}_i$ are the eigenvalues of $x$,
with the correct multiplicities, and the flag $F$ provides a
preferred ordering of them. This means that the map $\tilde\chi:
\tilde\gg \rightarrow \hh \iso \C^{n-1} \subset \C^n$,
$\tilde{\chi}(x,F) = (\tilde{t}_1,\dots,\tilde{t}_n)$, fits into a
commutative diagram
\[
\begin{CD}
 {\tilde{\gg}} @>>> \gg \\
 @V{\tilde{\chi}}VV @V{\chi}VV \\
 {\hh} @>>> {\hh/W.}
\end{CD}
\]

\begin{Lemma} \label{th:simultaneous}
(i) $\tilde\chi$ is a submersion, so each fibre
$\tilde\chi^{-1}(\tilde{t})$ is smooth. (ii) if $t$ lies in
$Conf_n^0(\C)$ and $\tilde{t} \in \hh$ is any point lying over it,
then $\tilde\chi^{-1}(\tilde{t}) \iso \chi^{-1}(t)$.
\end{Lemma}

\begin{Proof}
Both properties are elementary. For (i), since everything is
invariant under the $G$-action by conjugation, we may assume that
$F$ is the standard flag, so $x$ is a diagonal matrix. Then, by
changing only the diagonal coefficients, we get a subspace of
$T_x\tilde{\gg}$ which projects isomorphically to $\hh$. For (ii),
note that if the eigenvalues are pairwise different, any ordering of
them determines a unique compatible flag.
\end{Proof}

\begin{Example} \label{ex:sl2-cont}
We wrote down the adjoint quotient map for $\gg = \sl_2$ in Example
\ref{ex:sl2}. The simultaneous resolution is a classical
algebro-geometric construction: it consists of doing the base
extension $t^2 = a^2 + b^2 + c^2$ and then taking a small resolution
of that, which replaces the singular point with a $\C\mathbb{P}^1$.
\end{Example}

\begin{Lemma} \label{th:fibre-bundle}
$\tilde{\chi}$ is naturally a differentiable fibre bundle.
\end{Lemma}

\proof The diagonal $\C^*$-action $\rho$ on $\gg$ obviously descends
to an action on $\hh/W$. Note that both actions have positive
weights, hence contract the relevant spaces to a point (the origin).
$\rho$ also lifts to an action $\tilde{\rho}$ on $\tilde{\gg}$,
which keeps the flags constant. Choose a hermitian inner product on
$\gg$, take the function
\[
\psi(y) = \textstyle\half ||y||^2,
\]
and pull it back to a function $\tilde{\psi}$ on $\tilde{\gg}$. By
homogeneity, $\tilde{\chi}^{-1}(0)$ intersects all the level sets
$\tilde{\psi}^{-1}(c)$, $c>0$, transversally. It follows that there
is a small ball $B \subset \hh$ around the origin, such that for all
$\tilde{t} \in B$, $\tilde{\chi}^{-1}(\tilde{t})$ intersects
$\tilde\psi^{-1}(1)$ transversally. Using the $\C^*$-action to
rescale things, one sees that $\tilde{\chi}^{-1}(\tilde{t})$
intersects $\tilde\psi^{-1}(c)$ transversally for all $\tilde{t} \in
B$ and $c \geq 1$. An obvious argument with the gradient flow of
$\tilde\psi$ on $\tilde{\chi}^{-1}(\tilde{t})$ shows that one can
write that manifold as a union of a compact piece with boundary,
which is $\tilde{\chi}^{-1}(\tilde{t}) \cap \tilde\psi^{-1}([0;1])$,
and an infinite cone over that boundary. This is true for all
$\tilde{t} \in B$, which means that $\tilde{\chi}^{-1}(B) \cap
\tilde\psi^{-1}([0;1])$ is a differentiable fibre bundle with
compact fibres, and that the whole of $\tilde{\chi}^{-1}(B)$ is
obtained from it by attaching infinite cones to the fibre
boundaries. Finally, using the $\C^*$-action once more, the fibre
bundle structure can be extended from $B$ to the whole of $\hh$.
\qed

From this and Lemma \ref{th:simultaneous}(ii) it follows that $\chi:
\gg \rightarrow \hh/W$ itself, when restricted to $Conf_n^0(\C)$,
also becomes a differentiable fibre bundle.

Take a homogeneous slice $\Slice$ for a nilpotent $x$, and let
$\tilde\Slice$ be its preimage in $\tilde\gg$. Take $y \in \Slice$
and a preimage $\tilde{y} \in \tilde\Slice$. Because the adjoint
action lifts to $\tilde\gg$, the space $T_y(Gy)$ is contained in the
image of the differential $T_{\tilde{y}}\tilde\Slice \rightarrow
T_y\Slice$. In view of Lemma \ref{th:global-slice}(i), this has the
following consequences: the projection $\tilde\gg \rightarrow \gg$
is transverse to $\Slice$; hence, $\tilde\Slice \subset \tilde{\gg}$
is a smooth submanifold and transverse to all $G$-orbits; finally,
$\tilde\chi|\tilde\Slice: \tilde\Slice \rightarrow \hh$ is a
submersion.

\begin{Lemma} \label{th:fibre-bundle-2}
$\tilde\chi|\tilde\Slice: \tilde\Slice \rightarrow \hh$ is naturally
a differentiable fibre bundle.
\end{Lemma}

\begin{Proof}
$\lambda$ lifts to a $\C^*$-action $\tilde\lambda$ on
$\tilde\Slice$, which contracts that space to the compact subset
lying over the point $x \in \Slice$. On the base spaces, the
corresponding $\C^*$-action on $\hh/W$ obviously lifts to $\hh$, and
contracts that space to the origin. Let $\xi: \Slice \rightarrow \R$
be an exhausting function which is $\lambda$-homogeneous,
$\xi(\lambda_r(y)) = r^{2\alpha} \xi(y)$ for some $\alpha>0$. After
pulling this back to a function $\tilde{\xi}$ on $\tilde\Slice$, one
finds that $\tilde{\xi}^{-1}(c)$ intersects $\tilde{\chi}^{-1}(0)
\cap \tilde\Slice$ transversally for all $c>0$. The rest is as
before.
\end{Proof}

As in Lemma \ref{th:simultaneous}, the map
$\tilde\chi^{-1}(\tilde{t}) \cap \tilde\Slice \iso \chi^{-1}(t) \cap
\Slice$ is an isomorphism for all $t \in \Conf_n^0(\C)$. Hence, the
restriction of $\chi|\Slice: \Slice \rightarrow \hh/W$ to
$\Conf_n^0(\C)$ is again a differentiable fibre bundle.

We will also need a variation on the idea of simultaneous
resolution, involving partial flags. Fix some $k<n$. Let
$\gg^{\mult}$ be the space of pairs $(x,E)$, where $x \in \gg$ and
$E \subset \C^n$ is a $k$-dimensional subspace, such that $x|E$ is
some multiple of the identity map. $\gg^{\mult}$ is a smooth
manifold (in fact a bundle over the Grassmannian with Lie algebra
fibres). Correspondingly, let $\hh^{\mult} \subset \hh$ be the
subspace of diagonal matrices whose first $k$ entries coincide, and
$W^{\mult} \subset W$ the subgroup of permutations which leave the
first $k$ entries fixed. The quotient is $\hh^{\mult}/W^{\mult} \iso
\C \times \C^{n-k-1}/S_{n-k} \iso \C^{n-k}$. There is a natural
holomorphic map $\chi^{\mult}: \gg^{\mult} \rightarrow
\hh^{\mult}/W^{\mult}$, where the first entry of
$\chi^{\mult}(x,E) \in \C\times\C^{n-k-1}/S_{n-k}$ is the (single,
common) eigenvalue of $x|E$. Now let
$\tilde{\gg}^{\mult}$ be the space of pairs $(x,F)$, where $x \in
\gg$ and $F = \{0 = F^0 \subset F^k \subset F^{k+1} \subset \cdots
\subset F^n = \C^n\}$ a partial flag, satisfying $x(F^i) \subset
F^i$ for all $i$, and $(x,F^k) \in \gg^{\mult}$. This fits into a
commutative diagram, where as usual the left $\downarrow$ is a
holomorphic submersion:
\begin{equation} \label{eq:mult}
\begin{CD}
 {\tilde{\gg}^{\mult}} @>>> \gg^{\mult} \\
 @V{\tilde{\chi}^{\mult}}VV @V{\chi^{\mult}}VV \\
 {\hh}^{\mult} @>>> {\hh^{\mult}/W^{\mult}.}
\end{CD}
\end{equation}
There are obvious $G$-actions on $\gg^{\mult}$ and
$\tilde{\gg}^{\mult}$, corresponding to the adjoint action on $\gg$.
Therefore, the discussion preceding Lemma \ref{th:fibre-bundle-2}
carries over to the present situation, which means the following.
Let $\Slice \subset \gg$ be a homogeneous slice at some nilpotent
element $x$, $\Slice^{\mult} \subset \gg^{\mult}$ the subspace of
those $(y,E)$ such that $y \in \Slice$, and $\tilde\Slice^{\mult}$
the corresponding subspace of $\tilde\gg^{\mult}$. By restricting
\eqref{eq:mult} one gets a diagram
\begin{equation} \label{eq:slice-mult}
\begin{CD}
 {\tilde{\Slice}^{\mult}} @>>> \Slice^{\mult} \\
 @V{\tilde{\chi}^{\mult}}VV @V{\chi^{\mult}}VV \\
 {\hh}^{\mult} @>>> {\hh^{\mult}/W^{\mult}.}
\end{CD}
\end{equation}
and the left $\downarrow$ in this is again a holomorphic submersion.
Moreover, the $\C^*$-action $\lambda$ on $\Slice$ also lifts to
$\Slice^{\mult}$ and $\tilde{\Slice}^{\mult}$, which means that one
can argue as in Lemma \ref{th:fibre-bundle-2}, showing that the
submersion is in fact a differentiable fibre bundle.

\begin{Remarks}
For a general semisimple Lie algebra $\gg$, one defines $\tilde\gg$
to be the space of pairs $(x,{\mathfrak b})$, where ${\mathfrak b}
\subset \gg$ is a Borel subalgebra containing $x$. This comes with a
map $\tilde\chi: \tilde\gg \rightarrow \hh$, which is a simultaneous
resolution of the adjoint quotient map $\chi$ \cite[Section
3.3]{slodowy82} (the corresponding result for algebraic groups
already appears in \cite{brieskorn70b}). With that at hand, Lemma
\ref{th:fibre-bundle} generalizes easily (see \cite{rossmann95} for
another version of the same argument, and a more detailed study of
the monodromy of the resulting fibration over $\hh^{\reg}/W$).

The use of simultaneous resolution for Jacobson-Morozov slices also
goes back to Brieskorn and Slodowy \cite{brieskorn70b,slodowy}, who
looked at subregular nilpotent elements in simply-laced $\gg$. For
those slices, one recovers the simultaneous resolutions of ADE
singularities that had been previously discovered by Brieskorn by
more elementary means; Example \ref{ex:sl2-cont} is the simplest
case.
\end{Remarks}

\section{$(m,m)$-type nilpotent slices\label{sec:slice}}

We now focus on the particular slices $\Slice_m \subset \sl_{2m}$
relevant to our main construction. Running the basic idea of
\eqref{eq:basic} and Lemma \ref{th:global-slice} in reverse, one
finds that the geometry of $\chi|\Slice_m$ is modelled on that of
the whole adjoint quotient map $\chi$. In particular, in a process
which goes back to Examples \ref{ex:sl2} and \ref{ex:sl3}, we will
see $(A_1)$ and $(A_2)$ singularities appearing. The other main
point is that $\Slice_{m-1}$ is embedded into $\Slice_m$ in a
natural way (Lemma \ref{th:induction}); this will form the basis for
several inductive arguments later on.

\subsection{Definition and first properties}

From now on $n = 2m$, which means that we will work with $\gg =
\sl_{2m}$. Let $x \in \gg$ be a nilpotent element with two Jordan
blocks of size $m$. We find it convenient to think of $\C^{2m} =
\C^2 \oplus \cdots \oplus \C^2$, and then to write $n^+$ as in
\eqref{eq:standard-jordan}, but where the scalar entries $0,1$ are
replaced by the corresponding 2x2 matrices (the zero matrix and the
identity matrix in ${\mathfrak gl}_2$) mapping the $\C^2$ summands
to each other. The equations \eqref{eq:standard-jm} and
\eqref{eq:standard-lambda}, when re-interpreted in the same sense,
then describe a JM triple $(n^+ = x,n^-,h)$ and the associated
$\C^*$-action $\lambda$ on $\gg$. Consider the affine space
$\Slice_m = n^+ + \ker(z \mapsto zn^-)$, which consists of matrices
of the form
\begin{equation} \label{eq:y-matrix}
 y = \begin{pmatrix}
 y_{11} & 1 &&& \\
 y_{21} && 1 && \\
 \dots &&& \dots & \\
 y_{m-1,1} &&&& 1 \\
 y_{m1} &&&& 0
 \end{pmatrix}
\end{equation}
with $y_{11} \in \sl_2$, and $y_{i1} \in \mathfrak{gl}_2$ for $i>1$.
In parallel with Example \ref{ex:regular1}, we have:

\begin{Lemma}
$\Slice_m$ is a homogeneous slice for the adjoint action at $x$.
\end{Lemma}

\begin{Proof}
Suppose that $z \in \gg$ is such that only the first two columns of
$[n^+,z]$ are nonzero. Inspection of $[n^+,z]$ shows that $z$ must
be upper triangular when written in our usual block form, and from
that one sees that $[n^+,z] = 0$. This shows that the tangent space
to the adjoint orbit intersects the tangent space to our slice
trivially. By looking at the eigenvalues of $ad(h)$, one sees that
the $\sl_2$-module $\gg$ splits into a direct sum of irreducible
representations of dimensions $1,3,\dots,2m-1$. More precisely,
there are three trivial summands of dimension $1$, and four summands
of all the other dimensions. Thinking of the Jacobson-Morozov
procedure, this means that a transverse slice must have dimension $3
+ 4(m-1)$. This shows that the tangent space to our slice is indeed
complementary to the orbit directions. $\lambda$-invariance is
obvious from \eqref{eq:standard-lambda}.
\end{Proof}

\begin{Lemma} \label{th:eigenspace}
For any $y \in \Slice_m$ and any $\mu \in \C$, projection to the
first two coordinates yields an injective map $ker(\mu \cdot 1 - y)
\rightarrow \C^2$. In particular, that eigenspace is at most
two-dimensional.
\end{Lemma}

\begin{Proof}
Suppose that the contrary is true, which means that $ker(\mu \cdot 1
- y)$ has nonzero intersection with $\{0\}^2 \times \C^{2m-2}$.
Using the $\C^*$-action, one sees that the same holds for
$ker(r^2\mu \cdot 1 - \lambda_r(y))$, and as $r \rightarrow 0$, one
obtains a nonzero element in $ker(n^+) \cap (\{0\}^2 \times
\C^{2m-2})$, which is a contradiction.
\end{Proof}

\begin{Lemma} \label{th:induction}
The subspace of those $y \in \Slice_m$ such that $ker(y)$ is
two-dimensional can be canonically identified with $\Slice_{m-1}$.
This identification is compatible with the adjoint quotient map: if
$y$ has eigenvalues $(0,0,\mu_3,\dots,\mu_{2m})$, then the
corresponding element $\bar{y} \in \Slice_{m-1}$ has eigenvalues
$(\mu_3,\dots,\mu_{2m})$.
\end{Lemma}

\begin{Proof}
In the previous proof, we saw that a vector in $ker(y)$ is uniquely
determined by its first two entries. For there to be two linearly
independent such vectors, it is necessary and sufficient that
$y_{m1} = 0$. The identification of this subspace with
$\Slice_{m-1}$ is the straightforward one. One can see it as
restriction of linear maps to the subspace $\C^{2m-2} \times \{0\}^2
\supset im(y)$, and then the second statement becomes obvious.
\end{Proof}

\begin{Remark}
By Lemma \ref{th:invariant-slice}, $\Slice_m$ is orbit-preservingly
isomorphic to the Jacobson-Morozov slice at $n^+$, and
correspondingly for $\Slice_{m-1}$. As a consequence, Lemma
\ref{th:induction} also holds for Jacobson-Morozov slices, but the
isomorphism obtained in this way is no longer quite canonical (nor
as simple as before), and we have not found a more direct
construction. This is what makes $\Slice_m$ more convenient for our
purpose.
\end{Remark}

\subsection{Two eigenvalues coincide\label{subsec:2-coincide}}

Let $t \in \hh/W$ be a point corresponding to a collection of $2m$
pairwise different eigenvalues $(\mu_1,\dots,\mu_{2m})$. By Example
\ref{ex:different-ev} and Lemma \ref{th:global-slice}(ii), this is a
regular value of $\chi|\Slice_m$, so the fibre
$(\chi|\Slice_m)^{-1}(t)$ is a smooth complex manifold of dimension
$2m$.

Next, take the case where $\mu_1 = \mu_2$, with no other
coincidences between the eigenvalues. Then $\chi^{-1}(t)$ is the
union of two orbits: the regular orbit $\OO^{\reg}$ (of matrices
with an indecomposable Jordan block of size two for the eigenvalue
$\mu_1$), which is open and dense in $\chi^{-1}(t)$; and the
subregular orbit $\OO^{\sub}$ (of matrices having two independent
$\mu_1$-eigenvectors), which is closed. From Example
\ref{th:all-regular} we know that elements of $\OO^{\reg}$ are
regular points of $\chi$, whereas those of $\OO^{\sub}$ have
singularities of $(A_1)$ type in transverse direction to the orbit.
In view of Lemma \ref{th:global-slice}, the intersections
$\OO^{\reg} \cap \Slice_m$, $\OO^{\sub} \cap \Slice_m$ have the same
properties with respect to the map $\chi|\Slice_m$. We need a more
precise global version of the latter statement, in which one can see
an entire neighbourhood of $\OO^{\sub} \cap \Slice_m$.

At every point $y \in \OO^{\sub} \cap \Slice_m$, choose a subspace
of $T_y\Slice_m$ which is complementary to $T_y(\OO^{\sub} \cap
\Slice_m)$ and depends holomorphically on $y$. This splitting
problem has a positive solution because $\OO^{\sub} \cap \Slice_m$
is affine, so that the relevant $Ext^1$ obstruction group is zero.
Translate those subspaces by adding $y$, to obtain a family of
affine subspaces $\Slice_y'$ which form a tubular neighbourhood of
$\OO^{\sub} \cap \Slice_m$ inside $\Slice_m$. Because $\Slice_m$
intersects $\OO^{\sub}$ transversally, each $\Slice_y'$ is also a
local slice at $y$ for the adjoint action on $\gg$. On the other
hand, $y$ is semisimple, so we are precisely in the situation
discussed in Example \ref{th:2-block}, which means that we can
construct another local slice at $y$ by setting $\Slice_y = y +
(\hat\gg \oplus {\mathfrak z})$. Recall that $\hat\gg = \sl(E_y)$,
where $E_y$ is the $\mu_1$-eigenspace of $y$. By Lemma
\ref{th:eigenspace}, projection yields preferred isomorphisms $E_y
\iso \C^2$ and hence $\hat\gg \iso \sl_2$. The remaining part
${\mathfrak z}$ can moreover be identified with $\C^{2m-2}$. By Lemma
\ref{th:stupid}(iii), we can find a local isomorphism $\Slice_y \iso
\Slice_y'$ which moves points only inside their adjoint orbits,
hence relates $\chi|\Slice_y$ to $\chi|\Slice_y'$. Strictly
speaking, this isomorphism \eqref{eq:stupid-id} requires a choice of
local submanifold $K_y \subset G$ complementary to $G_y$. One can
view this as another splitting problem, which can be solved in a way
that depends holomorphically on $y$ for cohomological reasons.
Alternatively, an elementary argument starting from $\gg =
[\gg,y]\oplus \gg_y$ shows that in this case, $K_y$
can be explicitly taken to be $exp([\gg,y])$. In either way, the outcome is
that we get a family of isomorphisms (defined locally near $y$)
$\Slice_y' \iso y + (\sl_2 \oplus \C^{2m-2})$. Moreover, with
respect to that trivialization, the adjoint quotient map becomes
\eqref{eq:fibered-a1-trivial} on each slice. In particular, we
obtain:

\begin{Lemma} \label{th:2-coincide}
Let $D \subset \hh/W$ be a disc corresponding to eigenvalues
$(\mu_1-\sqrt{\epsilon},\mu_2+\sqrt{\epsilon},\mu_3,\dots,\mu_{2m})$
with $\epsilon$ small. Then there is a neighbourhood of $\OO^{\sub}
\cap \Slice_m$ inside $\chi^{-1}(D) \cap \Slice_m$, and an
isomorphism of that with a neighbourhood of $(\OO^{\sub} \cap
\Slice_m) \times \{0\}^3$ inside $(\OO^{\sub} \cap \Slice_m) \times
\C^3$. This isomorphism fits into a commutative diagram
\[
\begin{CD}
 \chi^{-1}(D) \cap \Slice_m @>{\text{local $\iso$ defined near $\OO^{\sub} \cap
\Slice_m$}}>> (\OO^{\sub} \cap \Slice_m) \times \C^3 \\
 @V{\chi}VV @V{a^2+b^2+c^2}VV \\
 D @>{\quad\qquad\qquad\qquad\qquad\quad}>> \C
\end{CD}
\]
where $a,b,c$ are coordinates on $\C^3$. \qed
\end{Lemma}

This means that the function $\chi|\chi^{-1}(D) \cap \Slice_m$ is
nondegenerate in transverse direction to the critical submanifold
$\OO^{\sub} \cap \Slice_m$. This is the analogue of the Morse-Bott
condition in real topology, and we will refer to it by saying that
the critical submanifold is of fibered $(A_1)$ type (note that by
using the explicit local slices given by Lemma \ref{th:make-slice},
we have reached this conclusion without appealing to any Morse
Lemma-type arguments). Our case has the additional feature that the
normal data along the critical submanifold, consisting of the normal
bundle and the nondegenerate quadratic form on its fibres, are
trivial. It is instructive to compare it to the behaviour of the
adjoint quotient map on $\chi^{-1}(D)$ (without the slice):
$\OO^{\sub}$ is still a critical submanifold of fibered $(A_1)$
type, but the normal data are no longer trivial.

As a variation on this theme, one can consider the case where $2r$
eigenvalues come together in pairs, as in Example
\ref{th:2-2-block}. The relevant fibre $\chi^{-1}(t)$ consists of
$2^r$ orbits, since the restriction of $x \in \chi^{-1}(t)$ to its
generalized $\mu_{2j-1}$-eigenspace may be either semisimple or not,
for each $j = 1,\dots,r$. The smallest orbit
$\OO^{\scriptscriptstyle \mathrm{min}}$, which consists of those $x$
that are actually semisimple, is closed. A straightforward
adaptation of the previous argument yields the following description
of the local structure near $\OO^{\scriptscriptstyle \mathrm{min}}
\cap \Slice_m$:

\begin{Lemma} \label{th:2k-coincide}
Let $P \subset \hh/W$ be a $k$-dimensional polydisc corresponding
to eigenvalues
$(\mu_1-\sqrt{\epsilon_1},\mu_2+\sqrt{\epsilon_1},\dots,
\mu_{2k-1}-\sqrt{\epsilon_k},\mu_{2k}+\sqrt{\epsilon_k},
\mu_{2k+1},\dots,\mu_{2m}$) with the $\epsilon$'s small. Then there
is a neighbourhood of $\OO^{\scriptscriptstyle \mathrm{min}} \cap
\Slice_m$ inside $\chi^{-1}(P) \cap \Slice_m$, and an isomorphism of
that with a neighbourhood of $(\OO^{\scriptscriptstyle \mathrm{min}}
\cap \Slice_m) \times \{0\}^{3k}$ inside $(\OO^{\scriptscriptstyle
\mathrm{min}} \cap \Slice_m) \times \C^{3k}$. The isomorphism fits
into a commutative diagram
\[
\begin{CD}
 \chi^{-1}(P) \cap \Slice_m @>{\text{local $\iso$ defined near
 $\OO^{\scriptscriptstyle \mathrm{min}} \cap \Slice_m$}} >>
 (\OO^{\scriptscriptstyle \mathrm{min}} \cap \Slice_m) \times \C^{3k}
\\
 @V{\chi}VV @V{(a_1^2+b_1^2+c_1^2,\dots,a_k^k+b_k^2+c_k^2)}VV \\
 P @>{\hspace{10em}}>> \C^r
\end{CD}
\]
where $a_j,b_j,c_j$ are coordinates on $\C^{3k}$. \qed
\end{Lemma}

\subsection{Three eigenvalues coincide\label{subsec:3-coincide}}

Now take a point $t \in \hh/W$ which corresponds to a set of
eigenvalues $(\mu_1, \dots,\mu_{2m})$ of which the first three
coincide, $\mu_1 = \mu_2 = \mu_3$, and which are otherwise pairwise
distinct. The adjoint fibre $\chi^{-1}(t)$ contains three orbits:
the regular orbit $\OO^{\reg}$ has an indecomposable Jordan block of
size 3 for the eigenvalue $\mu_1$; the subregular orbit $\OO^{\sub}$
has two Jordan blocks of sizes $1,2$; and the minimal orbit consists
of matrices with three independent $\mu_1$-eigenvectors. However,
that last orbit does not intersect $\Slice_m$, due to Lemma
\ref{th:eigenspace}. As a consequence, $\OO^{\sub} \cap \Slice_m$ is
closed in $\Slice_m$, which makes the situation fairly similar to
the one we looked at before.

Take $y \in \OO^{\sub} \cap \Slice_m$, and let $E_y$ be the
$\mu_1$-eigenspace of its semisimple part $y_s$. As a first step, we
want to choose a JM triple for $y_n$ inside $\sl(E_y)$. As explained
in Example \ref{ex:sl3-again}, such triples correspond bijectively
to splittings of the flag
\begin{equation} \label{eq:flag}
 0 \subset F^1_y = (\mu_1 \cdot 1 - y)(E_y) \subset
 F^2_y = ker(\mu_1 \cdot 1 - y) \subset E_y.
\end{equation}
Because $\OO^{\sub} \cap \Slice_m$ is affine, one can use the
vanishing of $Ext^1$ to find such splittings which vary
holomorphically with $y$, and thence holomorphically varying JM
slices $\hat\Slice_{y_n} \subset \sl(E_y)$. From these, Lemma
\ref{th:make-slice} produces slices $\Slice_y \subset \gg$ for $y$
itself, without any further choices. To find a better global
picture, we appeal again to Kostant's uniqueness theorem. This says
that we can find an isomorphism
\begin{equation} \label{eq:three-d}
E_y \iso \C^3
\end{equation}
such that the induced map $\sl(E_y) \iso \sl_3$ identifies our JM
triple with the one explicitly given in Example \ref{ex:sl3-again}.
The isomorphism is unique up to the action of the subgroup $\C^*
\times \C^* \subset GL_3$ consisting of diagonal matrices
$diag(\zeta,\tau,\zeta)$: namely, the isomorphism of Lie algebras is
unique up to the subgroup of diagonal matrices in \eqref{eq:g}, and
lifting that to an isomorphism of vector spaces adds the central
$\C^* \subset GL_3$. By construction, \eqref{eq:three-d} takes
$F_y^2$ to $\C^2 \times \{0\}$. Another isomorphism $F_y^2 \iso
\C^2$ is given by Lemma \ref{th:eigenspace}, and we can constrain
\eqref{eq:three-d} by asking that the two maps have equal
determinant. This reduces the ambiguity to a single $\C^*$ factor
$diag(\zeta,\zeta^{-1},\zeta)$, which means that \eqref{eq:three-d}
is completely determined by the induced identification $F_y^1
\rightarrow \C \times \{0\}^2$. The adjoint action of the $\C^*$ on
the model slice, in the coordinates \eqref{eq:a2-map}, is
\begin{equation} \label{eq:slice-action}
(a,b,c,d) \longmapsto (a,\zeta^{-2}b,\zeta^2 c,d).
\end{equation}
With that in mind, the argument from Lemma \ref{th:2-coincide}
carries over up to some easy modifications, with the following
outcome:

\begin{Lemma} \label{th:3-coincide}
Let $\FF \rightarrow \OO^{\sub} \cap \Slice_m$ be the line bundle
whose fibres are the spaces $F_y^1$ from \eqref{eq:flag}. Consider
the associated vector bundle
\begin{equation} \label{eq:associated-c4}
(\FF \setminus 0) \times_{\C^*} \C^4 = \C \oplus \FF^{-2} \oplus
\FF^2 \oplus \C
\end{equation}
with respect to the action from \eqref{eq:slice-action}; here $0$
denotes the zero-section and $\C$ the trivial line bundle over
$\OO^{\sub} \cap \Slice_m$. Let $P \hookrightarrow \hh/W$ be a small
bidisc parametrized by $(d,z)$, corresponding to the set of
eigenvalues
\[
(\mu_1+\{\text{all solutions of $\lambda^3 - d \lambda  + z = 0
$}\},\mu_4,\dots,\mu_{2m}).
\]
There is a neighbourhood of $\OO^{\sub} \cap \Slice_m$ inside
$\chi^{-1}(P) \cap \Slice_m$, and an isomorphism of that with a
neighbourhood of the zero-section inside $(\FF \setminus 0)
\times_{\C^*} \C^4$, which fits into a commutative diagram
\[
\begin{CD}
 \chi^{-1}(P) \cap \Slice_m
 @>{\text{local $\iso$ defined near $\OO^{\sub} \cap
 \Slice_m$}}>> {(\FF \setminus 0) \times_{\C^*} \C^4} \\
 @V{\chi}VV @V{p}VV \\
 P @>{\quad\qquad\quad(d,z)\quad\qquad\quad}>> \C^2
\end{CD}
\]
where $p$ is given by \eqref{eq:a2-map} on each $\C^4$ fibre.
\end{Lemma}

In the statement we have used $d$ both as a coordinate on $\C^4$ and
on $P$, but that should not be troublesome since the first component
of $p$ maps one identically to the other. Note also that the second
component $a^3 - ad + bc$ makes sense as a holomorphic function on
\eqref{eq:associated-c4} because $b$ and $c$ are sections of inverse
line bundles.

\subsection{A partial Grothendieck resolution\label{subsec:partial}}

Let's apply the construction from \eqref{eq:slice-mult}, with $k = 2$
and $n = 2m>2$, to our slices $\Slice_m$. Denote the resulting
spaces by $\Slice_m^{\mult}$, $\tilde{\Slice}_m^{\mult}$. Consider
the open subset $\hh^{\mult,\reg} \subset \hh^{\mult}$ of those
$(\mu_1,\mu_2 = \mu_1,\mu_3,\dots,\mu_{2m})$ where the
$(\mu_3,\dots,\mu_{2m})$ are pairwise disjoint among themselves (but
any of them may agree with $\mu_1 = \mu_2$). It is clear from the
definitions that if $(y,E) \in \Slice_m^{\mult}$ is such that
$\chi^{\mult}(y,E) \in \hh^{\mult,\reg}/W^{\mult}$, then a choice of
preimage $(y,F) \in \tilde{\Slice}_m^{\mult}$ is the same as an
ordering of the eigenvalues $\mu_3,\dots,\mu_{2m}$. In other words,
the restriction of \eqref{eq:slice-mult} to
$\hh^{\mult,\reg}/W^{\mult}$ is a pullback diagram. As a
consequence, the restriction of
\begin{equation} \label{eq:mult-fibr}
\chi^{\mult}|\Slice_m^{\mult}: \Slice_m^{\mult} \longrightarrow
\hh^{\mult}/W^{\mult}
\end{equation}
to $\hh^{\mult,\reg}/W^{\mult}$ becomes a fibre bundle. Fix some
$(\mu_1,\dots,\mu_{2m}) \in \hh^{\mult,\reg}$, let $t^{\mult}$ be
its image in $\hh^{\mult,\reg}/W^{\mult}$, and $t$ its image in
$\hh/W$. By definition, the fibre of \eqref{eq:mult-fibr} over
$t^{\mult}$ is the set of pairs $(y,E)$, where $y \in \Slice_m \cap
\chi^{-1}(t)$ is such that $y|E$ is a multiple of the identity. That
multiple must be the unique multiple eigenvalue $\mu_1$, and then
necessarily $E = ker(\mu_1 \cdot 1 - y)$ by Lemma
\ref{th:eigenspace}. In other words, the fibre of
\eqref{eq:mult-fibr} can be identified with the subspace of
$\chi^{-1}(t) \cap \Slice_m$ consisting of those $y$ which have two
independent $\mu_1$-eigenvalues. From the results in Sections
\ref{subsec:2-coincide} and \ref{subsec:3-coincide}, we see that
this subspace is just the set of singular points of $\chi^{-1}(t)
\cap \Slice_m$, which we will denote by $\CC_{m,t}$ from now on. The
upshot of this discussion is that the $\CC_{m,t}$ form a
differentiable fibre bundle
\begin{equation} \label{eq:critical-fibration}
\CC_m \longrightarrow \hh^{\mult,\reg}/W^{\mult}.
\end{equation}

\section{Parallel transport, Floer cohomology\label{sec:symplectic}}

This section deals with the necessary K{\"a}hler and symplectic
geometry, largely in the context of Stein fibre bundles. The main
objectives are the definition of relative vanishing cycles, and two
technical statements about Floer cohomology groups in situations
where the Lagrangian submanifolds concerned are constructed as such
cycles (Lemma \ref{th:kunneth} and Lemma \ref{th:thom}). Our main
technical trick involves deforming K{\"a}hler forms to make them
agree with the standard forms on certain subsets where we have
preferred holomorphic coordinates. One has to worry whether the
resulting families of Lagrangian submanifolds remain inside those
subsets, and addressing that requires some technical estimates of
parallel transport vector fields.

\subsection{Parallel transport\label{subsec:parallel}}

Let $\pi: Y \rightarrow T$ be a holomorphic map between complex
manifolds, which is a submersion with fibres $Y_t$. Suppose that $Y$
carries a K{\"a}hler metric, and equip the fibres with the induced
metrics, which in particular makes them into symplectic manifolds.
Take a path $\gamma: [0;1] \rightarrow T$ on the base. The parallel
transport vector field $H_\gamma$ is a vector field on the pullback
$\gamma^*Y \rightarrow [0;1]$: it consists of the unique sections of
$TY|Y_{\gamma(s)}$ which project to $\dot\gamma(s)$, and which are
orthogonal to the tangent space along the fibres. In the case where
$T = \C$, one can write explicitly
\begin{equation} \label{eq:horizontal}
 H_\gamma = \frac{\nabla \pi}{||\nabla \pi||^2} \dot\gamma(s).
\end{equation}
When $\pi$ is proper, integrating $H_\gamma$ yields a symplectic
isomorphism between fibres, the parallel transport map $h_\gamma:
Y_{\gamma(0)} \rightarrow Y_{\gamma(1)}$. If properness fails, the
integral lines may not exist for all times. In some cases, explicit
estimates of $H_\gamma$ may allow one to show that $h_\gamma$ is
still defined everywhere, or at least on a subset which is
sufficiently large to contain the geometric objects (Lagrangian
submanifolds, in our applications) that one wants to apply parallel
transport to. Alternatively, one can try to modify the vector field,
so as to make the domain of definition larger. We will now explain a
basic argument of the second kind, for the case when the fibres
$Y_t$ are Stein manifolds with finite topology (to simplify the
description, we impose slightly sharper technical conditions than
strictly necessary).

Suppose that there is a proper bounded below function $\psi: Y
\rightarrow \R$ such that
\begin{itemize} \itemsep0em
\item $-dd^c \psi$ is the given K{\"a}hler form on $Y$. \item
Outside a compact subset of $Y$, $||\nabla \psi||^2 \leq \rho \psi$
for some $\rho>0$. \item The fibrewise critical set of $\psi$,
consisting of those points where $d\psi|ker(D\pi)$ is zero, maps
properly to $T$.
\end{itemize}
Fix $t \in T$. The function $\psi_t = \psi|Y_t$ is proper and has a
compact set of critical points. Let $Z_t$ be its gradient vector
field. We have $Z_t . \psi_t = ||\nabla \psi_t||^2 \leq ||\nabla
\psi||^2 \leq \rho \psi_t$ outside a compact subset, which ensures
that the flow of $Z_t$ is defined for all times. In symplectic
geometry terms, this is a Liouville vector field, so $Y_t$ is a
manifold with an infinite convex contact-type cone. As before, let
$\gamma$ be a path in $T$. Because of the properness condition on
the fibrewise critical point set, we can choose $c>0$ in such a way
that the critical values of $\psi$ on the fibres $Y_{\gamma(s)}$ all
lie in $[0;c)$. Using that and the properness of $\psi$, one can
find a $\sigma>0$ such that the modified parallel transport vector
field
\[
 \bar{H}_\gamma = H_\gamma - \sigma Z_{\gamma(s)}
\]
satisfies $d\psi(\bar{H}_\gamma) < 0$ at all points $y \in
Y_{\gamma(s)}$ with $\psi(y) = c$. The integral lines of that vector
field necessarily stay within $\psi^{-1}([0;c])$, hence give rise to
a well-defined map $Y_{\gamma(0)} \cap \psi^{-1}([0;c]) \rightarrow
Y_{\gamma(1)}$. This is only conformally symplectic, but one can
repair that by composing with the time $\sigma$ map of the Liouville
flow on $Y_{\gamma(1)}$. The result is a symplectic embedding
\[
 h_\gamma^{\resc}: Y_{\gamma(0)} \cap \psi^{-1}([0;c]) \longrightarrow
 Y_{\gamma(1)}
\]
called rescaled symplectic parallel transport (to make the
distinction clear, we will sometimes refer to the maps $h_\gamma$
obtained by simply integrating horizontal vector fields as naive
parallel transport). It is independent of the choice of $\sigma$ up
to isotopy within the class of symplectic embeddings. One can take
$c$ arbitrarily large, and thereby define $h^{\resc}_\gamma$ on
arbitrarily big compact subsets of $Y_{\gamma(0)}$. Passing from
some value of $c$ to a larger one yields a map whose restriction to
the smaller domain is isotopic to the previous one. As a
consequence, the image $h^{\resc}_\gamma(L)$ of a closed Lagrangian
submanifold $L \subset Y_{\gamma(0)}$ is well-defined up to
Lagrangian isotopy, which is the most important fact for our
purpose.

\begin{Remark} \label{th:transport}
Under the same assumptions, one can construct rescaled parallel
transport maps which are defined on the whole fibre and are
symplectic isomorphisms, modelled at infinity on contactomorphisms
(this means that they commute with the Liouville flows outside
compact subsets). The first step is Gray's stability theorem, which
shows that the hypersurfaces $\psi^{-1}(c) \cap Y_{\gamma(s)}
\subset Y_{\gamma(s)}$ for various $s$ are isomorphic as contact
manifolds. Hence, the cone-like ends of the fibres are
symplectically isomorphic. One modifies the parallel transport
vector field to be compatible with these isomorphisms outside a
compact subset; for details see \cite[Section 6]{khovanov-seidel98}.
\end{Remark}

\subsection{Relative vanishing cycles\label{subsec:relative-vanishing}}

We now consider the local K{\"a}hler geometry around a fibered
$(A_1)$ type critical set. Since our main intended application is
provided by Lemma \ref{th:2-coincide}, we will concentrate on the
case where the normal data to the critical point set are trivial
(for a partial loosening of this restriction, see Remark
\ref{th:twisted-bundle}). Therefore, take any complex manifold $X$
and consider
\begin{equation} \label{eq:quadratic}
 \pi: Y = X \times \C^3 \longrightarrow \C, \quad
 \pi(x,a,b,c) = a^2 + b^2 + c^2.
\end{equation}
The critical point set is $Crit(\pi) = \{a = b = c = 0\}$, hence can
be identified with $X$. We equip $Y$ with any K{\"a}hler metric, the
fibres $Y_t$ with the induced metrics, and $X$ with the restriction
of the metric to $Crit(\pi)$. The symplectic form is denoted by
$\Omega \in \Omega^2(Y)$. Since the second derivative of $\pi$ in
transverse direction to $X$ is nondegenerate, the real part
$re(\pi)$ is a Morse-Bott function. Define its stable manifold $W
\subset Y$ to be the set of points $y$ such that the flow line of
$-\nabla re(\pi)$ starting at $y$ exists for all times $s \geq 0$,
and converges to a critical point in the limit $s \rightarrow
\infty$. Obviously $X$ itself is contained in $W$. Moreover, since
the negative gradient flow of $re(\pi)$ is also the Hamiltonian
vector field of $im(\pi)$, it leaves $im(\pi)$ invariant, hence $W$
lies inside $\pi^{-1}(\R^{\scriptscriptstyle \geq 0})$.

\begin{Lemma}
(i) $W \subset Y$ is a local real submanifold of codimension 3, and
its tangent space along $Crit(\pi)$ is $TX \times \R^3$. (ii) The
map $l: W \rightarrow X$ which assigns to a point its limit under
the negative gradient flow is a smooth submersion. (iii) $\Omega|W$
is equal to the pullback of $\Omega|X$ under $l$.
\end{Lemma}

The first two statements are standard Morse-Bott theory. One
possible approach, carried out in detail in \cite[Appendix
A]{austin-braam95}, goes roughly as follows. One first shows that
the convergence of gradient flow lines towards critical points
happens with exponential speed. Then, defining $W$ as a subspace of
the Banach manifold of all paths in $Y$ converging exponentially
towards a critical point, one finds that it is smooth by using the
implicit function theorem. The dimension can be computed from the
index of a suitably linearized problem, and the map $l$ is smooth by
construction. As for (iii), the gradient flow is symplectic, and if
we restrict it to $W$, then the limit of its derivative as $s
\rightarrow \infty$ gives $Dl$.

\begin{Lemma}
Let $K \subset X$ be a compact Lagrangian submanifold. Then for
sufficiently small $t>0$, $L_t = l^{-1}(K) \cap Y_t$ is a Lagrangian
submanifold of $Y_t$ diffeomorphic to $K \times S^2$.
\end{Lemma}

\begin{Proof}
$\pi|W: W \rightarrow \R$ is a function with a Morse-Bott type
nondegenerate minimum along $X \subset W$. The same holds if we
restrict it to the submanifold $l^{-1}(K)$ of points whose limit
lies in $K$. The Morse-Bott Lemma, together with the fact that the
normal bundle of $X$ in $W$ is trivial, imply that the sets $L_t$
for small $t$ are trivial $S^2$-bundles over $K$. The Lagrangian
property follows immediately from the previous Lemma.
\end{Proof}

We call $L_t$ the relative vanishing cycle associated to $K$. Of
course, by multiplying $\pi$ with some constant in $S^1$, one can
define stable manifolds which lie over other half-lines in $\C$, and
relative vanishing cycles $L_t \subset Y_t$ for all sufficiently
small $t \in \C^*$. There is also an equivalent formulation in terms
of parallel transport. Take the path $\gamma: [0;1] \rightarrow \C$,
$\gamma(s) = (1-s)t$, which runs straight into the critical value.
Then $L_t$ is the set of those $y \in Y_t$ such that the (naive)
parallel transport maps $h_{\gamma|[0;s]}$ are well-defined near $y$
for all $s<1$, and such that as $s \rightarrow 1$,
$h_{\gamma|[0;s]}(y)$ converges to a point of $K$. The two
definitions are equivalent essentially because for $t > 0$,
$H_{\gamma}$ and $-\nabla re(\pi)$ agree up to a positive scalar
factor, see \eqref{eq:horizontal}.

\begin{Remark} \label{th:twisted-bundle}
A variant of this geometry is where one has a holomorphic line
bundle $\LL \rightarrow X$ and looks at $Y = \C \oplus \LL^{-1}
\oplus \LL$ with the function $\pi: Y \rightarrow \C$, $\pi(a,b,c) =
a^2 + bc$ where $(a,b,c)$ are the fibre coordinates. The
construction of relative vanishing cycles $L_t \subset Y_t$ from $K
\subset X$ goes through as before, the only difference being that
topologically $L_t$ is a possibly nontrivial $S^2$-bundle over $K$,
in fact the projectivization $P(\LL \oplus \C)|K$.
\end{Remark}

To supplement the previous discussion, and (more importantly) as a
warmup exercise for \ref{subsec:a2} below, we will now explicitly
estimate $H_\gamma$ and thereby bound the position of the relative
vanishing cycles. Fix a relatively compact open subset $U \subset X$
and a ball $B \subset \C^3$ around the origin, and set $V = U \times
B \subset Y$.

\begin{Lemma} \label{th:bott-estimate}
There is a constant $\nu>0$ such that on $V$,
\[
 ||\nabla \pi||^2 \geq \nu^{-1} |\pi|.
\]
\end{Lemma}

\begin{Proof}
This would hold everywhere, with $\nu = 1/4$, if our metric was the
product of some metric on $X$ and the standard metric on $\C^3$.
Since $V$ is relatively compact, the statement is independent of the
choice of metric.
\end{Proof}

Now take a compact Lagrangian submanifold $K \subset U \subset X$.
Consider it as lying in the critical set of $\pi$, and denote by
$\delta>0$ its distance from $\partial V$ with respect to the given
metric. We want to show that the relative vanishing cycle $L_t$ is
well-defined and lies in $V$ for all
\begin{equation} \label{eq:bound-s}
0 < |t| < (1/100)\nu^{-1}\delta^2.
\end{equation}
Assume that $t>0$, and think in terms of parallel transport along
$\gamma(s) = s$. From the definition \eqref{eq:horizontal} and the
Lemma above, one sees that inside $V$, the horizontal vector field
on $Y_{\gamma(s)}$ is bounded by
\begin{equation} \label{eq:hgamma-estimate}
 ||H_{\gamma}|| \leq \nu^{1/2} s^{-1/2}.
\end{equation}
Suppose that we have a flow line of this vector field defined for $s
\in (0;t)$, and which converges to a point of $K$ as $s \rightarrow
0$. Supposing that $t$ satisfies \eqref{eq:bound-s}, then by
integrating \eqref{eq:hgamma-estimate} one finds that the whole flow
line lies at distance at most $2 \nu^{1/2} t^{1/2} < \delta/2$ from
$K$, hence it extends to $s = t$. With that in mind, the
well-definedness of the vanishing cycles and the fact that they lie
in $V$ is clear.

A similar estimate shows that if we take $t$ as in
\eqref{eq:bound-s} and consider the circle $\gamma_t: [0;2\pi]
\rightarrow \C^*$, $\gamma_t(s) = t \,\textit{exp}(is)$, then
parallel transport $h_{\gamma_t}(y)$ is well-defined and lies in $V$
for all $y \in L_t$. A priori we now have two Lagrangian
submanifolds in $Y_t$, namely the relative vanishing cycle $L_t$ and
its monodromy image $h_{\gamma_t}(L_t)$, however:

\begin{Lemma}
$L_t$ is Lagrangian isotopic to $h_{\gamma_t}(L_t)$ inside $V \cap
Y_t$.
\end{Lemma}

It may be helpful to first consider the case when the metric on $Y$
is the product of some K{\"a}hler metric on $X$ and the standard
metric on $\C^3$. Then the relative vanishing cycles are
\begin{equation} \label{eq:lzero}
 L_t = K \times \sqrt{t}\, S^2 \subset Y_t
\end{equation}
where $\sqrt{t} S^2 \subset \sqrt{t} \R^3 \subset \C^3$. The
monodromy is $id_X$ times the standard Picard-Lefschetz (Dehn twist)
monodromy, see \cite{seidel01} for an explicit computation, and the
Lemma is trivially true, since $h_{\gamma_t}(L_t) = L_t$. For
general metrics one argues as follows. The estimates made above show
that for all $\tau = \gamma_t(s)$ the relative vanishing cycle
$L_\tau$ is well-defined, and so is the (naive) parallel transport
along $\gamma_t|[s;2\pi]$ at least near $L_\tau$. Clearly
$h_{\gamma_t| [s;2\pi]}(L_\tau)$ is a family of Lagrangian
submanifolds connecting $h_{\gamma_t}(L_t)$ with $L_t$.

Guided by Lemma \ref{th:2k-coincide}, we also want to look at the
situation where the construction of vanishing cycles can be
iterated. Namely take $Y = X \times \C^{3k}$, with the function
\[
 \pi: Y \longrightarrow \C^k, \quad \pi(x,a_1,b_1,c_1,\dots) =
 (a_1^2+b_1^2+c_1^2,\dots,a_k^2+b_k^2+c_k^2),
\]
and a compact Lagrangian submanifold $K$ of $X = X \times \{0\}^{3k}
\subset Y$. One starts with the first component of $\pi$, restricted
suitably to a function
\[
\pi_1 : X \times \C^3 \times \{0\}^{3k-3} \longrightarrow \C.
\]
This yields a relative vanishing cycle $L_{t_1} \in \pi_1^{-1}(t_1)$
for small $t_1 \neq 0$. Fix some value of that parameter, and
consider the next component
\[
 \pi_2 : \pi^{-1}(\{t_1\} \times \C^{k-1}) \cap
 (X \times \C^6 \times \{0\}^{3k-6}) \longrightarrow \C.
\]
This has $\pi_1^{-1}(t_1)$ as its critical locus, and by writing
down things explicitly one sees that the relative vanishing cycle
construction can be applied to $L_{t_1}$ yielding an $L_{t_1,t_2}
\in \pi_2^{-1}(t_2)$. By repeating this one finally obtains an
iterated relative vanishing cycle, which is a Lagrangian submanifold
$L_{t_1,\dots,t_k} \subset Y_{t_1,\dots,t_k}$ diffeomorphic to $K
\times (S^2)^k$. A priori this may appear to work only for $0 <
|t_k| \ll |t_{k-1}| \ll |t_{k-2}| \ll \dots \ll |t_1|$, but an
inspection of the relevant parallel transport vector fields shows
that there is a uniform bound for all coordinates, meaning that
there is a $\sigma>0$ such that $L_{t_1,\dots,t_k}$ is defined
whenever $0<|t_j|<\sigma$ for all $j$. With this in mind, it makes
sense to state:

\begin{Lemma} \label{th:order}
If one changes the order in which the components of $\pi$ are used
to construct the iterated vanishing cycle, the outcome is the same
up to Lagrangian isotopy, at least as long as all the $|t_j|$ are
sufficiently small.
\end{Lemma}

\Proof Let $\Omega = \Omega^{(0)}$ be the given K{\"a}hler form on
$Y$. By restricting it to $X \times \{0\}^{3k}$ and taking the sum
of that and the standard form on $\C^{3k}$, one gets another
K{\"a}hler form $\Omega^{(1)}$. The statement of the Lemma would be
trivial for $\Omega^{(1)}$, since the corresponding iterated
vanishing cycles are simply
\begin{equation} \label{eq:squares}
 L^{(1)}_{t_1,\dots,t_k} = K \times \sqrt{t_1}S^2 \times \dots
 \times \sqrt{t_k}S^2.
\end{equation}
We will now use a Moser Lemma argument. Take the family
$\Omega^{(r)}$, $0 \leq r \leq 1$, of K{\"a}hler forms which
interpolate linearly between the two previously mentioned ones. By
integrating radially away from $X$, one can write $\Omega^{(1)} -
\Omega^{(0)} = d\Theta$ for some one-form $\Theta$ such that
$\Theta_y = 0$ for each $y = (x,0,\dots,0) \in X \times \{0\}^{3k}$.
Take a relatively compact open subset $V \subset Y$ which contains
$K$, and choose sufficiently small $t_1,\dots,t_k \neq 0$. Then the
iterated vanishing cycle $L^{(r)}_{t_1,\dots,t_k}$ is well-defined
for all $r$, and moreover, the Moser vector fields constructed from
$\Theta|Y_{t_1,\dots,t_k}$ integrate to give a family of symplectic
embeddings
\[
\phi_{t_1,\dots,t_k}^{(r)}: (V \cap
\pi^{-1}(t_1,\dots,t_k),\Omega^{(r)}) \longrightarrow
(\pi^{-1}(t_1,\dots,t_k),\Omega^{(0)}).
\]
From this one gets a Lagrangian isotopy from $L_{t_1,\dots,t_k}$ to
$\phi_{t_1,\dots,t_k}^{(1)}(L^{(1)}_{t_1,\dots,t_k})$. The same can
be done for the iterated vanishing cycles constructed using a
different ordering of the components of $\pi$, and since the
endpoints of the two isotopies are the same by \eqref{eq:squares},
the result follows. \qed

\subsection{Fibered $(A_2)$ singularities\label{subsec:a2}}

Basically the same strategy can be applied to the geometric
situation which appears in Lemma \ref{th:3-coincide}. However, since
that is somewhat more complicated, we prefer to first explain the
argument in the simplest example, which is just the map $\pi = p: Y
= \C^4 \rightarrow \C^2$ from \eqref{eq:a2-map}. We write $(d,z)$
for the coordinates on the base $\C^2$. The critical point set is
$Crit(\pi) = \{b = c = 0, \; d = 3a^2\}$, and its image is the cusp
curve $\Sigma = \{4 d^3 = 27 z^2\}$. We want to view $d$ as an
auxiliary parameter, so we consider the restrictions of $\pi$ to
$Y_d = \C^3 \times \{d\}$ as a family of functions $\pi_d: Y_d
\rightarrow \{d\} \times \C = \C$, writing $Y_{d,z} = \pi^{-1}(d,z)
= \pi_d^{-1}(z)$ for their fibres. The critical point set can then
be written as $Crit(\pi_d) = \{b = c = 0, \; a = \pm \sqrt{d/3}\}$.
For $d \neq 0$ this consists of two nondegenerate critical points,
denoted by $Crit(\pi_d)^{\pm}$, which project to the critical values
$\zeta_d^\pm = \pm \sqrt{4d^3/27}$. Later on, we will mostly
consider the case when $d > 0$, and then the convention is that
$\zeta_d^+ > 0$. For $d = 0$ the two critical points coalesce into a
single more degenerate one, which is of course exactly how
singularity theorists came to study $\pi$.

Equip $Y$ with some K{\"a}hler form $\Omega$. For any $d>0$ and $0 <
\epsilon \ll d$ there is a natural Lagrangian two-sphere
\begin{equation} \label{eq:l0}
L_{d,\epsilon} \subset Y_{d,\zeta_d^- + \epsilon},
\end{equation}
namely the vanishing cycle of $\pi_d: Y_d \rightarrow \C$ associated
to the critical value $\zeta_d^-$ (this is the classical vanishing
cycle construction, which is the special case $X = point$ of the
discussion in \ref{subsec:relative-vanishing} above). Take the path
$\gamma_{d,\epsilon}$ in $\C \setminus \{\zeta_d^\pm\}$ which runs
from $\zeta_d^- + \epsilon$ to $0$ along the real axis, then makes a
positive full circle around $\zeta_d^+$, and finally goes back to
its starting point along the real axis, see Figure
\ref{fig:gamma-path}. We want to look at the image of the vanishing
cycle by parallel transport, which is another Lagrangian two-sphere
\begin{equation} \label{eq:monodromy-image}
h_{\gamma_{d,\epsilon}}(L_{d,\epsilon}) \subset Y_{d,\zeta_d^- +
\epsilon}.
\end{equation}
To show that this is well-defined, it is necessary to control the
size of the parallel transport vector field. Fix a ball $B \subset
\C^4$ around the origin.
\begin{figure}[ht]
\begin{centering}
\begin{picture}(0,0)%
\includegraphics{gamma-path.pstex}%
\end{picture}%
\setlength{\unitlength}{3947sp}%
\begingroup\makeatletter\ifx\SetFigFont\undefined%
\gdef\SetFigFont#1#2#3#4#5{%
  \reset@font\fontsize{#1}{#2pt}%
  \fontfamily{#3}\fontseries{#4}\fontshape{#5}%
  \selectfont}%
\fi\endgroup%
\begin{picture}(4392,2692)(451,-2594)
\put(3451,-1411){\makebox(0,0)[lb]{\smash{\SetFigFont{12}{14.4}{\familydefault}{\mddefault}{\updefault}{$\zeta_d^+$}%
}}}
\put(3301,-2536){\makebox(0,0)[lb]{\smash{\SetFigFont{12}{14.4}{\familydefault}{\mddefault}{\updefault}{$\gamma_{d,\epsilon}$}%
}}}
\put(1051,-1411){\makebox(0,0)[lb]{\smash{\SetFigFont{12}{14.4}{\rmdefault}{\mddefault}{\updefault}{$\zeta_d^-$}%
}}}
\put(1251,-901){\makebox(0,0)[lb]{\smash{\SetFigFont{12}{14.4}{\familydefault}{\mddefault}{\updefault}{$\zeta_d^-\!+\!\epsilon$}%
}}}
\end{picture}
\caption{\label{fig:gamma-path}}
\end{centering}
\end{figure}

\begin{Lemma} \label{th:bastard}
On each $B_d = B \cap Y_d$ one has
\begin{equation} \label{eq:estimate}
  ||\nabla \pi_d||^2 \geq \nu^{-1} |d|^{1/2}
  \min\big( |\pi_d - \zeta_d^+|, |\pi_d - \zeta_d^-| \big)
\end{equation}
where $\nu>0$ is a constant independent of $d$.
\end{Lemma}

\begin{Proof}
As when proving Lemma \ref{th:bott-estimate}, we may assume that the
K{\"a}hler form is standard. By using the $S^1$-action $(x,a,b,c,d)
\mapsto (x,r^2 a,r^3 b,r^3 c,r^4 d)$, we may also assume that $d>0$.
A simple computation shows that
\begin{align*}
 |\pi_d - \zeta_d^{\pm}| & \leq
 |b|^2 + |c|^2 + |a \mp \sqrt{d/3}|^2 \cdot
 |a \pm 2 \sqrt{d/3} |. \\
\intertext{ On the other hand, }
 ||\nabla \pi_d||^2 & = 4|b|^2 + 4|c|^2 +
 9|a - \sqrt{d/3}|^2 \cdot |a + \sqrt{d/3}|^2. \\
\intertext{ For $re(a) \geq 0$, $|a + \sqrt{d/3}| \geq \frac{1}{2}
|a + 2 \sqrt{d/3}|$ and $|a + \sqrt{d/3}| \geq \sqrt{d/3}$, so }
 ||\nabla \pi_d||^2 & \geq
 \sqrt{d} \big( \frac{4}{\sqrt{d}} |b|^2 + \frac{4}{\sqrt{d}} |c|^2
 + \frac{9}{2\sqrt{3}} |a-\sqrt{d/3}|^2 \cdot |a+ 2\sqrt{d/3}| \big) \\
 & \geq \nu^{-1} |d|^{1/2} |\pi_d - \zeta_d^+|
\end{align*}
where $\nu$ is $\geq 2\sqrt{3}/9$ and is also an upper bound for
$\sqrt{d}/4$ on $B$. The other part of the minimum in
\eqref{eq:estimate} takes care of the case $re(a) \leq 0$.
\end{Proof}

For sufficiently small $d>0$, the critical points of $\pi_d$ will
lie close to $0$, and so will the sphere $L_{d,\epsilon}$ for
$\epsilon \ll d$. Using Lemma \ref{th:bastard}, we can now estimate
the length of any flow line of the parallel transport vector field
along $\gamma_{d,\epsilon}$, as long as that line remains inside
$B$. For the straight pieces from $\zeta_d^-+\epsilon$ to the origin
and back, one gets
\begin{align*}
 \int ||\nabla \pi_d||^{-1} & \leq
 2\int_0^{\zeta_d^+} \nu^{1/2} d^{-1/4} s^{-1/2} \, ds
 \leq 100 \nu^{1/2} d^{1/2}, \\
\intertext{ and for the circle around $\zeta_d^+$, }
 \int ||\nabla \pi_d||^{-1} & \leq
 \int_0^{2\pi \zeta_d^+} \nu^{1/2} d^{-1/4} (\zeta_d^+)^{-1/2} \, ds
 \leq 100 \nu^{1/2} d^{1/2}.
\end{align*}
Arguing as in our discussion of vanishing cycles, we arrive at the
desired conclusion: for $0 < \epsilon \ll d$ small, parallel
transport $h_{\gamma_{d,\epsilon}}$ is well-defined near
$L_{d,\epsilon}$, and moreover the image \eqref{eq:monodromy-image}
still lies in $B$.

To obtain a more concrete picture, assume momentarily that $\Omega$
is the standard K{\"a}hler form on $\C^4$. Borrowing from
\cite{khovanov-seidel98,shapere-vafa98}, we consider projection to
the $a$-coordinate, $q_{d,z}: Y_{d,z} \rightarrow \C$. The fibres of
this are affine quadrics, three of which are singular, corresponding
to the solutions of
\begin{equation} \label{eq:a-equation}
a^3 - ad - z = 0.
\end{equation}
The $S^1$ part of the $\C^*$-action from \eqref{eq:slice-action} is
a Hamiltonian circle action which is fibrewise with respect to
$q_{d,z}$, and whose moment map is $\mu(a,b,c,d) = |c|^2-|b|^2$. The
intersection
\begin{equation} \label{eq:circle}
 C_{d,z,a} = \mu^{-1}(0) \cap q_{d,z}^{-1}(a) = \{(b,c) \;:\;
 |b|^2 = |c|^2, \; bc = -a^3+ad+z\}
\end{equation}
is a circle if $a$ is a regular value, and shrinks to a point for
the singular values. To any embedded path $\alpha: [0;1] \rightarrow
\C$ such that $\alpha(r)^3 - \alpha(r)d - z$ vanishes exactly for $r
= 0,1$ one can associate an embedded smooth Lagrangian sphere in
$Y_{d,z}$,
\begin{equation} \label{eq:fibered-spheres}
 \Lambda_\alpha = \bigcup_{r=0}^1 C_{d,z,\alpha(r)}.
\end{equation}
Suppose that $d>0$ and $z = \zeta_d^- + \epsilon$ for $0 < \epsilon
\ll d$, so that \eqref{eq:a-equation} has three real solutions, of
which the rightmost two are close to $a = \sqrt{d/3}$.
\begin{figure}[hb]
\begin{centering}
\begin{picture}(0,0)%
\includegraphics{matching.pstex}%
\end{picture}%
\setlength{\unitlength}{3947sp}%
\begingroup\makeatletter\ifx\SetFigFont\undefined%
\gdef\SetFigFont#1#2#3#4#5{%
  \reset@font\fontsize{#1}{#2pt}%
  \fontfamily{#3}\fontseries{#4}\fontshape{#5}%
  \selectfont}%
\fi\endgroup%
\begin{picture}(2866,964)(1118,-1019)
\put(3426,-901){\makebox(0,0)[lb]{\smash{\SetFigFont{12}{14.4}{\familydefault}{\mddefault}{\updefault}{$\alpha$}%
}}}
\put(1951,-901){\makebox(0,0)[lb]{\smash{\SetFigFont{12}{14.4}{\familydefault}{\mddefault}{\updefault}{$\beta$}%
}}}
\put(2101,-211){\makebox(0,0)[lb]{\smash{\SetFigFont{12}{14.4}{\familydefault}{\mddefault}{\updefault}{$t_{\beta}(\alpha)$}%
}}}
\end{picture}
\caption{\label{fig:matching}} \end{centering}
\end{figure}

\begin{Lemma} \label{th:bc}
Assuming that the K{\"a}hler form on $Y$ is standard, the vanishing
cycle \eqref{eq:l0} and its monodromy image
\eqref{eq:monodromy-image} are the Lagrangian spheres
\eqref{eq:fibered-spheres} associated to paths in $\C$ which are
isotopic to $\alpha$ and $t_{\beta}(\alpha)$, respectively. Here
$\alpha$ and $\beta$ are as in Figure \ref{fig:matching}, and
$t_\beta$ denotes the positive half-twist around $\beta$.
\end{Lemma}

\begin{Proof}
The parallel transport vector fields on $\pi_d: Y_d \rightarrow \C$
are invariant with respect to the $S^1$-action
\eqref{eq:slice-action}, and $d\mu$ vanishes on them. Since the
critical points are fixed points of the action and lie in
$\mu^{-1}(0)$, it follows that all vanishing cycles are
$S^1$-invariant and also lie in $\mu^{-1}(0)$. One sees easily that
any Lagrangian sphere in $Y_{d,z}$ with these properties is
necessarily of the form $\Lambda_\alpha$ for some path $\alpha$ as
described above. Concerning $L_{d,\epsilon}$, one knows in addition
that it must lie close to the critical point $Crit(\pi_d)^-$ which
has $a = \sqrt{d/3}$, hence the corresponding $\alpha$ must stay
close to that value in $\C$, which determines its isotopy class
uniquely. The same argument as before proves that
\eqref{eq:monodromy-image} is of the form $\Lambda_{\alpha'}$ for
some path $\alpha'$. As one moves $z$ along $\gamma_{d,\epsilon}$,
the two leftmost solutions of \eqref{eq:a-equation} get exchanged,
and more precisely perform a positive half-twist around each other,
moving along a circle.
The path $\alpha'$ is isotopic to the image of $\alpha$ under the
resulting monodromy map of the three-pointed plane, which is
precisely $t_\beta$ (for more details on this last step see
\cite[Lemma 6.15]{khovanov-seidel98}).
\end{Proof}

We now turn to the realistic situation. Let $\FF \rightarrow X$ be a
holomorphic line bundle over some complex manifold, $Y = (\FF
\setminus 0) \times_{\C^*} \C^4$ be the associated bundle for the
$\C^*$-action \eqref{eq:slice-action}, and $\pi: Y \rightarrow \C^2$
the map which is equal to \eqref{eq:a2-map} on each $\C^4$ fibre.
$Y_d$, $\pi_d$ and $Y_{d,z}$ are defined in analogy with the
notation above. We equip $Y$ with an arbitrary K{\"a}hler form
$\Omega$, and restrict that to $X$ by identifying the latter space
with the zero-section of $Y$.

Let $K,K'$ be two closed Lagrangian submanifolds of $X$. The first
step in the construction, which was trivial in the previously
considered case $X = point$, goes as follows. Normalize the cusp
curve of critical values by the map $n: \C \rightarrow \Sigma$,
$n(w) = (3w^2,2w^3)$. The pullback of $Crit(\pi) \rightarrow \Sigma$
is the projection
\begin{equation} \label{eq:normalize}
n^*Crit(\pi) \iso X \times \C \longrightarrow \C,
\end{equation}
and while the pullback of $\Omega$ is no longer positive, it is
still nondegenerate on each fibre, which is sufficient to define
symplectic parallel transport. Starting with our original Lagrangian
submanifolds, which lie in the fibre over zero of
\eqref{eq:normalize}, we get two smooth families of Lagrangian
submanifolds in the fibres nearby. Changing back to the original
parameter, and supposing that $d > 0$ is sufficiently close to zero,
one now has Lagrangian submanifolds
\[
 K_d,K_d' \subset Crit(\pi_d)^-
\]
which as $d \rightarrow 0$ converge to $K$, $K'$ respectively. We
take the associated relative vanishing cycles for the map $\pi_d$,
which are Lagrangian submanifolds
\begin{equation} \label{eq:d-epsilon}
L_{d,\epsilon}^{\,},\; L_{d,\epsilon}' \subset Y_{d,\zeta_d^- +
\epsilon}
\end{equation}
for $0< \epsilon \ll d$ (this is actually the variant vanishing
cycle construction from Remark \ref{th:twisted-bundle} with line
bundle $\LL = \FF^2$, so \eqref{eq:d-epsilon} could be nontrivial
$S^2$-bundles over $K,K'$; however, that won't happen in our
application, see Remark \ref{th:top-trivial} below).

Take a relatively compact open subset $U \subset X$ containing
$K,K'$, and an open subset of $V \subset Y$ which, with respect to
some metric on the vector bundle $Y \rightarrow X$, is the unit ball
bundle over $U$. For $0 < \epsilon \ll d$, \eqref{eq:d-epsilon} will
lie inside $V$, and one has the same estimates for $\nabla \pi$ on
$V$ as (\ref{th:bastard}) (they can be derived for instance by
covering $\bar{U}$ with finitely many open subsets over which $\FF$
is trivial, and applying the previous argument to each of them).
Hence, parallel transport along the path $\gamma_{d,\epsilon}$ is
well-defined near $L_{d,\epsilon}'$, and yields another Lagrangian
submanifold
\begin{equation} \label{eq:d-epsilon-two}
 L_{d,\epsilon}'' = h_{\gamma_{d,\epsilon}}(L_{d,\epsilon}') \subset Y_{d,
 \zeta_d^- + \epsilon} \cap V.
\end{equation}

\begin{Background} \label{th:background}
At this point, we need to recall some general facts about symplectic
associated bundles. Let $(X,\omega)$ be a symplectic manifold, and
$\LL \rightarrow X$ a complex line bundle with a hermitian metric
and compatible connection. Let $P \subset \LL$ be the unit circle
bundle, and $\alpha \in \Omega^1(P)$ the connection one-form. Our
normalization is that if $R$ is the rotational vector field on $P$
(whose orbits are $2\pi$-periodic), then $\alpha(R) \equiv 1$. Let
$(M,\eta)$ be any symplectic manifold with a Hamiltonian circle
action, whose Killing field is $Z$ and whose moment map is $\mu$, so
$i_Z \eta = -d\mu$. The associated symplectic fibre bundle is
\[
 Y = P \times_{S^1} M \longrightarrow X.
\]
Take the two-form $\Omega = \omega + \eta + d(\alpha \mu)$ on $P
\times M$, where $\omega$ is pulled back via $P\rightarrow X$. This
satisfies $\Omega((R,-Z),\cdot) = 0$, hence it descends in a unique
way to a two-form on $Y$, also called $\Omega$. Take $V_1,V_2 \in
T_xX$ and lift them in the unique way to horizontal tangent vectors
$V_1^\natural, V_2^\natural$ on the circle bundle, so
$\alpha(V_k^\natural) = 0$. Take also $W_1,W_2 \in T_mM$, and
project $(V_k^\natural,W_k)$ to tangent vectors in the quotient $Y$.
Then
\[
 \Omega((V_1^\natural,W_1),(V_2^\natural,W_2))
 = \eta(V_1,V_2) + \omega(W_1,W_2) + \mu(m) d\alpha(V_1^\natural,V_2^\natural).
\]
This shows that $\Omega$ is nondegenerate, hence a symplectic form,
in a neighbourhood of $P \times_{S^1} \mu^{-1}(0) \subset Y$.
Moreover, if $K$ is a Lagrangian submanifold of $X$, and $L$ an
$S^1$-invariant Lagrangian submanifold of $M$ such that $\mu|L
\equiv 0$, then the associated $L$-bundle over $K$,
\begin{equation} \label{eq:fibrewise-lagrangian}
\Lambda = (P|K) \times_{S^1} L,
\end{equation}
is a Lagrangian submanifold of $(Y,\Omega)$.

Now assume that $X$ and $M$ are K{\"a}hler, $\LL$ is a holomorphic
line bundle and the connection is compatible with this structure,
and the circle action on $M$ is part of a holomorphic $\C^*$-action
$\rho$. One can then identify $Y$ with the holomorphic associated
bundle $(\LL \setminus 0) \times_{\C^*} M$, and $\Omega$ is
K{\"a}hler near $P \times_{S^1} \mu^{-1}(0)$. To see this, take a
nowhere vanishing holomorphic section $e$ of $\LL$ over some open
subset $U$. This defines a holomorphic trivialization of $(\LL
\setminus 0) \times_{\C^*} M$ over $U$. Its normalized form
$e/||e||$ gives a corresponding trivialization of $P \times_{S^1}
M$, with respect to which the connection one-form is $\alpha = -d^c
\log ||e||$. The difference between the two trivializations is the
map
\begin{equation} \label{eq:rescale}
 U \times M \longrightarrow U \times M, \quad 
 (x,m) \longmapsto (x,\rho_{\textit{exp}(h(x))}(m)),
\end{equation}
involving the radial part of the $\C^*$-action and $h = -\log
||e||$. The pullback of $\Omega$ by \eqref{eq:rescale} is $\eta +
\rho_{\textit{exp}(h)}^*\omega - dh \wedge
\rho_{\textit{exp}(h)}^*(d^c \mu) + d^c h \wedge
\rho_{\textit{exp}(h)}^* d\mu - dd^c h \cdot
\rho_{\textit{exp}(h)}^*\mu$, which is obviously of type $(1,1)$.
\end{Background}

Returning to our discussion: starting with a K{\"a}hler form on $X$,
a hermitian metric on $\FF$ and the standard K{\"a}hler form on
$\C^4$, we construct an associated form $\Omega$ on $Y = \FF
\times_{\C^*} \C^4$ which is K{\"a}hler at least in a neighbourhood
of the set $\mu^{-1}(0) = \{|b| = |c|\}$. Given $(d,z) \in \C^2
\setminus \Sigma$, a Lagrangian submanifold $K \subset X$, and a
path $\alpha$ in $\C$ of the same kind as before, one can define a
Lagrangian submanifold
\begin{equation} \label{eq:sphere-bundles}
 \Lambda_{d,z,K,\alpha} \subset Y_{d,z}
\end{equation}
by starting with the Lagrangian sphere \eqref{eq:fibered-spheres}
and applying the construction \eqref{eq:fibrewise-lagrangian}. Under
the map $Y_{d,z} \rightarrow X$, this is an $S^2$-bundle over $K$.
The following result is the generalization of Lemma \ref{th:bc} and
has the same proof:

\begin{Lemma} \label{th:fibred-matching}
Suppose that we use the K{\"a}hler form obtained from a K{\"a}hler
form on $X$, a hermitian metric on $\FF$, and the standard
K{\"a}hler form on $\C^4$. Let $K,K' \subset X$ be closed Lagrangian
submanifolds, and consider the Lagrangian submanifolds from
\eqref{eq:d-epsilon}, \eqref{eq:d-epsilon-two}. Then
\[
 L_{d,\epsilon} = \Lambda_{d,\zeta_d^- + \epsilon,K,\alpha},
 \qquad
 L_{d,\epsilon}'' = \Lambda_{d,\zeta_d^- + \epsilon,K',\alpha'}
\]
for paths $\alpha$ isotopic to the one from Figure
\ref{fig:matching}, and $\alpha'$ isotopic to $t_\beta(\alpha)$.
\qed
\end{Lemma}

\subsection{Floer cohomology: background\label{subsec:background}}

Let $(M,\omega)$ be a K{\"a}hler manifold such that $\omega$ is
exact and the underlying complex structure is Stein, meaning that
there is an exhausting plurisubharmonic function $\psi$. We assume
that $c_1(M) = 0$ and $H^1(M) = 0$. Let $L,L'$ be closed connected
Lagrangian submanifolds of $M$ with $H_1(L) = H_1(L') = 0$ and
$w_2(L) = w_2(L') = 0$. Then there is a well-defined Floer
cohomology group
\[
HF(L,L') = H(CF(L,L'),d_J)
\]
which is a finitely generated, relatively graded abelian group. One
can of course replace $\Z$ with any other abelian coefficient group,
and the universal coefficient theorem holds as usual. Relatively
graded means that there is a $\Z$-grading which is unique up to an
overall constant shift. We recall briefly the definition, which is
essentially Floer's original one \cite{floer88c} except for the
orientations of moduli spaces, which come from \cite{fooo}. First,
after a small Lagrangian perturbation, one may assume that the
intersection $L \cap L'$ is transverse. One then defines the Floer
chain complex to be the abelian group
\[
 CF(L,L') = \bigoplus_{x \in L \cap L'} O_x
\]
where $O_x$ is the orientation group of $x$. Formally, this is an abelian
group canonically associated to $x$ and generated by two elements -- labelled the possible ``coherent
orientations'' of $x$ --
with the relation that the sum of these orientations is zero. Hence,
$O_x \iso \Z$ but not canonically so. The association to $x$ of a
well-defined orientation group proceeds essentially as in the case of
Hamiltonian Floer cohomology 
\cite{floer-hofer93}, but there is a ``family index anomaly'' due to
which the consistency of the definition requires $Spin$ structures
on $L$ and $L'$ \cite{fooo}. Of course, our standing topological
assumptions on $L$ and $L'$ imply that
such structures exist and are unique. Similarly, one can associate
to any pair of intersection points $x,y$ a relative Maslov index
$\Delta\mu(x,y) \in \Z$, which satisfies $\Delta\mu(x,y) +
\Delta\mu(y,z) = \Delta\mu(x,z)$ and establishes the relative
grading. The differential $d_J$ is defined by considering solutions
of Floer's equation
\begin{equation} \label{eq:floer}
\left\{
\begin{aligned}
 & u: \R \times [0;1] \rightarrow M, \\
 & u(s,0) \in L, \quad u(s,1) \in L', \\
 & \partial_s u + J_t(u) \partial_t u = 0, \\
 & \lim_{s \rightarrow \pm\infty} u(s,\cdot) = x_{\pm}
\end{aligned}
\right.
\end{equation}
where $J = (J_t)_{0 \leq t \leq 1}$ is a generic smooth family of
$\omega$-compatible almost complex structures, which all agree with
the given complex structure outside a compact subset, and $x_\pm \in
L \cap L'$. More precisely, $d_J(x_+) = \sum_{x_-} n_{x_+,x_-} x_-$
where $n_{x_+,x_-}$ counts isolated solutions $u$ of
\eqref{eq:floer} (mod translation in the $s$-variable), with a sign
that can be canonically encoded as an isomorphism $\delta_u: O_{x_+}
\iso O_{x_-}$. 

The geometric assumptions set out at the beginning of the section enter in two
crucial ways (one elementary, and one going back to Floer
\cite{floer88c}). First, the  fact that our almost complex structures are
standard at infinity, together with the Stein property of $M$,
allows one to apply the maximum principle to $\psi \circ u$; as a
consequence, all solutions of \eqref{eq:floer}
remain within a fixed 
compact subset of $M$. Secondly, the exactness of $\omega$ and the
fact that $H^1(L) = H^1(L') = 0$ ensure that there is a well-defined
action functional on the space of paths from $L$ to $L'$, which
implies that the moduli spaces of \eqref{eq:floer} have well-behaved
compactifications (bounded energy, no bubbling).  The other
assumptions are of lesser importance, though of some relevance given
Conjecture \ref{th:conj}.  Existence of spin structures on $L$ and
$L'$ enables one to define Floer theory with integral rather than mod
2 coefficients; vanishing
of $c_1(M)$ and $H^1(M)$ are respectively relevant to the existence
and uniqueness of a $\Z$-grading in Floer cohomology, as discussed in
Section \ref{sec:misc}\ref{subsec:gradings}. 

\begin{Remark}
Even if the K{\"a}hler form is defined only on some open subset $U
\subset M$ which is holomorphically weakly convex (meaning that
holomorphic discs cannot touch $\partial U$ from the inside, unless
they are completely contained in $\partial U$), one can still define
Floer cohomology groups for Lagrangian submanifolds $L,L' \subset U$
satisfying the conditions set out above. The reason is that to
achieve transversality of the moduli spaces, it is sufficient to
consider almost complex structures $J_t$ which agree with the given
complex structure outside any given neighbourhood $U'$ of $L \cup
L'$. By taking a $U'$ whose closure is inside $U$ and applying weak
convexity, one sees that all solutions of \eqref{eq:floer} remain in
$U$, so that their energy can be estimated by the symplectic area
and hence the action functional. The most obvious application of
this is to sublevel sets $U = \{\psi < C\}$.
\end{Remark}

The most important aspect of Floer cohomology is its strong invariance
properties. A Hamiltonian isotopy of $M$ is a path $(g_t)$ of
symplectomorphisms starting at $g_0=id$ and defined by the flow of the vector
field associated to some smooth time-dependent function $H_t$, for
which $H:M\times I\rightarrow \R$ has compact support.  Floer
\cite{floer88c} proved invariance  of $HF^*(L,L')$
under Hamiltonian 
isotopies of either $L$ or $L'$.   Now let
$\{L_t=\phi_t(L)\}_{t\in[0,1]}$ be an arbitrary Lagrangian
isotopy of a closed 
Lagrangian submanifold $L$ in an exact symplectic manifold
$(M,\omega,\theta)$, where $d\theta=\omega$. The isotopy is said to be
exact if $[\theta|_{L_t}]\in H^1(L_t)\cong H^1(L)$ is constant. Note
that this makes sense, in that there is a canonical identification
$H^1(L_t)\cong H^1(L_0)$ for all $t$. By explicitly writing down and
integrating the appropriate vector fields, one sees that such an isotopy can be
embedded in a global ambient 
Hamiltonian isotopy of $M$; this is analogous to the characterisation
of Hamiltonian symplectomorphisms as those having zero flux,
cf. \cite[Chapter 10]{mcduff-salamon00}. That is,
there is some 
$(g_t)$ with $g_t(L)=L_t$. Under our standing
assumptions $H_1(L)=H_1(L')=0$ exactness is automatic, and it follows
that the Floer cohomology $HF^*(L,L')$ is
invariant under arbitrary Lagrangian isotopies. 

For a second invariance property, suppose that we have an isotopy
$(\omega^s)$ of symplectic forms on $M$, together with closed
submanifolds $L$, $L'$ which are $\omega^s$-Lagrangian for every $s$.
Suppose as usual that $H_1(L)=H_1(L')=0$ and that the other geometric
assumptions required for well-definition of Floer cohomology hold.  In
particular, suppose that all the $\omega^s$ are K\"ahler forms making $M$
geometrically bounded at infinity, for instance making $M$ Stein.
Then we claim $HF^*(L,L')$ is independent of the particular symplectic
form $\omega^s$ (given that the Floer differential counts solutions to
an equation defined without explicit mention of $\omega^s$ this is
perhaps not as surprising as it first seems). The result is proved
using a parametrised version  of the
Floer equation \eqref{eq:floer}, as in Floer's original
\cite{floer88c}.  The main technical difficulties stem from the
parameter values where birth-death processes occur for the
intersection points of the Lagrangian submanifolds; a careful
treatment of these issues has been given in \cite{lee01}. Note that
the discussion in \cite{lee01} analyses bifurcations occuring in
rather general ``one-parameter homotopies of Floer data'', and applies
equally to a 
parametrized Floer equation in which the almost complex structures are
compatible with a smoothly varying family of symplectic forms.  Another
approach would be to combine parametrized moduli spaces with the
continuation map technique, noting the energy bounds
required for compactness of spaces of solutions to the continuation
map equation carry over essentially as usual to this case. A closely
related but more difficult statement -- proving invariance of
symplectic homology of 
Stein manifolds under continuous variation of the symplectic form -- was proved
by Viterbo in \cite{viterbo99}.


Combining the two statements, if one has a holomorphic submersion $Y
\rightarrow T$ 
whose fibres are Stein, a path $\gamma: [0;1] \rightarrow T$ and
families $L_r$, $L'_r$ of closed Lagrangian submanifolds in the
fibres $Y_{\gamma(r)}$, with the required additional conditions to
make $HF^*(L_r,L'_r)$ well-defined, then it is the same for all $r$.
A helpful, allbeit informal, general principle is thus that Floer cohomology is
invariant under smooth deformation of all geometric objects involved,
as long as one remains within the class where it is well-defined.

\subsection{Floer cohomology: computations}

We will need two simple Floer cohomology computations for the
geometric situations studied in Sections
\ref{sec:slice}\ref{subsec:2-coincide} and \ref{subsec:3-coincide}.
They have the flavour of a ``K\"unneth formula'' and ``Thom
isomorphism'' respectively. The first computation takes place in the
following context:
\begin{itemize}
\item
 $Y$ is a complex manifold with a holomorphic function $\pi:
Y \rightarrow \C$. We have a complex submanifold $X \subset Y$ and
an isomorphism between a neighbourhood of that submanifold and a
neighbourhood of $X \times \{0\}^3 \subset X \times \C^3$, such that
the following diagram commutes:
\begin{equation} \label{eq:local-quad}
\begin{CD}
 Y @>{\text{local $\iso$ defined near $X$}}>> X \times \C^3 \\
 @V{\pi}VV @V{a^2+b^2+c^2}VV \\
 \C @>{\quad\quad\qquad\qquad\quad}>> \C
\end{CD}
\end{equation}
\item
 $Y$ is Stein, carries an exact K{\"a}hler form $\Omega$, and
satisfies $c_1(Y) = 0$ (which implies that $c_1(Y_t)$ and $c_1(X)$
also vanish). Moreover, $H^1(Y_t) = 0$ for small $t \neq 0$, and
$H^1(X) = 0$.
\end{itemize}

We equip $X$ and the smooth fibres $Y_t$ with the restrictions of
$\Omega$. Let $K,K'$ be closed Lagrangian submanifolds of $X$ which
have the properties necessary to define $HF(K,K')$, and consider for
sufficiently small $t \neq 0$ the associated relative vanishing
cycles $L_t,L_t' \subset Y_t$. Since these are products of $K,K'$
with $S^2$, their Floer cohomology $HF(L_t,L_t')$ is again
well-defined, and is independent of $t$ by the basic invariance
principle discussed above.

\begin{Lemma} \label{th:kunneth}
$HF(L_t,L_t') \iso HF(K,K') \otimes H^*(S^2)$, where $H^*(S^2)$
carries its standard grading.
\end{Lemma}

\begin{Proof}
One can find finitely many holomorphic functions
$\phi_1,\dots,\phi_l: Y \rightarrow \C$ whose common vanishing set
is $X$ (this is a general result about complex submanifolds of Stein
manifolds \cite{forster67}; however, note that all our applications
will be in the affine algebraic context where the counterpart is
trivial, so we are appealing to it only to keep the current
exposition general). Take a sublevel set $U = \{\psi(y) < C\} \cap X
\subset X$ which contains both $K$ and $K'$, and consider the open
subset $V = \{\psi(y) < C, \;\; |\phi_1(y)| < \delta, \dots,
|\phi_l(y)| < \delta\}$ for some $\delta>0$. Since $U$ is relatively
compact, one can make $\delta$ sufficiently small so as to ensure
that $V$ is contained in the neighbourhood of $X$ where the
isomorphism \eqref{eq:local-quad} is defined. By taking $t$ small,
one can achieve that the relative vanishing cycles $L_t,L_t'$ lie in
$Y_t \cap V$. By definition $V$ is holomorphically weakly convex, so
for a suitable choice of almost complex structure, the definition of
$HF(L_t,L_t')$ is local, meaning that all solutions of
\eqref{eq:floer} stay inside $Y_t \cap V$.

On $V$ there is another K{\"a}hler form $\Omega^{(1)}$, which is
obtained by taking the product of $\Omega|X$ and the standard form
on $\C^3$, and pulling that back by \eqref{eq:local-quad}. We can
consider the linear family of K{\"a}hler forms $\Omega^{(s)}$
interpolating between $\Omega^{(0)} = \Omega|V$ and $\Omega^{(1)}$.
For each $s \in [0;1]$ there are relative vanishing cycles
constructed from $K,K'$ using $\Omega^{(s)}$. These will be
well-defined inside $Y_t \cap V$ for sufficiently small $t \neq 0$
(one can see this explicitly by estimating the parallel transport
vector field, using Lemma \ref{th:bott-estimate}), and the Floer
cohomology is local in the same sense as before. It follows that
$HF(L_t,L_t')$ is isomorphic to the Floer cohomology of the
corresponding vanishing cycles for $\Omega^{(1)}$. But in the
coordinates given by the isomorphism \eqref{eq:local-quad}, these
cycles are simply $K \times \sqrt{t}S^2$ and $K' \times
\sqrt{t}S^2$, compare \eqref{eq:lzero}. At this point, the
K{\"u}nneth isomorphism in Floer cohomology (which has the same
essentially trivial proof as in ordinary Morse theory) finishes the
argument.
\end{Proof}

The second situation we are concerned with is as follows:
\begin{itemize}
\item
 $Y$ is a complex manifold with a holomorphic map $\pi: Y
\rightarrow \C^2$. We have a complex submanifold $X \subset Y$
equipped with a holomorphic line bundle $\FF$, and an isomorphism
between a neighbourhood of that submanifold and a neighbourhood of
the zero-section inside $(\FF \setminus 0) \times_{\C^*} \C^4$,
where the associated bundle is formed with respect to
\eqref{eq:slice-action}. This should fit into a commutative diagram
\begin{equation} \label{eq:local-a2}
\begin{CD}
 Y
 @>{\text{local $\iso$ defined near $X$}}>>
 {(\FF \setminus 0) \times_{\C^*} \C^4} \\
 @V{\pi}VV @V{p}VV \\
 {\C^2} @>{\quad\qquad\qquad\qquad\quad}>> \C^2,
\end{CD}
\end{equation}
where $p$ is given on each fibre by \eqref{eq:a2-map}.

\item
 $Y$ is Stein, carries an exact K{\"a}hler form $\Omega$, and
satisfies $c_1(Y) = 0$ (which implies that $c_1(Y_{d,z})$, for a
regular value $(d,z)$ of $\pi$, and $c_1(X)$ also vanish). We also
require $H^1(Y_{d,z}) = 0$ for small regular values $(d,z)$, and
$H^1(X) = 0$. Finally, $\FF$ should be a subbundle of the trivial
$\C^2$-bundle over $X$.
\end{itemize}

Take closed Lagrangian submanifolds $K,K' \subset X$, satisfying the
conditions above so in particular $HF(K,K')$ is well-defined. The
construction in \ref{subsec:a2} above associates to these new
Lagrangian submanifolds as in \eqref{eq:d-epsilon},
\eqref{eq:d-epsilon-two} inside the fibre $Y_{d,\zeta_d^- +
\epsilon}$ for $0 < \epsilon \ll d$.

\begin{Remark} \label{th:top-trivial}
The last-made assumption in the list above implies that there is a
short exact sequence of vector bundles $0 \rightarrow \FF
\rightarrow \C^2 \rightarrow \FF^{-1} \rightarrow 0$. Because $X$ is
Stein, it follows that $\FF \oplus \FF^{-1}$ is trivial, and
therefore so is the vector bundle
\[
(\FF \setminus 0) \times_{\C^*} \C^4 = \C \oplus \FF^{-2} \oplus
\FF^2 \oplus \C = (\FF \oplus \FF^{-1})^{\otimes 2}.
\]
Another consequence is that $P(\FF^2 \oplus \C)$ is trivial, which
by construction means that the Lagrangian submanifolds we
constructed will be diffeomorphic to $K \times S^2$, $K' \times S^2$
respectively. In particular they again satisfy $w_2 = 0$, so the
Floer cohomology $HF(L_{d,\epsilon},L''_{d,\epsilon})$ is
well-defined.
\end{Remark}

\begin{Lemma} \label{th:thom}
$HF(L_{d,\epsilon},L''_{d,\epsilon}) \iso HF(K,K')$.
\end{Lemma}

\begin{Proof}
Suppose that $K,K'$ intersect transversally. The first steps of the
proof are the same as in Lemma \ref{th:kunneth}: one can achieve
that the Lagrangian submanifolds concerned lie inside the set where
the local isomorphism \eqref{eq:local-a2} is defined, and moreover
Floer cohomology can be localized to that subset. This means that
from now on, our computations will all take place in the fibre
bundle
\begin{equation} \label{eq:fibre-bundle}
 Y_{d,\zeta_d^- + \epsilon}' = (\FF \setminus 0) \times_{\C^*}
 \{a^3 - ad + bc = \zeta_d^- + \epsilon\} \subset
 (\FF \setminus 0) \times_{\C^*} \C^4 \rightarrow X.
\end{equation}
Moreover, one can replace the given K{\"a}hler form $\Omega$ by an
$\Omega^{(1)}$ which is constructed from $\Omega|X$, a hermitian
metric on $\FF$, and the standard form on $\C^4$ ($\Omega^{(1)}$ is
strictly speaking defined only in a neighbourhood of the
zero-section, but all our arguments will take place inside that
neighbourhood). The Lagrangian submanifolds can then be described
explicitly as in Lemma \ref{th:fibred-matching}, in fact after a
suitable Lagrangian isotopy we may assume that
\begin{equation} \label{eq:ls}
 L_{d,\epsilon} = \Lambda_{d,\zeta_d^- + \epsilon,K,\alpha}, \quad
 L_{d,\epsilon}'' = \Lambda_{d,\zeta_d^- + \epsilon,K',t_\beta(\alpha)}
\end{equation}
where the paths $\alpha$, $t_\beta(\alpha)$ are as in Figure
\ref{fig:matching}, as opposed to merely lying in the same isotopy
class. Note that the two paths intersect only in one endpoint, which
is the rightmost solution $a_+$ of $a^3 - ad = \zeta_d^- +
\epsilon$. By construction \eqref{eq:fibered-spheres} this implies
that in the fibre over each point $x \in K \cap K'$, there is a
unique (and transverse) intersection point $y \in L_{d,\epsilon}
\cap L_{d,\epsilon}''$, given by $b = c = 0$ and $a = a_+$.

Consider solutions $u$ of Floer's equation for the Lagrangian
submanifolds \eqref{eq:ls} inside \eqref{eq:fibre-bundle} using the
standard complex structure as $J$. Under projection to the
$a$-variable, $u$ gets mapped to a finite energy holomorphic map $\R
\times [0;1] \rightarrow \C$ whose boundary lies on the paths
$\alpha$, $t_\beta(\alpha)$. This is necessarily constant equal to
$a_+$, and that implies that $u$ lies in the subset of
\eqref{eq:fibre-bundle} where $a = a_+$, and because the only points
of our Lagrangian submanifolds which satisfy $a = a_+$ also have $b
= c = 0$, we find that the $(b,c)$ components of $u(s,t)$ vanish for
$t = 0,1$. The $b$ component, for instance, can be viewed as a
holomorphic section of the pullback bundle
\begin{equation} \label{eq:pullback-line}
 u^*\FF^{-2} \rightarrow \R \times [0;1]
\end{equation}
which is zero along $\R \times \{0;1\}$, so unique continuation shows
that it is identically zero; similarly for the $c$-component. This
sets up a bijective correspondence between solutions of Floer's
equation in our fibre bundle \eqref{eq:fibre-bundle} with boundary
conditions \eqref{eq:ls}, and those of the corresponding equation in
$X$ with boundary condition $K,K'$ (the correspondence is given by
projection in one direction, and conversely by lifting to the
submanifold $X_+ \subset Y$ defined by $a = a_+$, $b = c = 0$). This
would be a way of proving our result if we could assume that the
standard complex structure was regular in the sense of Floer theory
for the pair $(K,K')$. That is not a realistic assumption, but the
way to repair the argument is a standard exercise in Floer theory,
so we will only sketch it (one could even claim that this step is
trivial, at the price of relying on somewhat complicated virtual
perturbation theory).

Fix small neighbourhoods $U_1 \subset \bar{U}_1 \subset U_2$ of
$X^+$ inside \eqref{eq:fibre-bundle}. On the base $X$, choose
$\bar{J} = (\bar{J}_t)$ which is a small compactly supported
perturbation of the standard complex structure, making the solutions
of Floer's equation for $(K,K')$ regular. From $\bar{J}$ and the
given connection on $\FF$ one gets an induced $t$-dependent almost
complex structure on \eqref{eq:fibre-bundle} which is compatible
with our symplectic form. We choose a $J = (J_t)$ which agrees with
this induced almost complex structure on $U_1$, and is equal to the
standard complex structure outside $U_2$. This can be done in such a
way that $J$ is still everywhere a small perturbation of the
standard complex structure, and a Gromov compactness argument will
tell us that all solutions of Floer's equation for this $J$ are
contained in $U_1$, so that for all practical purposes this is the
almost complex structure induced from $\bar{J}$ and the connection.
As a consequence, we still have all the properties used in the
argument above ($J_t$-holomorphicity of the projection to the
$a$-variable, the fact that the $b$ component can be viewed as a
holomorphic section of \eqref{eq:pullback-line} and similarly for
$c$, and $(J_t,\bar{J}_t)$-holomorphicity of projection to the base
$X$), but now with the added benefit of regularity.
\end{Proof}

\section{Symplectic geometry of $\chi|\Slice_m$}

With the symplectic techniques at hand, we now return to the
specifics of $\chi|\Slice_m$. In particular we define the Lagrangian
submanifolds $L_{\wp}$, the monodromy maps $h^{\resc}_\beta$, and
therefore $\Kh_{\symp}$. Its invariance under the Markov moves (and
its behaviour under adding unlinked unknotted components) is a swift
consequence of the preceding material. As in the previous section we
deal with Floer cohomology as a relatively graded group only, but
this will be remedied in the following section.

\subsection{Open braids and parallel transport maps\label{subsec:braid}}

Take the affine transverse slice $\Slice_m$ from Section
\ref{sec:slice} inside $\gg = \sl_{2m}$. Recall that points $t \in
\Conf_0^{2m}(\C) \subset \hh/W$, which correspond to pairwise
different eigenvalues $(\mu_1,\dots,\mu_{2m})$ with
\begin{equation} \label{eq:normal}
\mu_1 + \cdots + \mu_{2m} = 0,
\end{equation}
are regular values of $\chi|\Slice_m$. We denote by $\Y_{m,t} =
(\chi|\Slice_m)^{-1}(t)$ the fibre over $t$. Occasionally we will
extend this notation to points of the whole configuration space
$\Conf_{2m}(\C)$, with the understanding that in that case,
$\Y_{m,t}$ is really the fibre over the normalized configuration
$(\mu_1 - (\sum_k \mu_k)/2m, \dots, \mu_{2m} - (\sum_k \mu_k)/2m)$.

\begin{Lemma} \label{th:1st-betti}
$H^1(\Y_{m,t}) = 0$.
\end{Lemma}

This can be read off from the literature. Because of the
simultaneous resolution, $\Y_{m,t}$ is diffeomorphic to the fibre
over $0$ of $\tilde{\Slice}_m \rightarrow \hh$. By \cite[p.\
50]{slodowy82} that fibre deformation retracts onto its ``compact
core'', which is the preimage of the nilpotent $x = n^+$ under
$\tilde{\Slice}_m \rightarrow \Slice_m$, or equivalently by
definition the set of flags stabilized by $n^+$. The inclusion of
that set into the full flag variety induces a surjective map on
cohomology \cite[p.\ 60]{slodowy82}, and of course that variety has
$H^1 = 0$. Alternatively one can follow \cite{khovanov02c} and
appeal to the presentation of the entire algebra $H^*(\Y_{m,t})$
given in \cite{deconcini-procesi81}.

When it comes to choosing a K{\"a}hler metric on $\Slice_m$, we are
guided by the requirements of the (rescaled) parallel transport
construction and the proof of Lemma \ref{th:fibre-bundle-2}. Fix
some real number $\alpha > m$. For each $i = 2,4,\dots,2m$ take the
functions $\xi_i(z) = |z|^{2\alpha/i}$ on $\C$. By applying $\xi_i$
to each coordinate of $\Slice_m$ on which $\lambda$ acts with weight
$i$, which in terms of \eqref{eq:y-matrix} means each entry of
$y_{1i}$, and summing up all these terms, one gets a proper
$C^1$-function $\xi$ on $\Slice_m$. Now find compactly supported
functions $\eta_k$ on $\C$ such that $\psi_i = \eta_i + \xi_i$ is
$C^{\infty}$, and add the $\psi_i$ in the same way as before to form
another function $\psi$ on $\Slice_m$.

\begin{Lemma} \label{th:rescaling}
$\psi$ is asymptotically homogeneous for the radial part of the
$\C^*$-action $\lambda$, in the sense that
\[
 \lim_{r \rightarrow \infty}
 \frac{\psi \circ \lambda_r}{r^{2\alpha}} = \xi
\]
where the convergence is uniform in $C^1$-sense.
\end{Lemma}

\begin{Proof}
As $r \rightarrow \infty$, the rescaled functions $\eta_i(r^i z)$ on
$\C$ are supported on progressively smaller neighbourhoods of the
origin. Their $C^0$-norms are of course uniformly bounded, and their
derivatives grow like $r^i$. Since $i \leq 2m < 2\alpha$, the limit
$\eta_i(r^i z)/r^{2\alpha}$ goes to zero uniformly.
\end{Proof}

Let $\tilde{\xi},\tilde{\psi}$ be the lifts of $\xi,\psi$ to the
simultaneous resolution $\tilde{\Slice}_m$.

\begin{Lemma}
The union of the critical points of $\tilde{\psi}$ on the fibres of
$\tilde{\chi}|\tilde{\Slice}_m$ forms a subset of $\tilde{\Slice}_m$
which projects properly to $\hh$.
\end{Lemma}

\begin{Proof}
Suppose on the contrary that there is a sequence of points
$\tilde{y}_j \in \tilde{\Slice}_m$ which are critical points in the
respective fibres of $\tilde{\chi}|\tilde{\Slice}_m$, such that
$\tilde{y}_j$ goes to infinity but $\tilde{\chi}(\tilde{y}_j)$
remains bounded. After rescaling with a suitable sequence
$\lambda_{r_j}$ and applying Lemma \ref{th:rescaling}, one obtains a
limiting point $\tilde{y} \in \tilde{\Slice}_m$ whose projection to
$\Slice_m$ lies on the unit sphere for the obvious identification
$\Slice_m \iso \C^{4m-1}$, such that $\tilde{\chi}(\tilde{y}) = 0$,
and which is a critical point for $\tilde{\xi}$ on
$\tilde{\chi}^{-1}(0) \cap \tilde{\Slice}_m$. But that is is
impossible due to the homogeneity of $\xi$.
\end{Proof}

We will assume from now on that the $\eta_k$ have been chosen in
such a way that $-dd^c \psi_k > 0$ everywhere, and equip $\Slice_m$
with the metric defined by the K{\"a}hler form $\Omega = -dd^c
\psi$.

\begin{Lemma}
Outside a compact subset, we have an inequality $|| \nabla \psi ||^2
< \rho \psi$ for some $\rho>0$.
\end{Lemma}

\begin{Proof}
An explicit computation shows that $|| \nabla \psi_k ||^2 \leq c +
\rho\psi$ for some $c,\rho>0$. Since our metric is the product of
K{\"a}hler metrics on each coordinate, we can add this up and
suitably adjust the constants, to get $|| \nabla \psi ||^2 \leq c +
(\rho/2)\psi$. But $(\rho/2)\psi > c$ outside a compact subset.
\end{Proof}

In view of the two Lemmas above, the argument from Section
\ref{sec:symplectic}\ref{subsec:parallel} shows that the family
$\chi^{-1}(\Conf_{2m}^0(\C)) \cap \Slice_m \rightarrow
\Conf_{2m}^0(\C)$ has well-defined rescaled parallel transport maps
(defined on arbitrarily large compact subsets of the fibres, or even
on the entire fibres if one is willing to take the slightly more
complicated route indicated in Remark \ref{th:transport}). If
$\beta$ is a piecewise smooth path $[0;1] \rightarrow
\Conf_{2m}^0(\C)$, the associated rescaled parallel transport is
denoted by
\begin{equation} \label{eq:y-transport}
h^{\resc}_\beta: \Y_{m,\beta(0)} \longrightarrow \Y_{m,\beta(1)}.
\end{equation}
As before, we extend this notation to arbitrary open braids, which
are paths in $\Conf_{2m}(\C)$, with the understanding that one
translates each $\beta(s)$ by an $s$-dependent amount so that
\eqref{eq:normal} is again satisfied.

We will also need a version of this discussion for the critical
point set fibration \eqref{eq:critical-fibration}. The base space
$\hh^{\mult,\reg}/W^{\mult}$ of that can be identified with
$\Conf_{2m-2}(\C)$ by forgetting the first two eigenvalues. The
total space comes with a natural map $\CC_m \rightarrow \Slice_m$
which is an embedding on each fibre, and we pull back the K{\"a}hler
form by it. As explained in Section
\ref{sec:slice}\ref{subsec:partial}, the arguments from Lemma
\ref{th:fibre-bundle-2} can be easily adapted to show that
\eqref{eq:critical-fibration} is a differentiable fibre bundle. A
combination of these arguments and the ones used above proves that
there are well-defined parallel transport maps
\begin{equation} \label{eq:c-transport}
 h^{\resc}_{\bar{\beta}}: \CC_{m,\bar{\beta}(0)} \rightarrow
 \CC_{m,\bar{\beta}(1)}
\end{equation}
for any path $\bar\beta$ in $Conf_{2m-2}(\C)$. Lemma
\ref{th:induction} says that the fibre $\CC_{m,\bar{t}}$ over a
point $\bar{t} \in \Conf_{2m-2}^0(\C) \subset \Conf_{2m-2}(\C)$ can
be identified with the corresponding space $\Y_{m-1,\bar{t}}$. This
is compatible with our choice of symplectic forms (provided one
takes the same $\alpha$ and functions $\psi_k$ for both $m$ and
$m-1$, see Lemma \ref{th:alpha} below). Hence, if $\bar\beta$ lies
in $Conf_{2m-2}^0(\C)$ then the parallel transport
\eqref{eq:c-transport} is the same as the corresponding map
\eqref{eq:y-transport} for $m-1$. Note that even though there is no
canonical isomorphism $\CC_{m,\bar{t}} \iso \Y_{m-1,\bar{t}}$ for
general $\bar{t} \in Conf_{2m-2}(\C)$, one can partially remedy this
by moving $\bar{t}$ into the subset $Conf_{2m}^0(\C)$ by
translation, and then combining the isomorphism defined there with
parallel transport \eqref{eq:c-transport} along the family of
translated configurations to get back to the original fibre.

\subsection{Lagrangian submanifolds from matchings}

Take $t = (\mu_1,\dots,\mu_{2m}) \in \Conf_{2m}(\C)$. A crossingless
matching $\wp$ with endpoints $t$ is a collection of $m$ disjoint
embedded unoriented arcs $(\delta_1,\dots,\delta_m)$ in $\C$ which
join together the points of $t$ in pairs. For the moment, we include
an ordering of the arcs as part of the data (although that will be
dropped at some point later on), and order the configuration
correspondingly, so that $\delta_k$ has endpoints
$\mu_{2k-1},\mu_{2k}$. We will associate to each such $\wp$ a
Lagrangian submanifold
\begin{equation} \label{eq:matching-lagrangian}
 L_{\wp} \subset \Y_{m,t}
\end{equation}
which is diffeomorphic to $(S^2)^m$ and unique up to Lagrangian
isotopy. Choose a path $[0;1) \rightarrow Conf_{2m}(\C)$ starting at
$t$ which moves the points as follows: the endpoints of
$\delta_2,\dots,\delta_m$ remain fixed, and the two endpoints of
$\delta_1$ move towards each other along that arc, colliding in the
limit $s \rightarrow 1$. For simplicity we assume that the arc is a
straight line near its midpoint, and that the colliding points move
towards each other with the same speed for $s$ close to $1$. After
translating to meet the normalization condition \eqref{eq:normal} at
all times, we get a path $\gamma: [0;1] \rightarrow \hh/W$ such that
the point $\gamma(1)$ corresponds to a collection of eigenvalues
$(\mu_1',\dots,\mu_{2m}')$ where $\mu_1' = \mu_2'$, and $\mu_k' =
\mu_k + \mu_1/(m-1) - \mu_1'/(m-1)$ for $k \geq 3$. Note that all
eigenvalues except the first two are pairwise distinct.

For $m = 1$ the construction is straightforward. $\chi: \Slice_1
\rightarrow \hh/W = \C$ has a single nondegenerate critical point in
the fibre over $\gamma(1) = 0$, hence in the nearby fibres
$\gamma(1-s)$ for small $s$ we have an associated vanishing cycle,
which is a Lagrangian two-sphere. We then use reverse parallel
transport along $\gamma|[0;1-s]$ to move this back to the fibre
$\Y_{1,t}$, which gives us \eqref{eq:matching-lagrangian}. In the
general case, one proceeds by induction on $m$. Let $\bar{\wp}$ be
the crossingless matching obtained from $\wp$ by removing the
component $\delta_1$, and $\bar{t} \in Conf_{2m-2}(\C)$ its
endpoints. By assumption, there is a well-defined Lagrangian
submanifold $L_{\bar{\wp}} \in \Y_{m-1,\bar{t}}$. Lemma
\ref{th:induction} says that one can identify $\Y_{m-1,\bar{t}}$
with the fibre of \eqref{eq:critical-fibration} over the point
$(\mu_3 + \mu_1/(m-1),\dots,\mu_{2m} + \mu_1/(m-1)) \in
Conf_{2m-2}(\C) \iso \hh^{\mult,\reg}/W^{\mult}$. Using the parallel
transport maps \eqref{eq:c-transport} over a path which translates
this configuration, one can move the Lagrangian submanifold to the
fibre of \eqref{eq:critical-fibration} over
$(\mu_3',\dots,\mu_{2m}')$, which by our discussion in Section
\ref{sec:slice}\ref{subsec:partial} is the singular locus of
$\Y_{m,\gamma(1)}$. The local model from Lemma \ref{th:2-coincide}
shows that one can apply the relative vanishing cycle construction,
which yields a Lagrangian submanifold in the fibre
$\Y_{m,\gamma(1-s)}$ for small $s$. As before, reverse parallel
transport is then used to move this to the original fibre, which
gives rise to \eqref{eq:matching-lagrangian}. Topologically, the
relative vanishing cycle procedure takes the product of a given
Lagrangian submanifold with $S^2$, hence the outcome is
diffeomorphic to $(S^2)^m$ as claimed. While the construction
involves many choices, none of them carries any nontrivial topology
(the space of possible choices in each step is path-connected, and
indeed weakly contractible), so that the outcome is well-defined up
to Lagrangian isotopy.

There is another property which follows directly from the
definition. Namely, suppose that we have an open braid $\beta: [0;1]
\rightarrow Conf_{2m}(\C)$ and a smooth family $\wp(s)$ of
crossingless matchings with endpoints $\beta(s)$. Then one can
construct the $L_{\wp(s)} \subset \Y_{m,\beta(s)}$ in such a way
that they depend smoothly on $s$. By using parallel transport over
$\beta|[s;1]$ to carry them into a common fibre, one finds that
there is a Lagrangian isotopy
\begin{equation} \label{eq:transport-matchings}
 L_{\wp(1)} \htp h^{\resc}_{\beta}(L_{\wp(0)}).
\end{equation}
As a particular obvious special case, a smooth family of
crossingless matchings with the same endpoints leads to a family of
isotopic Lagrangian submanifolds.

\begin{figure}[ht]
\begin{centering}
\begin{picture}(0,0)%
\includegraphics{slide.pstex}%
\end{picture}%
\setlength{\unitlength}{3947sp}%
\begingroup\makeatletter\ifx\SetFigFont\undefined%
\gdef\SetFigFont#1#2#3#4#5{%
  \reset@font\fontsize{#1}{#2pt}%
  \fontfamily{#3}\fontseries{#4}\fontshape{#5}%
  \selectfont}%
\fi\endgroup%
\begin{picture}(4736,2647)(1110,-2034)
\put(1989,-149){\makebox(0,0)[lb]{\smash{{\SetFigFont{12}{14.4}{\rmdefault}{\mddefault}{\updefault}$\delta_2$}}}}
\put(5689,-936){\makebox(0,0)[lb]{\smash{{\SetFigFont{12}{14.4}{\rmdefault}{\mddefault}{\updefault}$\delta_2$}}}}
\put(2189,-1361){\makebox(0,0)[lb]{\smash{{\SetFigFont{12}{14.4}{\rmdefault}{\mddefault}{\updefault}$\delta_1$}}}}
\put(1389,-1974){\makebox(0,0)[lb]{\smash{{\SetFigFont{12}{14.4}{\rmdefault}{\mddefault}{\updefault}matching
$\wp$}}}}
\put(4764,-1099){\makebox(0,0)[lb]{\smash{{\SetFigFont{12}{14.4}{\rmdefault}{\mddefault}{\updefault}$\delta_1$}}}}
\put(4451,-1974){\makebox(0,0)[lb]{\smash{{\SetFigFont{12}{14.4}{\rmdefault}{\mddefault}{\updefault}matching
$\wp'$}}}}
\end{picture}%
\caption{\label{fig:slide}}
\end{centering}
\end{figure}

The reverse of the previous statement is false: non-isotopic
crossingless matchings can also sometimes lead to isotopic
Lagrangian submanifolds. Let $\wp,\wp'$ be crossingless matchings
with the same endpoints, which are related to each other as in
Figure \ref{fig:slide}: we choose an embedded path joining the first
and second arc of $\wp$ and avoiding all other components (shown
dashed in the picture), and then define $\wp'$ by passing the second
component over the first as indicated by that path.

\begin{Lemma} \label{th:slide}
$L_{\wp}, L_{\wp'}$ are Lagrangian isotopic.
\end{Lemma}

\begin{Proof}
In view of \eqref{eq:transport-matchings} we may assume that the
endpoints $t = (\mu_1,\dots,\mu_{2m})$ of $\wp$ satisfy
\eqref{eq:normal} as well as $\mu_1 = -\mu_2$, and that the point $0
\in \C$ lies on $\delta_1$. We can then choose $\gamma$ so that the
point $\gamma(1)$ corresponds to $(0,0,\mu_3,\dots,\mu_{2m})$, and
the construction of $L_{\wp}$ simplifies slightly in that the step
involving \eqref{eq:c-transport} becomes trivial. Namely, one
considers $L_{\bar{\wp}} \in \Y_{m-1,\bar{t}}$ where $\bar{t} =
(\mu_3,\dots,\mu_{2m})$, identifies the latter space with the
singular set of $\Y_{m,\gamma(1)}$, takes the associated relative
vanishing cycle in $\Y_{m,\gamma(1-s)}$, and then carries it back to
$\Y_{m,t}$ by parallel transport. The definition of $L_{\wp'}$ is
the same except that we start with $L_{\bar{\wp}'}$. But $\bar{\wp}$
and $\bar{\wp}'$ are isotopic as crossingless matchings with fixed
endpoints, hence $L_{\bar{\wp}} \htp L_{\bar{\wp}'}$, and that
carries over to the associated relative vanishing cycles.
\end{Proof}

\begin{Lemma}
Up to Lagrangian isotopy, $L_{\wp}$ is independent of the ordering
of the components of $\wp$.
\end{Lemma}

\begin{Proof}
Because of the recursive nature of the definition, we only need to
show that exchanging $\delta_1$ and $\delta_2$ does not affect the
Lagrangian submanifold. In view of \eqref{eq:transport-matchings} we
may suppose that $\delta_1$ is a straight short line segment
$[-\sqrt{e_1};\sqrt{e_1}] \subset \C$, and similarly
$\delta_2=[\lambda-\sqrt{e_2};\lambda+\sqrt{e_2}]$, for small
$e_1,e_2 \neq 0$ (in fact, we will only see in the course of the
argument what the precise bounds are, but that is not a problem).

By taking the paths short, we remove the need for using parallel
transport \eqref{eq:y-transport} in the definition of $L_\wp$, at
least for the last two steps in the recursive procedure. What
remains of these steps is the following: one starts with an already
defined Lagrangian submanifold inside the singular point set of
$\Y_{m-1,\bar{t}_1}$, where $\bar{t}_1 =
(\lambda,\lambda,\mu_5,\dots,\mu_{2m})$. The relative vanishing
cycle procedure associates to this a Lagrangian submanifold inside a
nearby smooth fibre $\Y_{m-1,\bar{t}}$, such as $\bar{t} =
(\mu_3,\dots,\mu_{2m})$ if $e_2$ has been chosen sufficiently small.
$\Y_{m-1,\bar{t}}$ can in turn be identified with the singular point
set of $\Y_{m,t_1}$ for $t_1 = (0,0,\mu_3,\dots,\mu_{2m})$, and
forming the relative vanishing cycle again gives a Lagrangian
submanifold in the nearby fibre $\Y_{m,t}$, which is $L_{\wp}$.

We will now reformulate this as follows. Let $w: P \hookrightarrow
\hh/W$ be a small embedded bidisc, so that $w(z_1,z_2)$ corresponds
to the set of eigenvalues $(-\sqrt{z_1},\sqrt{z_1},
\lambda-\sqrt{z_2},\lambda+\sqrt{z_2}, \mu_5,\dots,\mu_{2m})$. Using
Lemma \ref{th:induction} one can identify the singular set of
$\Y_{m-1,\bar{t}_1}$ with the intersection $\Y_{m,w(0,0)} \cap
\OO^{\min}$; here $\OO^{\min} \subset \gg$ is the orbit consisting
of matrices where both the kernel and the $\lambda$-eigenspace are
two-dimensional. Our construction starts with a Lagrangian
submanifold inside this intersection, forms the relative vanishing
cycle inside the critical set of $\Y_{m,w(0,e_2)}$, and then takes
the relative vanishing cycle of that inside the whole smooth fibre
$\Y_{m,w(e_1,e_2)}$. The local structure of $\Slice_m \cap
\chi^{-1}(w(P))$ near $\OO^{\min}$ is described by Lemma
\ref{th:2k-coincide}, and our iteration of the vanishing cycle
procedure is precisely that discussed at the end of Section
\ref{sec:symplectic}\ref{subsec:relative-vanishing}. Lemma
\ref{th:order} allows one to reverse the order in such a procedure,
which in our case corresponds exactly to exchanging $\delta_1$ and
$\delta_2$.
\end{Proof}

We may therefore drop the ordering of the components in the
definition of a crossingless matching.

\begin{Remark}
As the reader may have noticed, the construction of Lagrangian
submanifolds also goes through if one starts with a matching with
ordered components which may intersect each other (of course, the
endpoints must still be distinct, and in addition disjoint
from the interiors of the arcs). However, the result is not really
more general: by the same argument as in the proof of Lemma
\ref{th:slide} one can slide the intersections of the components off
each other to make the matching into a crossingless one, and the
Lagrangian submanifold will not be affected.
\end{Remark}

\subsection{Definition of the invariant}

\begin{figure}[ht]
\begin{centering}
\begin{picture}(0,0)%
\includegraphics{plusminus.pstex}%
\end{picture}%
\setlength{\unitlength}{3947sp}%
\begingroup\makeatletter\ifx\SetFigFont\undefined%
\gdef\SetFigFont#1#2#3#4#5{%
  \reset@font\fontsize{#1}{#2pt}%
  \fontfamily{#3}\fontseries{#4}\fontshape{#5}%
  \selectfont}%
\fi\endgroup%
\begin{picture}(5156,2616)(630,-1903)
\put(4700,-286){\makebox(0,0)[lb]{\smash{{\SetFigFont{12}{14.4}{\rmdefault}{\mddefault}{\updefault}$\wp_-$}}}}
\put(1580,-986){\makebox(0,0)[lb]{\smash{{\SetFigFont{12}{14.4}{\rmdefault}{\mddefault}{\updefault}$\wp_+$}}}}
\end{picture}%
\caption{\label{fig:plusminus}}
\end{centering}
\end{figure}
Fix a $t_0 \in \Conf_{2m}(\C)$ which is a configuration of points on
the real line. We denote by $\wp_+,\wp_-$ the crossingless matchings
with endpoints $t_0$ which consist of a family of concentric arcs in
the upper, respectively lower, half plane (Figure
\ref{fig:plusminus}). These are unique up to isotopy, hence so are
the associated submanifolds \eqref{eq:matching-lagrangian}. A
repeated application of Lemma \ref{th:slide} shows that $L_{\wp_+}$,
$L_{\wp_-}$ are actually isotopic, so we will usually just write
$L_{\wp_{\pm}}$ instead. Next, take an oriented link $\kappa$
presented as a braid closure with $2m$ strands (Figure
\ref{fig:braid}), with the left side of the diagram being $b \in
Br_m$. Via the standard inclusion $Br_m \times Br_m \hookrightarrow
Br_{2m}$ we turn this into a $2m$-stranded braid $b \times 1^m$, and
represent that by a loop $\beta: [0;1] \rightarrow Conf_{2m}(\C)$
starting and ending at $t_0$. As already explained in the
Introduction, we set
\begin{figure}[ht]
\begin{centering}
\begin{picture}(0,0)%
\includegraphics{braid.pstex}%
\end{picture}%
\setlength{\unitlength}{3947sp}%
\begingroup\makeatletter\ifx\SetFigFont\undefined%
\gdef\SetFigFont#1#2#3#4#5{%
  \reset@font\fontsize{#1}{#2pt}%
  \fontfamily{#3}\fontseries{#4}\fontshape{#5}%
  \selectfont}%
\fi\endgroup%
\begin{picture}(2322,3066)(314,-2713)
\put(250,-1211){\makebox(0,0)[lb]{\smash{{\SetFigFont{12}{14.4}{\rmdefault}{\mddefault}{\updefault}$b$}}}}
\end{picture}%
\caption{\label{fig:braid}}
\end{centering}
\end{figure}

\begin{Definition}
\begin{equation} \label{eq:def}
Kh_{\symp}(\kappa) =
HF(L_{\wp_{\pm}},h^{\resc}_\beta(L_{\wp_{\pm}})).
\end{equation}
\end{Definition}

The Floer cohomology is taken inside $\Y_{m,t_0}$ which is certainly
Stein, with the exact K{\"a}hler form $\Omega = -dd^c \psi$. Since
$\Y_{m,t_0}$ is a regular fibre of a holomorphic map $\Slice_m
\rightarrow \hh/W$ between affine spaces, its Chern classes are
zero, and moreover $H^1(\Y_{m,t_0}) = 0$ by Lemma
\ref{th:1st-betti}. Finally, $L_{\wp_{\pm}}$ is diffeomorphic to
$(S^2)^m$, hence has $H_1 = 0$ and is spin, so the Floer cohomology
group above really is well-defined. Explicitly, the choices made in
constructing $L_{\wp_{\pm}}$ -- including the choice of representative
$\beta$ for $b\times 1^m$ and of the choices entering into the
definition of $h^{\resc}_\beta$ -- affect $L_{\wp_{\pm}}$ up to
Lagrangian isotopy, and 
any such isotopy is exact. 

\begin{Lemma} \label{th:alpha}
The Floer cohomology group \eqref{eq:def} is independent of the
choices made in the definition of the K{\"a}hler form $\Omega$.
\end{Lemma}

The statement is independence of the K{\"a}hler form only within the
very restricted class that we are considering, which means
independence of the choice of $\alpha$ and of compactly supported
functions $\eta_k$. But one can linearly interpolate between any two
such forms, and get corresponding smooth families of Lagrangian
submanifolds to use in \eqref{eq:def}, so the result follows from
the invariance properties of Floer cohomology described in Section
\ref{sec:symplectic}\ref{subsec:background}  (proved  by
bifurcation analysis).

The fact that \eqref{eq:def} is an oriented link invariant will
follow from its invariance under Markov moves. The type $I$ move
replaces the braid $b$ by $s_k^{-1} b s_k$ where $s_1,\dots,s_{m-1}$
are the standard generators of $Br_m$ (if one sees braids as
diffeomorphisms of the punctured disc, the $s_k$ are positive
half-twists; it is an unfortunate consequence of the standard
convention that these are represented by braids with a negative
crossing). We will use the same notation for the generators of
$Br_{2m}$, and choose representatives $\sigma_k$ for them which are
loops in $(Conf_{2m}(\C),t_0)$.

\begin{Lemma} \label{th:twist-untwist}
Up to Lagrangian isotopy, $L_{\wp_{\pm}}$ is invariant under
parallel transport along $\sigma_{2m-k}^{-1} \circ \sigma_k$.
\end{Lemma}

\begin{Proof}
We use \eqref{eq:transport-matchings}. Moving the crossingless
matching $\wp_+$ smoothly, so that its endpoints follow
$\sigma_{2m-k}^{-1} \circ \sigma_k$, yields another crossingless
matching shown on the right in Figure \ref{fig:movematching}. That
is clearly obtained from $\wp_+$ by an operation as in Figure
\ref{fig:slide}, so the associated Lagrangian submanifold is
isotopic to $L_{\wp_{\pm}}$ because of Lemma \ref{th:slide}.
\end{Proof}

\begin{figure}[ht]
\begin{centering}
\begin{picture}(0,0)%
\includegraphics{movematching.pstex}%
\end{picture}%
\setlength{\unitlength}{3947sp}%
\begingroup\makeatletter\ifx\SetFigFont\undefined%
\gdef\SetFigFont#1#2#3#4#5{%
  \reset@font\fontsize{#1}{#2pt}%
  \fontfamily{#3}\fontseries{#4}\fontshape{#5}%
  \selectfont}%
\fi\endgroup%
\begin{picture}(5168,1337)(630,-959)
\end{picture}%
\caption{\label{fig:movematching}}
\end{centering}
\end{figure}

\begin{Proposition} \label{th:markov-1}
Up to isomorphism of relatively graded abelian groups, the Floer
cohomology \eqref{eq:def} is invariant under type $I$ Markov moves.
\end{Proposition}

\begin{Proof}
By symplectomorphism invariance of Floer cohomology, and the
previous Lemma, we have
\begin{align*}
 & HF(L_{\wp_{\pm}},h^{\resc}_{\sigma_k^{-1}} h^{\resc}_\beta
 h^{\resc}_{\sigma_k}(L_{\wp_{\pm}})) \\
 & \iso HF(h^{\resc}_{\sigma_k}(L_{\wp_{\pm}}),
 h^{\resc}_\beta h^{\resc}_{\sigma_k}(L_{\wp_{\pm}})) \\
 & \iso HF(h^{\resc}_{\sigma_{2m-k}}(L_{\wp_{\pm}}),
 h^{\resc}_\beta h^{\resc}_{\sigma_k}(L_{\wp_{\pm}})) \\
 & \iso HF(L_{\wp_{\pm}},
 h^{\resc}_{\sigma_{2m-k}^{-1}} h^{\resc}_\beta
 h^{\resc}_{\sigma_k}(L_{\wp_{\pm}})) \\
\intertext{and since $s_{2m-k}$ and $b \times 1^m$ commute in
$Br_{2m}$, this is}
 & \iso HF(L_{\wp_{\pm}},h^{\resc}_\beta
 h^{\resc}_{\sigma_{2m-k}^{-1} \circ \,\sigma_k}(L_{\wp_{\pm}})) \\
 & \iso HF(L_{\wp_{\pm}},h^{\resc}_\beta(L_{\wp_{\pm}})).
\end{align*}
\end{Proof}

\subsection{Markov II}

Before going on to the remaining Markov move we want to deal with a
different property of our Floer groups, whose proof is simpler but
somewhat analogous, hence can serve as a warmup exercise. Namely,
suppose that our oriented link has an unknotted and unlinked
component, which appears in the braid presentation as shown in
Figure \ref{fig:unknot}. This means that $b = \bar{b} \times 1$,
where $\bar{b} \in Br_{m-1}$ and we embed that into $Br_m$ by
considering the leftmost $m-1$ strands.

\begin{figure}[ht]
\begin{centering}
\begin{picture}(0,0)%
\includegraphics{unknot.pstex}%
\end{picture}%
\setlength{\unitlength}{3947sp}%
\begingroup\makeatletter\ifx\SetFigFont\undefined%
\gdef\SetFigFont#1#2#3#4#5{%
  \reset@font\fontsize{#1}{#2pt}%
  \fontfamily{#3}\fontseries{#4}\fontshape{#5}%
  \selectfont}%
\fi\endgroup%
\begin{picture}(2314,2926)(226,-2560)
\put(226,-1186){\makebox(0,0)[lb]{\smash{{\SetFigFont{12}{14.4}{\rmdefault}{\mddefault}{\updefault}$\bar{b}$}}}}
\end{picture}%
\caption{\label{fig:unknot}}
\end{centering}
\end{figure}%

For simplicity, assume that $t_0 = (\mu_1,\dots,\mu_{2m}) \in
\Conf_{2m}^0(\C)$, and that the middle two points of the
configuration are $\mu_m,\mu_{m+1} = \pm \sqrt{e}$ for some small
$e>0$. Define $t_1 \in \hh/W$ by replacing these two with $(0,0)$.
Similarly we define $\bar{t}_0 \in \Conf_{2m-2}^0(\C)$ by deleting
the same points from the configuration, and get a crossingless
matching $\bar{\wp}_{\pm}$ with endpoints $\bar{t}_0$ by removing
the corresponding component from $\wp_{\pm}$. We choose a
representative $\beta$ of $b$ in which the points $\pm \sqrt{e}$
remain fixed, so that by deleting these points one gets a
representative $\bar\beta$ of $\bar{b}$.

Let's start by restating part of the definition of
$h_\beta^{\resc}(L_{\wp_{\pm}})$. One starts with $L_{\bar{\wp}_\pm}
\in \Y_{m-1,\bar{t}_0}$, identifies the latter space with the set of
singular points of $\Y_{m,t_1}$ using Lemma \ref{th:induction}, and
then uses the relative vanishing cycle construction to obtain a
Lagrangian submanifold in in the nearby fibre $\Y_{m,t_0}$, which by
definition is just $L_{\wp_{\pm}}$. One then applies the monodromy
along $\beta$ to that submanifold. Our first claim is that the order
of the two last steps can be inverted, in the sense that if one
applies monodromy along $\bar\beta$ to $L_{\bar{\wp}_{\pm}}$ and
then takes the associated vanishing cycle of the result, an isotopic
Lagrangian submanifold is obtained. The reason is that one can
interpolate continuously between the two processes, by starting with
$L_{\bar{\wp}_{\pm}}$ and applying the monodromy along
$\bar\beta|[0;s]$ for some $s$, then taking the relative vanishing
cycle of the result, and applying the monodromy along $\beta|[s;1]$
to that. The outcome of our discussion is that the pair of
Lagrangian submanifolds
$(L_{\wp_{\pm}},h_{\beta}^{\resc}(L_{\wp_{\pm}}))$ in $\Y_{m,t_0}$
is obtained by taking the pair
$(L_{\bar{\wp}_{\pm}},h_{\bar{\beta}}^{\resc}(L_{\bar{\wp}_{\pm}}))$
inside $\Y_{m-1,\bar{t}_0}$ and applying the relative vanishing
cycle construction to both.

We want to apply Lemma \ref{th:kunneth} to this situation. The total
space $Y$ will be $\chi^{-1}(D) \cap \Slice_m$, where $D \subset
\hh/W$ is a small disc corresponding to eigenvalues
$(\mu_1,\dots,\mu_{m-1},-\sqrt{z},\sqrt{z},\mu_{m+2},\dots,\mu_{2m})$.
$X \subset Y$ is the subset of matrices which have a two-dimensional
kernel. The local structure around $X$ is described by Lemma
\ref{th:2-coincide}, and the other assumptions of Lemma
\ref{th:kunneth} are satisfied for obvious reasons, so
\begin{equation} \label{eq:add-unknot}
 HF(L_{\wp_{\pm}}, h_{\beta}^{\resc}(L_{\wp_{\pm}})) \iso
 HF(L_{\bar{\wp}_{\pm}}, h_{\bar{\beta}}^{\resc}(L_{\bar{\wp}_{\pm}}))
 \otimes H^*(S^2).
\end{equation}
By definition, the first factor on the right is the Floer group
associated to the braid presentation obtained from our original one
by removing the unknotted component.

\begin{figure}[t]
\begin{centering}
\begin{picture}(0,0)%
\includegraphics{markov-2.pstex}%
\end{picture}%
\setlength{\unitlength}{3947sp}%
\begingroup\makeatletter\ifx\SetFigFont\undefined%
\gdef\SetFigFont#1#2#3#4#5{%
  \reset@font\fontsize{#1}{#2pt}%
  \fontfamily{#3}\fontseries{#4}\fontshape{#5}%
  \selectfont}%
\fi\endgroup%
\begin{picture}(3211,4106)(276,-3748)
\put(676,-1449){\makebox(0,0)[lb]{\smash{{\SetFigFont{12}{14.4}{\rmdefault}{\mddefault}{\updefault}$\bar{b}$}}}}
\put(2209,-2536){\makebox(0,0)[lb]{\smash{{\SetFigFont{12}{14.4}{\rmdefault}{\mddefault}{\updefault}$s_{m-1}$}}}}
\end{picture}
\caption{\label{fig:markov-2}}
\end{centering}
\end{figure}

The basic setup for the type $II^+$ Markov move is that one starts
with $\bar{b} \in Br_{m-1}$, and then adds a single strand plus a
half-twist of that strand with its neighbour, $b = s_{m-1} (\bar{b}
\times 1) \in Br_m$. Suppose that $m \geq 3$, and assume that the
base point $t_0 = (\mu_1,\dots,\mu_{2m})$ lies in
$\Conf^0_{2m}(\C)$, with the $\mu_k$ ordered in the obvious way on
the real line, and with $\mu_{m-1},\mu_m,\mu_{m+1}$ small and
satisfying $\mu_{m-1} + \mu_m + \mu_{m+1} = 0$ (we will constrain
the choice of these during the course of the argument, which is not
a problem). Let $\bar{t}_0 \in \Conf_{2m-2}^0(\C)$ be the
configuration $(\mu_1,\dots,\mu_{m-2},0,\mu_{m+2},\dots,\mu_{2m})$.
We choose loops $\beta$ and $\beta'$ in $\Conf_{2m}^0(\C)$ based at
$t_0$, representing $b \times 1^m$ and $\bar{b} \times 1^{m+1} =
s_{m-1}^{-1} (b \times 1^m)$ respectively, and similarly a loop
$\bar\beta$ in $\Conf_{2m-2}^0(\C)$ based at $\bar{t}_0$ which
represents $\bar{b} \times 1^{m-1}$ (Figure \ref{fig:markov-2}).

As in Lemma \ref{th:3-coincide}, we consider the embedding $w: P
\rightarrow \hh/W$ of a small bidisc, so that $w(d,z)$ corresponds
to eigenvalues $(\mu_1,\dots,\mu_{m-2},\text{all solutions of }
\lambda^3 - d\lambda + z = 0, \mu_{m+2},\dots,\mu_{2m})$. Start with
the pair
\begin{equation} \label{eq:k-k}
 K = L_{\bar{\wp}_{\pm}}, \quad
 K' = h^{\resc}_{\bar{\beta}}(L_{\bar{\wp}_{\pm}})
\end{equation}
of Lagrangian submanifolds in $\Y_{m-1,\bar{t}_0}$, and identify
$\Y_{m-1,\bar{t}_0}$ itself with the singular point set of
$(\chi|\Slice_m)^{-1}w(0,0)$. The key to the argument is Lemma
\ref{th:3-coincide} which describes the local structure of
$(\chi|\Slice_m)^{-1}w(P)$ near that subset, since that allows us to
apply the results of the discussion from Section
\ref{sec:symplectic}\ref{subsec:a2}; we will use terms from that
discussion freely from now on. By using parallel transport in the
family of singular sets over the cusp curve in $P$ (which
corresponds to sets of eigenvalues where two coincide), one moves
both of \eqref{eq:k-k} to Lagrangian submanifolds $K_d,K_d'$ in the
singular set of a nearby fibre $(\chi|\Slice_m)^{-1}w(d,\zeta_d^-)$.
Here $d>0$ is small and $\zeta_d^-$ is the negative solution of
$4d^3 + 27z^2 = 0$. Then, by the relative vanishing cycle procedure
applied to the family of fibres $(\chi|\Slice_m)^{-1}(d,z)$ with
fixed $d$, one gets associated Lagrangian submanifolds
$L_{d,\epsilon}, L'_{d,\epsilon}$ in
$(\chi|\Slice_m)^{-1}w(d,\zeta_d^- + \epsilon)$ for $0 < \epsilon
\ll d$. It is no problem to assume that our base point $t_0$ was in
fact chosen so that $\mu_{d-1},\mu_d,\mu_{d+1}$ are the solutions of
$\lambda^3 - d\lambda + (\zeta^d + \epsilon) = 0$, and then the
Lagrangian submanifolds which we constructed lie precisely in
$\Y_{m,t_0}$. The same argument as in the proof of
\eqref{eq:add-unknot}, inverting the order of parallel transport and
relative vanishing cycle procedures, allows us to identify the
outcome up to Lagrangian isotopy:
\begin{equation} \label{eq:first-lagrangians}
 L_{d,\epsilon} \htp L_{\wp_{\pm}}, \quad
 L_{d,\epsilon}' \htp h^{\resc}_{\beta'}(L_{\wp_{\pm}}).
\end{equation}
The second Lagrangian submanifold in \eqref{eq:first-lagrangians} is
not yet quite the one which would appear in the formula
\eqref{eq:def} for the link diagram from Figure \ref{fig:markov-2}.
What we are missing is the generator $s_{k-1}$, which corresponds to
moving $\mu_{m-1}, \mu_m$ around each other in a positive
half-circle. Reversing an argument made in Section
\ref{sec:symplectic}\ref{subsec:a2}, we find that in terms of
coordinates on $P$ this can be achieved by fixing $d$ and moving $z$
along the loop $\gamma_{d,\epsilon}$ from Figure
\ref{fig:gamma-path}. We therefore obtain the following modified
version of \eqref{eq:first-lagrangians}: if $L_{d,\epsilon}''$ is
defined as in \eqref{eq:d-epsilon-two} then
\begin{equation}
 L_{d,\epsilon} \htp L_{\wp_{\pm}}, \quad
 L_{d,\epsilon}'' \htp h^{\resc}_\beta(L_{\wp_{\pm}}).
\end{equation}
Applying Lemma \ref{th:thom}, whose assumptions are easily verified
due to Lemma \ref{th:3-coincide}, one finds that
$HF(L_{\wp_{\pm}},h^{\resc}_\beta(L_{\wp_{\pm}}))$ is isomorphic to
$HF(K,K')$, which by definition is the Floer group associated to the
link presentation with $2m-2$ strands and braid $\bar{b}$. Up to
now, we have excluded the lowest strand case $m = 2$, since then the
conditions imposed above on $t$ are impossible to satisfy without
violating \eqref{eq:normal}. Concretely, this means that one has to
bring the three eigenvalues $\mu_{m-1},\mu_m,\mu_{m+1}$ together at
some nonzero point. That adds a small and entirely harmless
intermediate step to the proof, which is the use of parallel
transport for a suitable family $\CC_{1,\gamma(s)}$ to bring that
point back to zero, and there to make the identification with
$\Y_{1,\bar{t}_0}$. With that taken into account, we have shown:

\begin{Proposition} \label{th:markov-2}
Up to isomorphism of relatively graded abelian groups, the Floer
cohomology \eqref{eq:def} is invariant under type $II^+$ Markov
moves. \qed
\end{Proposition}

The discussion above can be easily adapted to the Markov move of
type $II^-$, where $s_{m-1}^{-1}$ occurs instead of $s_{m-1}$. The
geometry of this situation is very similar to the previous one, the
difference being that the path $\gamma_{d,\epsilon}$ has to be taken
with reversed orientation, and consequently that $t_\beta(\alpha)$
in Figure \ref{fig:matching} should be replaced by
$t_\beta^{-1}(\alpha)$. This still intersects $\alpha$ only in the
rightmost endpoint, so the proof of Lemma \ref{th:thom} goes through
exactly as before. We omit the details. By combining this with
Proposition \ref{th:markov-1} and \ref{th:markov-2}, one gets that
$Kh_{\symp}(\kappa)$ (as a relatively graded group) is an invariant
of the oriented link $\kappa$.

\section{Miscellany\label{sec:misc}}

This section lists some concluding observations. First we explain
how to equip $\Kh_{\symp}$ with a suitable absolute grading, which
rounds off the construction of the invariant in the form presented
in the Introduction. Secondly, we see why orientation reversal of
all components leaves it unchanged. Finally, we compute it for the
trefoil knot. Since these are somewhat peripheral topics (even
though they have some relevance in view of Conjecture
\ref{th:conj}), we will give less details than usual.

\subsection{Gradings\label{subsec:gradings}}

Let $M$ be a Stein manifold with an exact K{\"a}hler form,
satisfying $c_1(M) = 0$ and $H^1(M) = 0$. Pick a differentiable
trivialization of the canonical bundle, which is a nowhere zero
complex volume form $\eta_M$. Any Lagrangian submanifold $L \subset
M$ comes with a canonical circle-valued squared phase function
$\alpha_L : L \rightarrow S^1$, defined by
\begin{equation} \label{eq:phase}
 \alpha_L(x) = \frac{\eta_M(\xi_1,\dots,\xi_n)^2}
 {|\eta_M(\xi_1,\dots,\xi_n)|^2}
\end{equation}
for any orthonormal basis $\xi_1,\dots,\xi_n$ of $TL_x$. By
definition a grading of $L$ is a lift of this to a real-valued
function $\tilde{\alpha}_L$, let's say for concreteness
$\textit{exp}(2\pi i \tilde{\alpha}_L) = \alpha_L$. This is always
possible in the context we used for defining Floer cohomology, since
$H^1(L) = 0$ was part of the assumptions. For a pair of Lagrangian
submanifolds $L_0,L_1$ equipped with gradings, the relative grading
of Floer cohomology can be improved to an absolute one
\cite{kontsevich94,seidel99}. If we denote by $L \mapsto L[1]$ the
process which subtracts the constant $1$ from the grading, then
\begin{equation} \label{eq:both-sided-shift}
 HF^*(L_0,L_1[1]) = HF^*(L_0[-1],L_1) = HF^{*+1}(L_0,L_1).
\end{equation}
It may seem that this theory depends on the choice of $\eta_M$, but
in fact all that matters is its homotopy class as smooth
trivialization, which is unique since $H^1(M) = 0$.

In our application, we start by choosing arbitrary trivializations
$\eta_{\Slice_m}$ and $\eta_{\hh/W}$ on the total space and base space of
$\chi|\Slice_m$. There is an induced family of trivializations of
the canonical bundles of the regular fibres, characterized by
\begin{equation} \label{eq:division}
 \eta_{\Y_{m,t}} \wedge \chi^*\eta_{\hh/W} = \eta_{\Slice_m}
 \qquad \text{on $\Y_{m,t}$.}
\end{equation}
Choose a grading for $L_{\wp_{\pm}} \subset \Y_{m,t_0}$. Given a
path $\beta: [0;1] \rightarrow \Conf_{2m}(\C)$ starting at $t_0$,
one can continue the given grading uniquely to a smooth family of
gradings of the images $h^{\resc}_{\beta|[0;s]}(L_{\wp_{\pm}})$, in
particular the monodromy images which appear in \eqref{eq:def} carry
induced gradings, so the Floer cohomology group in that definition
is now absolutely graded. Shifting the original choice of grading
affects both Lagrangian submanifolds involved in the same way, and
the effect on Floer cohomology cancels out due to
\eqref{eq:both-sided-shift}. We take this absolutely graded group
and apply a final shift to the grading, which depends on the number
of strands $m$ and writhe $w$ of the braid presentation, thus
arriving at the final definition:
\begin{equation} \label{eq:graded-kh}
 \Kh^*(\kappa) =
 HF^{*+m+w}(L_{\wp_{\pm}},h^{\resc}_\beta(L_{\wp_{\pm}})).
\end{equation}
The isomorphisms in the proof of Proposition \ref{th:markov-1} are
compatible with the absolute gradings, and the writhe does not
change since we add one positive and one negative crossing, so
\eqref{eq:graded-kh} is invariant under Markov $I$.

We next look at the role of absolute gradings in
\eqref{eq:add-unknot}. The basic situation in the proof of that
isomorphism is encoded in the diagram
\begin{equation} \label{eq:grading-diagram}
\begin{CD}
 \chi^{-1}(w(D)) \cap \Slice_m
 @>{\text{local $\iso$ defined near $\Y_{m-1,\bar{t}_0}$}}>>
 & \Y_{m-1,\bar{t}_0} \times \C^3 \\
 @V{w^{-1} \circ \chi}VV & @V{a^2+b^2+c^2}VV \\
 D @>{\qquad\qquad\qquad\qquad\qquad}>> & \C
\end{CD}
\end{equation}
where $w: D \rightarrow \hh/W$, for $D \subset \C$ a small disc
around the origin, maps $z$ to the set of eigenvalues
$(\mu_1,\dots,\mu_{m-1},-\sqrt{z},\sqrt{z},\mu_{m+2},\dots,\mu_{2m})$;
and $\Y_{m-1,\bar{t}_0}$ is identified with the singular set of
$\Y_{m,w(0)}$. We can assume that on the subset of $\chi^{-1}(w(D))
\cap \Slice_m$ which is the domain of the $\iso$ in
\eqref{eq:grading-diagram}, $\eta_{\Slice_m}$ is the wedge product
of a previously defined $\eta_{\Y_{m-1,\bar{t}_0}}$, the standard
form $da \wedge db \wedge dc$ on $\C^3$, and the form $d\mu_1 \wedge
\dots \wedge d\mu_{m-1} \wedge d\mu_{m+2} \wedge \dots \wedge
d\mu_{2m}$. Similarly, we may assume that on the image of $w$,
$\eta_{\hh/W} = dz \wedge d\mu_1 \wedge d\mu_{m-1} \wedge d\mu_{m+2}
\wedge \dots \wedge d\mu_{2m}$ where $z$ is the parameter of $D$.
This is because we have complete freedom in the choice of these
volume forms, so we can prescribe them arbitrarily on any subset
with zero first Betti number (this condition is imposed to ensure
extendibility to the whole space).

As in the discussion preceding \eqref{eq:add-unknot}, we start with
the Lagrangian submanifolds $K = L_{\bar{\wp}_{\pm}}$ and $K' =
h^{\resc}_{\bar{\beta}}(L_{\bar{\wp}_{\pm}})$ inside
$\Y_{m-1,\bar{t}_0}$, which have already been equipped with gradings
following the prescription given above, and then take the associated
relative vanishing cycles $L_z,L'_z$ in $\Y_{m,w(z)}$ for some small
$z \neq 0$. For this we may use the local isomorphism
\eqref{eq:grading-diagram} and a K{\"a}hler form which is the
product of the given one on $\Y_{m-1,\bar{t}_0}$ and the standard
form on $\C^3$, since that is how the Floer cohomology computation
in Lemma \ref{th:kunneth} is carried out anyway. Then $L_z = K
\times \sqrt{z} S^2$, $L_z' = K' \times \sqrt{z} S^2$, and a
straightforward computation of the relevant phases shows that
\[
 \alpha_{L_z} = \alpha_K \cdot \frac{z}{|z|}, \quad
 \alpha_{L_z'} = \alpha_{K'} \cdot \frac{z}{|z|}.
\]
To equip the relative vanishing cycles with gradings, what one has
to do is therefore to choose a branch of $\textit{arg}(z)$. The main
thrust of the proof of \eqref{eq:add-unknot} is that these vanishing
cycles are related to $L_{\wp_{\pm}}$ and
$h^{\resc}_{\beta}(L_{\wp_{\pm}})$, respectively. Inspection of the
argument shows that in order to make this relation work on the level
of Lagrangian submanifolds equipped with gradings, the same branch
of $\textit{arg}(z)$ has to be used for both $L_z$ and $L_z'$. In
that case, the K{\"u}nneth formula from Lemma \ref{th:kunneth} holds
as an isomorphism of absolutely graded groups where $H^*(S^2)$
carries its natural grading. Taking into account the additional
shift that comes from the number of strands, we find that:

\begin{Proposition}
Under disjoint sum with an unlinked unknot, $Kh_{\symp}^*(L \sqcup
U) \iso Kh_{\symp}^*(L) \otimes H^{*+1}(S^2)$. \qed
\end{Proposition}

The role of the grading in Markov $II^+$ is fundamentally very
similar, with an additional contribution to the phase function
coming from the $S^2$ factor added when taking relative vanishing
cycles, and a corresponding correction term to the degree of
intersection points. This correction will be the same in all cases,
so it is enough to look at the toy model example studied at the
beginning of Section \ref{sec:symplectic}\ref{subsec:a2}. Hence, let
$Y_{d,\zeta_d^- + \epsilon}$ be the fibre of the map $p: \C^4
\rightarrow \C^2$ from \eqref{eq:a2-map} for some $0 < \epsilon \ll
d$. We consider the Lagrangian spheres $L_{d,\epsilon}$ and
$h_{\gamma_{d,\epsilon}}(L_{d,\epsilon})$ from \eqref{eq:l0} and
\eqref{eq:monodromy-image}, equipping the first one with an
arbitrary grading and the second one with the induced grading,
coming from the fact that it is a monodromy image of the first. We
take the standard K{\"a}hler form, and then apply a Lagrangian
isotopy if necessary, so that following Lemma \ref{th:bc} our
Lagrangian spheres are $\Lambda_\alpha$ and
$\Lambda_{t_\beta(\alpha)}$ respectively. As shown in the proof of
that Lemma, the monodromy corresponds to the half-twist around
$\beta$ in the base, so the gradings we have chosen will have the
property that they are approximately the same at the unique
intersection point, which we call $q$ (assuming that as represented
in Figure \ref{fig:matching}, the angle between $\alpha$ and
$t_\beta(\alpha)$ at the common endpoint is small). Near $q$ one can
locally write
\[
\Lambda_\alpha = \textit{graph}(dh) \subset
T^*\Lambda_{t_\beta(\alpha)},
\]
where $h$ has a nondegenerate local minimum at $q$. In view of the
fact about the gradings mentioned above, standard properties of the
Maslov index \cite{robb-sal93} imply that the Maslov index of $q$
reduces to the Morse index of $(D^2h)_q$, which is $0$. The same
argument works for Markov $II^-$ except that the second path is now
$t_\beta^{-1}(\alpha)$, which runs to the left of $\alpha$, and so
the function $h$ has a local maximum, leading to a Maslov index of
$2$. These are the desired correction terms, and so the graded
versions of the isomorphism arising from Lemma \ref{th:thom} and its
analogue are as follows:

\begin{Lemma}
Take $\bar{b} \in Br_{m-1}$ and set $b^{\pm} = s_{m-1}^{\pm}(\bar{b}
\times 1) \in Br_m$. Let $\bar\beta$ be a path in $\Conf_{2m-2}(\C)$
representing $\bar{b} \times 1^{m-1}$, and similarly $\beta^{\pm}$
paths in $\Conf_{2m}(\C)$ representing $b^{\pm} \times 1^m$. Then
there are isomorphisms of graded abelian groups,
\begin{align*}
 & HF^*(L_{\wp_{\pm}},h^{\resc}_{\beta^+}(L_{\wp_{\pm}})) \iso
 HF^*(L_{\bar{\wp}_{\pm}},h^{\resc}_{\bar\beta}(L_{\bar{\wp}_{\pm}})),
 \\
 & HF^*(L_{\wp_{\pm}},h^{\resc}_{\beta^-}(L_{\wp_{\pm}})) \iso
 HF^{*-2}(L_{\bar{\wp}_{\pm}},h^{\resc}_{\bar\beta}(L_{\bar{\wp}_{\pm}})).
\end{align*}
\qed
\end{Lemma}

Inspection of \eqref{eq:graded-kh} shows that these precisely cancel
out against the changes in $m+w$. This proves the invariance under
Markov $II$ of $\Kh_{\symp}^*$ as a graded group, and thereby
completes our proof of Theorem \ref{th:main}.

\subsection{Orientation-reversal}

Complex conjugation $c(y) = \bar{y}$ acts on $\Slice_m$, and induces
the obvious map, also denoted by $c$, on the base of the adjoint
quotient map $\hh/W$. In particular, if $t \in \Conf_{2m}^0(\C)$
consists of real eigenvalues, like our base point $t_0$, we have an
induced involution of $\Y_{m,t}$ which reverses the sign of the
K{\"a}hler form. As an obvious consequence of its behaviour on
$\hh/W$, and the definition of parallel transport and of the
Lagrangian submanifolds associated to crossingless matchings, we
have
\begin{equation} \label{eq:involution}
\begin{aligned}
 & c \circ h_{\beta}^{\resc} \circ c = h_{c(\beta)}^{\resc}, \\
 & c(L_{\wp}) = L_{c(\wp)}.
\end{aligned}
\end{equation}
In particular, since $c(\wp_{\pm}) = \wp_{\mp}$, our basic
Lagrangian submanifold $L_{\wp_\pm}$ is invariant under this
involution up to isotopy. Since $c$ is antisymplectic, it induces
isomorphisms on Floer cohomology groups which exchange the two
factors involved,
\begin{equation} \label{eq:reversal}
 HF(c(L_0),c(L_1)) \iso HF(L_1,L_0).
\end{equation}
A choice of grading for $L_k$ induces a grading for $c(L_k)$, and with
that in mind \eqref{eq:reversal} becomes an isomorphism of graded
groups. This is just the fact that complex conjugation acts on the
first cohomology of the Grassmannian of Lagrangian subspaces in
a symplectic vector space by multiplication by $-1$; suitably unwound
(and coupled with the definition of the absolute grading \cite{seidel99}),
this implies 
that $c(L_k[-1]) = (c(L_k))[1]$, in other words increasing the
absolute Maslov index of an intersection point of the $L_k$
decreases the index of the point viewed as an intersection of
the $c(L_k)$.  By combining this with \eqref{eq:involution} and
symplectomorphism invariance of Floer cohomology, one finds that
\begin{align*}
 & HF^*(L_{\wp_{\pm}},h^{\resc}_\beta(L_{\wp_{\pm}})) \\
 & \iso HF^*((c \circ h^{\resc}_\beta)(L_{\wp_{\pm}}),
 c(L_{\wp_{\pm}})) \\
 & \iso
 HF^*(L_{\wp_{\pm}},(h^{\resc}_{c(\beta)})^{-1}(L_{\wp_{\pm}})).
\end{align*}
If $\beta$ represents $b \times 1^m$ for some $b \in Br_m$, then
$c(\beta)^{-1}$ represents $c(b) \times 1^m$, where the braid $c(b)$
is obtained from $b$ by the antiautomorphism of $Br_m$ which inverts
the order of the letters in a word with respect to the standard
presentation. If $b$ gives a braid presentation for an oriented link
$\kappa$, then $c(b)$ corresponds to the same link with the
orientation of all components reversed. Since both presentations
have the same number of strands and the same writhe, the computation
above shows:

\begin{Proposition}
Up to isomorphism of graded groups, $\Kh^*_{\symp}(\kappa)$ remains
unchanged if we reverse the orientation of all components of
$\kappa$. \qed
\end{Proposition}

\subsection{The trefoil}

We now look at the left-handed trefoil knot $\kappa$ (more precisely, the knot
coming from the braid closure with $b = s_1^3 \in Br_2$). The first
part of the proof is to reduce things to an open subset of
$\Y_{2,t_0}$ where one has nice holomorphic coordinates, and then
deform the K{\"a}hler form to a more standard one. This runs
entirely parallel to the corresponding argument for Markov $II$, so
we will omit it and simply state the outcome.

Let $X \subset \C^3$ be the quadric $u^2 + v^2 + w^2 = z$ for some
$z \neq 0$. Choose some $\sqrt{z}$, and consider the line bundle
$\FF \rightarrow X$ whose fibre is the $i\sqrt{z}$-eigenspace of the
matrix
\[
\begin{pmatrix} iu & v+iw \\ -v+iw & -iu \end{pmatrix}.
\]
Inside $\C \oplus \FF^{-2} \oplus \FF^2 \rightarrow X$ with fibre
coordinates $(a,b,c)$, consider the sub-fibre bundle $Y$ defined by
$a^3 - ad + bc = \zeta_d^- + \epsilon$ for some small $0 < \epsilon
\ll d$, and where $\zeta_d^-$ is as in Section
\ref{sec:symplectic}\ref{subsec:a2}. We construct a K{\"a}hler form
on $Y$ (or more precisely on an open subset which is sufficiently
large for our purpose) by combining the standard form on the base
and fibre, and a hermitian metric on $\FF$, as set out in Background
\ref{th:background}. Define a Lagrangian submanifold $L \subset Y$
by taking $\sqrt{z} S^2 \subset X$ on the base, and fibrewise over
it the Lagrangian sphere $\Lambda_\alpha$ from
\eqref{eq:fibered-spheres} in the fibres; and another Lagrangian $L'
\subset Y$ in the same way using $\alpha' = t_\beta^3(\alpha)$
instead, see Figure \ref{fig:matching3} (the basic notation is
carried over from Figure \ref{fig:matching}).

\begin{figure}[h]
\begin{centering}
\begin{picture}(0,0)%
\includegraphics{matching3.pstex}%
\end{picture}%
\setlength{\unitlength}{3947sp}%
\begingroup\makeatletter\ifx\SetFigFont\undefined%
\gdef\SetFigFont#1#2#3#4#5{%
  \reset@font\fontsize{#1}{#2pt}%
  \fontfamily{#3}\fontseries{#4}\fontshape{#5}%
  \selectfont}%
\fi\endgroup%
\begin{picture}(3190,1382)(794,-1133)
\put(2026,50){\makebox(0,0)[lb]{\smash{{\SetFigFont{12}{14.4}{\familydefault}{\mddefault}{\updefault}$\alpha'
 = t_{\beta}^3(\alpha)$}}}}
\put(3376,-854){\makebox(0,0)[lb]{\smash{{\SetFigFont{12}{14.4}{\familydefault}{\mddefault}{\updefault}$\alpha$}}}}
\end{picture}%
\caption{\label{fig:matching3}}
\end{centering}
\end{figure}%
\begin{Lemma}
$L \cap L' \iso S^2 \sqcup \mathbb{RP}^3$.
\end{Lemma}

\begin{Proof}
The paths $\alpha$ and $\alpha'$ intersect in one endpoint and one
interior point, and the corresponding intersections of
$\Lambda_\alpha$ and $\Lambda_\alpha'$ consist of a single point and
a circle \eqref{eq:circle} respectively. This takes place in each
fibre over $(u,v,w) \in \sqrt{z}S^2$, leading to a total
intersection which is the disjoint union of a copy of the $S^2$ and
a circle bundle over it. The degree of the circle bundle equals the
multiplicity of the $S^1$-action \eqref{eq:slice-action} on the
circle \eqref{eq:circle}, which is $\pm 2$.
\end{Proof}

The intersection $L \cap L'$ is clean in the sense of
\cite{pozniak}, so one has a Morse-Bott type long exact sequence
\begin{equation} \label{eq:morse-bott}
 \cdots H^{*-2}(\mathbb{RP}^3) \rightarrow HF^*(L,L') \rightarrow
 H^*(S^2) \stackrel{\partial}{\rightarrow} H^{*-1}(\mathbb{RP}^3) \cdots
\end{equation}
We have given the gradings and resulting Maslov indices in this
sequence without proof, but the nontrivial contribution to them
comes from the geometry in the fibres of $Y \rightarrow X$, and can
be read off from \cite[Lemma 6.18]{khovanov-seidel98}. The
differential $\partial$ is necessarily zero, and taking into account
the shift factor $m + w = 2-3 = -1$, we have

\begin{Proposition}
$\Kh_{\symp}^*(\kappa) \iso H^{*-1}(S^2) \oplus
H^{*-3}(\mathbb{RP}^3)$. \qed
\end{Proposition}

This agrees with the computation of \cite[Section 7]{khovanov98}
after collapsing the bigrading according to the prescription of
Conjecture \ref{th:conj}.


\end{document}